\documentclass{article}
\usepackage{amsmath}
\usepackage{amssymb}
\begin{document}
\def\e#1\e{\begin{equation}#1\end{equation}}
\def\ea#1\ea{\begin{align}#1\end{align}}
\def\eq#1{{\rm(\ref{#1})}}
\newtheorem{thm}{Theorem}[section]
\newtheorem{lem}[thm]{Lemma}
\newtheorem{prop}[thm]{Proposition}
\newtheorem{cor}[thm]{Corollary}
\newenvironment{dfn}{\medskip\refstepcounter{thm}
\noindent{\bf Definition \thesection.\arabic{thm}\ }}{\medskip}
\newenvironment{ex}{\medskip\refstepcounter{thm}
\noindent{\bf Example \thesection.\arabic{thm}\ }}{\medskip}
\newenvironment{proof}[1][,]{\medskip\ifcat,#1
\noindent{\it Proof.\ }\else\noindent{\it Proof of #1.\ }\fi}
{\hfill$\square$\medskip}
\def\dim{\mathop{\rm dim}}
\def\mod{\mathop{\rm mod}}
\def\Re{\mathop{\rm Re}}
\def\Im{\mathop{\rm Im}}
\def\Image{\mathop{\rm Image}}
\def\Aff{\mathop{\rm Aff}}
\def\Vect{\mathop{\rm Vect}}
\def\vol{\mathop{\rm vol}}
\def\hcf{\mathop{\rm hcf}}
\def\GL{\mathop{\rm GL}}
\def\SO{\mathop{\rm SO}}
\def\U{\mathbin{\rm U}}
\def\SL{\mathop{\rm SL}}
\def\SU{\mathop{\rm SU}}
\def\sn{{\textstyle\mathop{\rm sn}}}
\def\cn{{\textstyle\mathop{\rm cn}}}
\def\dn{{\textstyle\mathop{\rm dn}}}
\def\sech{\mathop{\rm sech}}
\def\cosech{\mathop{\rm cosech}}
\def\ge{\geqslant}
\def\le{\leqslant}
\def\cal{\mathcal}
\def\bb{\mathbb}
\def\R{\mathbin{\mathbb R}}
\def\Z{\mathbin{\mathbb Z}}
\def\Q{\mathbin{\mathbb Q}}
\def\C{\mathbin{\mathbb C}}
\def\g{\mathbin{\mathfrak g}}
\def\h{\mathbin{\mathfrak h}}
\def\al{\alpha}
\def\be{\beta}
\def\ga{\gamma}
\def\de{\delta}
\def\ep{\epsilon}
\def\up{\upsilon}
\def\th{\theta}
\def\la{\lambda}
\def\ka{\kappa}
\def\vp{\varphi}
\def\si{\sigma}
\def\La{\Lambda}
\def\Om{\Omega}
\def\Si{\Sigma}
\def\om{\omega}
\def\d{{\rm d}}
\def\pd{\partial}
\def\ts{\textstyle}
\def\sst{\scriptscriptstyle}
\def\pha{\phantom}
\def\w{\wedge}
\def\lt{\ltimes}
\def\sm{\setminus}
\def\ov{\overline}
\def\iy{\infty}
\def\ra{\rightarrow}
\def\t{\times}
\def\ha{{\ts{\frac{1}{2}}}}
\def\op{\oplus}
\def\ot{\otimes}
\def\ms#1{\vert#1\vert^2}
\def\bms#1{\bigl\vert#1\bigr\vert^2}
\def\md#1{\vert #1 \vert}
\def\bmd#1{\bigl\vert #1 \bigr\vert}
\def\an#1{\langle#1\rangle}
\def\ban#1{\bigl\langle#1\bigr\rangle}
\title{Evolution equations for \\
special Lagrangian 3-folds in ${\mathbb C}^3$}
\author{Dominic Joyce \\ Lincoln College, Oxford}
\date{}
\maketitle

\section{Introduction}
\label{ev1}

This is the third in a series of papers \cite{Joyc2,Joyc3,Joyc4,Joyc5,Joyc6}
constructing examples of special Lagrangian submanifolds (SL $m$-folds) in 
$\C^m$. The principal motivation for these papers is to lay the foundations
for a study of the singularities of compact special Lagrangian $m$-folds 
in Calabi--Yau $m$-folds, particularly in low dimensions such as $m=3$. 
Special Lagrangian $m$-folds in $\C^m$ should provide local models for 
how singularities develop in special Lagrangian $m$-folds in Calabi--Yau 
$m$-folds.

Understanding such singularities will be essential in making 
rigorous the explanation of Mirror Symmetry of Calabi--Yau 3-folds 
$X,\hat X$ proposed by Strominger, Yau and Zaslow \cite{SYZ}, which 
involves dual `fibrations' of $X,\hat X$ by special Lagrangian 3-tori, 
with some singular fibres. It will also be important in resolving 
conjectures made by the author \cite{Joyc1}, which attempt to define 
an invariant of Calabi--Yau 3-folds by counting special Lagrangian
homology 3-spheres.

In this paper we construct and study several families of special
Lagrangian 3-folds in $\C^3$, using the `evolution' construction
method of \cite{Joyc3}. In \cite[\S 3--\S 4]{Joyc3} we described
a general construction of SL $m$-folds in $\C^m$, which will be
summarized in \S\ref{ev3} below. The construction requires a set
of {\it evolution data} $(P,\chi)$, including an
$(m\!-\!1)$-submanifold $P$ in $\R^n$ for $n\ge m$. Then $N$ is
the subset of $\C^m$ swept out by the image of $P$ under a
1-parameter family of linear or affine maps $\phi_t:\R^n\ra\C^m$,
which satisfy a first-order, nonlinear o.d.e.\ in~$t$.

In \cite{Joyc3} we restricted our attention to evolution data
$(P,\chi)$ in which $n=m$ and $P$ is a quadric in $\R^m$. In
\S\ref{ev4} we shall establish a correspondence between sets of
evolution data with $m=3$ and homogeneous symplectic 2-manifolds
$\Si$ with a transitive, faithful, Hamiltonian action of a Lie
group $G$. This enables us to write down various interesting new
examples of evolution data with~$m=3$.

We shall focus on two examples in particular, and study the
associated SL 3-folds in $\C^3$ in detail. The first, given in
Example \ref{ev4ex1}, comes from the action of $\SL(2,\R)\lt\R^2$
on $\R^2$ by affine transformations, and the corresponding family
of SL 3-folds are discussed in sections~\ref{ev5}--\ref{ev11}.

The second example, in Example \ref{ev4ex2}, comes from the action of
$K\lt U_k$ on $T^*\R$ for $k\ge 1$, where $K=\R_+\lt\R$ is the group
of oriented affine transformations of $\R$, and $U_k$ the vector space
of polynomial 1-forms on $\R$ of degree less than $k$, acting on $T^*\R$
by translation in the fibres. This gives a series of families of ruled
SL 3-folds in $\C^3$, which are discussed in~\S\ref{ev12}.

The construction of \S\ref{ev5}--\S\ref{ev11} involves a family of
quadratic maps $\Phi_t:\R^2\ra\C^3$ depending on $t\in\R$, which
satisfy the~o.d.e.
\begin{equation*}
\frac{\d\Phi_t}{\d t}=\frac{\pd\Phi_t}{\pd y_1}\t\frac{\pd\Phi_t}{\pd y_2},
\end{equation*}
where $(y_1,y_2)$ are the coordinates on $\R^2$, and `$\t$'
is an anti-bilinear cross product on $\C^3$. Defining
$\Phi:\R^3\ra\C^3$ by $\Phi(y_1,y_2,t)=\Phi_t(y_1,y_2)$, it
turns out that under certain conditions on the initial data
$\Phi_0$ the image $N=\Image\Phi$ is a special Lagrangian 3-fold.

For generic initial data this map $\Phi$ is an immersion, so that
$N$ is a nonsingular immersed 3-submanifold diffeomorphic to $\R^3$.
But for a set of initial data of real codimension one, $\Phi$ is not an
immersion, and $N$ has singular points. Also, under certain circumstances
$\Phi$ may be periodic in $t$, and then $N$ will be generically
diffeomorphic to ${\cal S}^1\t\R^2$ rather than~$\R^3$.

The basic details of the construction are explained in
\S\ref{ev5}, in terms of o.d.e.s for vector-valued functions
${\bf z}_1,\ldots,{\bf z}_6:\R\ra\C^3$. We discuss the symmetries
of the construction, and show that the family of SL 3-folds we have
constructed, up to automorphisms of $\C^3$, is 9-dimensional.
Section \ref{ev6} studies and describes the singularities of the
corresponding SL 3-folds, which we believe are of a new kind.

The goal of sections \ref{ev7}--\ref{ev11} is to solve the o.d.e.s
for ${\bf z}_1,\ldots,{\bf z}_6$ as explicitly as we are able to,
and so to write our examples of SL 3-folds as explicitly as possible.
To do this we split into several cases, and use the symmetries of
the problem to write each case in a convenient form. In \S\ref{ev7}
we divide into four cases (i)--(iv) of increasing complexity,
depending on the rank of the homogeneous quadratic part of $\Phi$.
Cases (i) and (ii) are easy and are dealt with at once. Case (iii)
is the subject of \S\ref{ev8}, and case (iv) is divided into
subcases and discussed in~\S\ref{ev9}--\S\ref{ev11}.

In Theorems \ref{ev8thm2}, \ref{ev10thm2} and \ref{ev11thm1} we are
able to write down three families of special Lagrangian 3-folds in
$\C^3$ very explicitly, and these are some of the main results of
the paper. Also, in \S\ref{ev112} we study the condition for the
family $\{\Phi_t:t\in\R\}$ to be periodic in a special case. This
has a surprisingly abundant and structured set of solutions, which
leads in Theorem \ref{ev11thm2} to a countable set of distinct
families of immersed SL 3-folds diffeomorphic to~${\cal S}^1\t\R^2$.

Section \ref{ev12} then studies the series of families of SL
3-folds arising from the sets of affine evolution data given
in Example \ref{ev4ex2}. The development follows parts of
\S\ref{ev5}--\S\ref{ev11} closely, and so we leave out some of
the details. For each $k=1,2,\ldots$ we construct a family of
immersed SL 3-folds in $\C^3$ diffeomorphic to $\R^3$, which
can be written down in parametric form entirely explicitly.

For $k=1$ they are isomorphic to the SL 3-folds of
\cite[Ex.~7.5]{Joyc3}, and for $k=2$ they are isomorphic to the
family studied in \S\ref{ev8}. When $k\ge 3$ these families
include many nontrivial periodic solutions, yielding families of
immersed SL 3-folds in $\C^3$ diffeomorphic to~${\cal S}^1\t\R^2$. 

\section{Special Lagrangian submanifolds in $\C^m$}
\label{ev2}

We begin by defining {\it calibrations} and {\it calibrated 
submanifolds}, following Harvey and Lawson~\cite{HaLa}.

\begin{dfn} Let $(M,g)$ be a Riemannian manifold. An {\it oriented
tangent $k$-plane} $V$ on $M$ is a vector subspace $V$ of
some tangent space $T_xM$ to $M$ with $\dim V=k$, equipped
with an orientation. If $V$ is an oriented tangent $k$-plane
on $M$ then $g\vert_V$ is a Euclidean metric on $V$, so 
combining $g\vert_V$ with the orientation on $V$ gives a 
natural {\it volume form} $\vol_V$ on $V$, which is a 
$k$-form on~$V$.

Now let $\vp$ be a closed $k$-form on $M$. We say that
$\vp$ is a {\it calibration} on $M$ if for every oriented
$k$-plane $V$ on $M$ we have $\vp\vert_V\le \vol_V$. Here
$\vp\vert_V=\al\cdot\vol_V$ for some $\al\in\R$, and 
$\vp\vert_V\le\vol_V$ if $\al\le 1$. Let $N$ be an 
oriented submanifold of $M$ with dimension $k$. Then 
each tangent space $T_xN$ for $x\in N$ is an oriented
tangent $k$-plane. We say that $N$ is a {\it calibrated 
submanifold} if $\vp\vert_{T_xN}=\vol_{T_xN}$ for all~$x\in N$.
\label{ev2def1}
\end{dfn}

It is easy to show that calibrated submanifolds are automatically
{\it minimal submanifolds} \cite[Th.~II.4.2]{HaLa}. Here is the 
definition of special Lagrangian submanifolds in $\C^m$, taken
from~\cite[\S III]{HaLa}.

\begin{dfn} Let $\C^m$ have complex coordinates $(z_1,\dots,z_m)$, 
and define a metric $g$, a real 2-form $\om$ and a complex $m$-form 
$\Om$ on $\C^m$ by
\begin{align*}
g=\ms{\d z_1}+\cdots+\ms{\d z_m},\quad
\om&=\frac{i}{2}(\d z_1\w\d\bar z_1+\cdots+\d z_m\w\d\bar z_m),\\
\text{and}\quad\Om&=\d z_1\w\cdots\w\d z_m.
\end{align*}
Then $\Re\Om$ and $\Im\Om$ are real $m$-forms on $\C^m$. Let
$L$ be an oriented real submanifold of $\C^m$ of real dimension 
$m$. We say that $L$ is a {\it special Lagrangian submanifold} 
of $\C^m$ if $L$ is calibrated with respect to $\Re\Om$, in the 
sense of Definition \ref{ev2def1}. We will often abbreviate 
`special Lagrangian' by `SL', and `$m$-dimensional submanifold' 
by `$m$-fold', so that we shall talk about SL $m$-folds in~$\C^m$. 
\end{dfn}

As in \cite{Joyc1,Joyc2} there is also a more general definition 
of special Lagrangian submanifolds involving a {\it phase} 
${\rm e}^{i\th}$, but we will not use it in this paper. Harvey 
and Lawson \cite[Cor.~III.1.11]{HaLa} give the following 
alternative characterization of special Lagrangian submanifolds.

\begin{prop} Let\/ $L$ be a real $m$-dimensional submanifold 
of $\C^m$. Then $L$ admits an orientation making it into an
SL submanifold of\/ $\C^m$ if and only if\/ $\om\vert_L\equiv 0$ 
and\/~$\Im\Om\vert_L\equiv 0$.
\end{prop}

Note that an $m$-dimensional submanifold $L$ in $\C^m$ is 
called {\it Lagrangian} if $\om\vert_L\equiv 0$. Thus special 
Lagrangian submanifolds are Lagrangian submanifolds satisfying 
the extra condition that $\Im\Om\vert_L\equiv 0$, which is how 
they get their name.

\section{Review of the `evolution' construction of \protect\cite{Joyc3}}
\label{ev3}

We now review the construction of special Lagrangian $m$-folds
in $\C^m$ given by the author in \cite[\S 3]{Joyc3}, which will be
used in \S\ref{ev5} and \S\ref{ev12} to construct the SL 3-folds we
are interested in. There are two versions of the construction,
{\it linear} and {\it affine} (linear plus constant), and we shall be
using the affine version. The construction depends on some {\it evolution
data}, which we now define, following~\cite[Def.~3.4]{Joyc3}.

\begin{dfn} Let $2\le m\le n$ be integers. A set of 
{\it affine evolution data} is a pair $(P,\chi)$, where 
$P$ is an $(m-1)$-dimensional submanifold of $\R^n$, and 
$\chi:\R^n\ra\La^{m-1}\R^n$ is an affine map, such that 
$\chi(p)$ is a nonzero element of $\La^{m-1}TP$ in 
$\La^{m-1}\R^n$ for each nonsingular $p\in P$. We suppose 
also that $P$ is not contained in any proper affine subspace 
$\R^k$ of~$\R^n$.

Let $\Aff(\R^n,\C^m)$ be the affine space of affine maps 
$\phi:\R^n\ra\C^m$, and define ${\cal C}_P$ to be the subset 
of $\phi\in\Aff(\R^n,\C^m)$ satisfying
\begin{itemize}
\item[(i)] $\phi^*(\om)\vert_P\equiv 0$, and
\item[(ii)] $\phi\vert_{T_pP}:T_pP\ra\C^m$ is injective 
for all $p$ in a dense open subset of~$P$.
\end{itemize}
Then ${\cal C}_P$ is nonempty, and is an open set in the
intersection of a finite number of quadrics in~$\Aff(\R^n,\C^m)$. 
\label{ev3def}
\end{dfn}

We may define {\it linear evolution data} in the same way, but
using linear rather than affine maps. With this definition, the
construction is contained in the following theorem, taken
from~\cite[Th.~3.5]{Joyc3}.

\begin{thm} Let $(P,\chi)$ be a set of affine evolution data, 
and $n,m,\Aff(\R^n\!,\C^m)$ and\/ ${\cal C}_P$ be as above.
Suppose $\phi\in{\cal C}_P$. Then there exists $\ep>0$ and a 
unique real analytic family $\bigl\{\phi_t:t\in(-\ep,\ep)\bigr\}$ 
in ${\cal C}_P$ with\/ $\phi_0=\phi$, satisfying the equation
\e
\left(\frac{\d\phi_t}{\d t}(x)\right)^b=(\phi_t)_*(\chi(x))^{a_1
\ldots a_{m-1}}(\Re\Om)_{a_1\ldots a_{m-1}a_m}g^{a_mb}
\label{ev3eq1}
\e
for all\/ $x\in\R^n$, using the index notation for tensors in $\C^m$.
Furthermore, $N=\bigl\{\phi_t(p):t\in(-\ep,\ep)$, $p\in P\bigr\}$ is 
a special Lagrangian submanifold in $\C^m$ wherever it is nonsingular. 
\label{ev3thm1}
\end{thm}

Here is a brief explanation of the theorem and its proof. Equation
\eq{ev3eq1} is a first-order o.d.e.\ upon $\phi_t$. The key point
to note is that as $\chi$ is affine, the right hand side of
\eq{ev3eq1} is affine in $x$, and so \eq{ev3eq1} makes sense as
an evolution equation for affine maps $\phi_t$. However, the right
hand side of \eq{ev3eq1} is a homogeneous polynomial of order $m-1$
in $\phi_t$, so for $m>2$ it is a {\it nonlinear}~o.d.e.

The special Lagrangian $m$-fold $N$ is the total space of a
1-parameter family of real $(m\!-\!1)$-dimensional submanifolds
$\phi_t(P)$ of $\C^m$, each of which is an affine image of the
$(m\!-\!1)$-manifold $P$ in $\R^n$. Thus we can think of
\eq{ev3eq1} as an evolution equation in a certain class of real
$(m\!-\!1)$-submanifolds of~$\C^m$.

The theorem is proved using the following result, taken
from~\cite[Th.~3.3]{Joyc2}.

\begin{thm} Let\/ $P$ be a compact, orientable, real analytic 
$(m-1)$-manifold, $\chi$ a real analytic, nonvanishing section 
of\/ $\La^{m-1}TP$, and\/ $\phi:P\ra\C^m$ a real analytic
embedding (immersion) such that\/ $\phi^*(\om)\equiv 0$ on 
$P$. Then there exists $\ep>0$ and a unique family 
$\bigl\{\phi_t:t\in(-\ep,\ep)\bigr\}$ of real analytic maps 
$\phi_t:P\ra\C^m$ with $\phi_0=\phi$, satisfying the equation
\e
\left(\frac{\d\phi_t}{\d t}\right)^b=(\phi_t)_*(\chi)^{a_1\ldots a_{m-1}}
(\Re\Om)_{a_1\ldots a_{m-1}a_m}g^{a_mb},
\label{ev3eq2}
\e
using the index notation for (real) tensors on $\C^m$. Define 
$\Phi:(-\ep,\ep)\t P\ra\C^m$ by $\Phi(t,p)=\phi_t(p)$. Then 
$N=\Image\Phi$ is a nonsingular embedded (immersed) special 
Lagrangian submanifold of\/~$\C^m$.
\label{ev3thm2}
\end{thm}

This constructs SL $m$-folds in $\C^m$ by evolving arbitrary
(compact) real analytic $(m\!-\!1)$-submanifolds $P$ of $\C^m$
with $\om\vert_P\equiv 0$. The trouble with this result is that
as the set of such submanifolds is infinite-dimensional, the
theorem is really an infinite-dimensional evolution problem,
and so is difficult to solve explicitly.

What we achieve in Definition \ref{ev3def} and Theorem \ref{ev3thm1} is
to find a special class $\cal C$ of real analytic $(m-1)$-submanifolds
$P$ of $\C^m$ with $\om\vert_P\equiv 0$, depending on finitely many
real parameters $c_1,\ldots,c_n$, such that the evolution equation
\eq{ev3eq2} stays within the class~$\cal C$.

In fact \eq{ev3eq2} reduces to \eq{ev3eq1}, which is basically the
same equation, but is now a first order o.d.e.\ on $c_1,\ldots,c_n$,
as functions of $t$. Thus we have reduced the infinite-dimensional
problem of evolving submanifolds in $\C^m$ to a {\it finite-dimensional}
o.d.e., which we may be able to solve explicitly.

\section{A geometric interpretation of evolution data with $m=3$}
\label{ev4}

In \cite[Th.~4.9]{Joyc3} the author showed that every set of
evolution data $(P,\chi)$ in $\R^n$ admits a locally transitive
symmetry group $G$ in $\GL(n,\R)$, and that when $m=3$ there is
a $G$-invariant surjective map $(\R^n)^*\ra\g$ with kernel 0 or
$\R$, where $\g$ is the Lie algebra of $G$. Motivated by this,
we shall now present a correspondence between sets of linear
or affine evolution data $(P,\chi)$ with $m=3$, and symplectic
2-manifolds $(\Si,\om)$ with a transitive Hamiltonian symmetry group.

Let $(\Si,\om)$ be a symplectic 2-manifold, not necessarily compact,
and $G$ a connected Lie group with Lie algebra $\g$ acting faithfully
and transitively on $\Si$. Suppose that $G$ preserves $\om$ and every
element of $\g$ admits a moment map; this is called a {\it Hamiltonian
action}, and holds automatically if $\Si$ is simply-connected. Define
$V$ to be the vector space of moment maps of elements of $\g$,
including constant functions.

That is, $V$ is the vector space of smooth maps $f:\Si\ra\R$ such that 
$\d f=x\cdot\om$ for some $x\in\g$. Then $V\cong\R\op\g$. Define 
$\psi:\Si\ra V^*$ by $\psi(x)\cdot f=f(x)$ for all $f\in V$ and 
$x\in \Si$. As $G$ acts transitively on $\Si$, one can show that $\psi$ 
is an immersion. Let $P=\psi(\Si)$, so that $P$ is an immersed 
2-submanifold in~$V^*$.

Now the {\it Poisson bracket} on $(\Si,\om)$ yields a natural 
bilinear, antisymmetric map $\{\,,\,\}:V\t V\ra V$ given in 
index notation by $\{f,f'\}=\om^{ab}(\d f)_a(\d f')_b$, where
$\om^{ab}$ is the inverse of $\om_{ab}$. This makes $V$ into a
Lie algebra, which is an extension of the Lie algebra $\g$ by~$\R$.

Thus $P$ is a submanifold in the dual of a Lie algebra. In fact 
$P$ is a {\it coadjoint orbit}, that is, an orbit of the coadjoint 
action on $V^*$ of the connected, simply-connected Lie group associated 
to $V$. It is well known that all coadjoint orbits have a natural 
symplectic structure.

As the Poisson bracket is bilinear and antisymmetric, we can extend 
it to a linear map $\{\,,\,\}:\La^2V\ra V$. Define $\chi:V^*\ra\La^2V^*$ 
to be the dual of this linear map. Using the fact that $G$ acts
transitively on $\Si$, it is not difficult to show that $\chi(p)$ is a 
nonzero element of $\La^2T_pP\subset\La^2V^*$ for each $p\in P$. Thus
$(P,\chi)$ is a set of {\it linear evolution data} in the vector space~$V^*$.

Actually, it is usually nicer to regard $(P,\chi)$ as {\it affine} 
evolution data, in the following way. Let $\bf 1$ be the constant 
function 1 on $\Si$. Then ${\bf 1}\in V$. Define $f:V^*\ra\R$ by 
$f(\al)=\al({\bf 1})$. Then $U=f^{-1}(1)$ is a hyperplane in $V^*$, 
which contains $P$. We can regard $U$ as an affine space modelled on 
$\g^*$. The restriction of $\chi$ to $U$ is an affine map 
$\chi:U\ra\La^2U$, and $(P,\chi)$ is a set of affine evolution data 
in the affine space~$U$.

We have shown that given a symplectic 2-fold $(\Si,\om)$ with a
faithful, transitive, Hamiltonian action of a Lie group $G$, we
can construct sets of linear and affine evolution data with $m=3$.
We now explain how to reverse this construction, so that starting
with a set of evolution data with $m=3$, we construct a symplectic
2-manifold $(\Si,\om)$ and group action. For simplicity we work with
the linear case, as affine evolution data in $\R^n$ can always be
reduced to linear evolution data in~$\R^{n+1}$.

Let $(P',\chi')$ be a set of linear evolution data with $m=3$ in a real
vector space $W$. That is, $P'$ is a connected 2-submanifold of $W$,
lying in no proper vector subspace of $W$, and $\chi':W\ra\La^2W$ a
linear map such that $\chi'(p)\in\La^2T_pP'\sm\{0\}$ for each $p\in P'$.
Let $\om'$ be the symplectic structure on $P'$ dual to $\chi'\vert_{P'}$.
Then $(P',\om')$ is a symplectic 2-manifold. We shall define an
antisymmetric, bilinear bracket $[\,,\,]:W^*\t W^*\ra W^*$ which
makes $W^*$ into a Lie algebra.

Regard $\chi'$ as an element of $W^*\ot\La^2 W$, and for each
$\al,\be\in W^*$, define $[\al,\be]=\chi'\cdot(\al\w\be)$, where
`$\,\cdot\,$' is the natural pairing between $\La^2W^*$ and $\La^2W$.
Now $W^*$ is the vector space of linear maps $W\ra\R$, so elements of
$W^*$ give real functions on $P'$ by restriction. As $P'$ is contained
in no proper subspace of $W$, this map from $W^*$ to functions on $P'$
is injective, so $W^*$ may be viewed as a vector space of functions
on $P'$. Thought of in this way, it is easy to show that the bracket
$[\,,\,]$ on $W^*$ is actually the {\it Poisson bracket} on functions
on $P'$, induced by the symplectic structure~$\om'$.

But the Poisson bracket automatically satisfies the Jacobi identity.
Thus $[\,,\,]$ makes $W^*$ into a {\it Lie algebra}. To each element
of $W^*$ we associate its Hamiltonian vector field, giving a map
$W^*\ra\Vect(P')$, where $\Vect(P')$ is the smooth vector fields
on $P'$. This is a {\it Lie algebra automorphism}, with respect to
the usual Lie bracket of vector fields on $\Vect(P')$. Let $\g$ be
the image of $W^*$ in $\Vect(P')$. Then $\g$ is a finite-dimensional
Lie subalgebra of~$\Vect(P')$.

It can be shown that either
\begin{itemize}
\item[(a)] $P'$ is contained in no affine hyperplane in $W$
and $\g\cong W^*$, or
\item[(b)] $P'$ is contained in an affine hyperplane in $W$
and $\g\cong W^*/\R$, where $\R$ is an ideal in $W^*$. Also,
$(P',\chi')$ reduces to a set of affine evolution data in one
fewer dimension.
\end{itemize}

By Cartan's theorems, there exists a unique connected,
simply-connected Lie group $G$ with Lie algebra $\g$. Choose
$p\in P'$, and let $\h$ be the vector subspace of vector fields
in $\g$ that vanish at $p$. Then $\h$ is a Lie subalgebra of
$\g$, and so corresponds to a unique connected Lie subgroup $H$
of~$G$.

As we can take two functions in $W^*$ to be coordinates near $p$,
the corresponding vectors in $T_pP'$ are linearly independent, and
so $\dim\h=\dim\g-2$. Also, there is a natural isomorphism
$T_pP'\cong\g/\h$. Since the vector fields in $\g$ are Hamiltonian,
they preserve $\om'$. Thus the adjoint action of $\h$ on $\g/\h$
preserves the nonzero 2-form $\om'_p$ on $T_pP'\cong\g/\h$.
As $H$ is connected, it follows that $H$ also preserves~$\om'_p$.

It can be shown that $H$ is closed in $G$. Then $\Si=G/H$ is a
connected 2-manifold with a natural $G$-action. The tangent
space $T_H\Si$ is isomorphic to $\g/\h\cong T_pP'$, and so has
a nonzero 2-form $\om'_p$. As $\om'_p$ is $H$-invariant, this
extends to a $G$-invariant, nonvanishing 2-form $\om$ on~$\Si$.

This makes $\Si$ into a symplectic 2-manifold with a transitive
$G$-action preserving $\om$. As $G$ is simply-connected, so is
$\Si$. Thus, every element of $\g$ admits a moment map for its
action on $\Si$. Therefore by the construction at the beginning
of this section, we can associate a set of linear evolution
data $(P,\chi)$ in a vector space $V^*$ to $(\Si,\om)$ and~$G$.

One can prove that $W$ is naturally isomorphic to $V^*$, and that
this isomorphism identifies $P'$ with an open subset of $P$, and
$\chi'$ with $\chi$. The details are left to the reader. To sum up,
we have proved the following theorem:

\begin{thm} Let\/ $(\Si,\om)$ be a symplectic $2$-manifold, and\/
$G$ a connected Lie group with a faithful, transitive, Hamiltonian
action on $\Si$. Then we can construct sets of linear and affine
evolution data $(P,\chi)$ with\/ $m=3$ and\/ $P\cong \Si$. Conversely,
every set of linear or affine evolution data with\/ $m=3$ locally
arises from this construction.
\label{ev4thm}
\end{thm}

All of the quadric examples of \cite[\S 4.1]{Joyc3} for $m=3$ can be easily
extracted from this construction. For example:
\begin{itemize}
\item[(i)] Let $\Si$ be ${\cal S}^2$, and $G$ be SO(3) acting by isometries.
Then $P$ is the sphere $x_1^2+x_2^2+x_3^2=1$ in~$\R^3$.
\item[(ii)] Let $\Si$ be the hyperbolic plane ${\cal H}^2$, and $G$ be
$SO(1,1)_+$ acting by isometries. Then $P$ is half of the hyperboloid
$x_1^2-x_2^2-x_3^2=1$ in~$\R^3$.
\item[(iii)] Let $\Si$ be $\R^2\sm\{0\}$ with the standard volume form, and
$G$ be $\SL(2,\R)$ acting as usual. Then $P$ is one of the pair of 
cones $x_1^2-x_2^2-x_3^2=0$ in~$\R^3$.
\item[(iv)] Let $\Si$ be $\R^2$ and $G$ be the group of Euclidean
transformations $\SO(2)\lt\R^2$. Then $P$ is the paraboloid
$x_1^2+x_2^2+x_3=0$ in~$\R^3$.
\end{itemize}

But we can also find interesting new sets of evolution data which are
not quadrics. For instance, when we take $\Si$ to be $\R^2$ with its
standard volume form, and $G$ to be $\SL(2,\R)\lt\R^2$ acting on $\Si$
by affine transformations, we get the following example.

\begin{ex} Consider the map $\psi:\R^2\ra\R^5$ given by
\begin{equation*}
\psi:(y_1,y_2)\mapsto \bigl(\ha(y_1^2+y_2^2),\ha(y_1^2-y_2^2),y_1y_2,
y_1,y_2\bigr).
\end{equation*}
The image of $\psi$ is
\begin{equation*}
P=\bigl\{(x_1,\ldots,x_5)\in\R^5:x_1=\ha(x_4^2+x_5^2),\;
x_2=\ha(x_4^2-x_5^2),\; x_3=x_4x_5\bigr\},
\end{equation*}
which is diffeomorphic to $\R^2$. Writing $e_j=\frac{\pd}{\pd x_j}$,
calculation shows that
\begin{equation*}
\psi_*\bigl({\ts\frac{\pd}{\pd y_1}}\bigr)=y_1e_1+y_1e_2+y_2e_3+e_4,\quad 
\psi_*\bigl({\ts\frac{\pd}{\pd y_2}}\bigr)=y_2e_1-y_2e_2+y_1e_3+e_5,
\end{equation*}
and therefore
\begin{align*}
\psi_*&\bigl({\ts\frac{\pd}{\pd y_1}\w\frac{\pd}{\pd y_2}}\bigr)=
(y_1^2+y_2^2)e_2\w e_3+(y_1^2-y_2^2)e_1\w e_3-2y_1y_2e_1\w e_2\\
&+y_1(e_1\w e_5+e_2\w e_5-e_3\w e_4)+y_2(-e_1\w e_4+e_2\w e_4+e_3\w e_5)\\
&+e_4\w e_5.
\end{align*}
Thus if we define an affine map $\chi:\R^5\ra\La^2\R^5$ by
\begin{align*}
&\chi(x_1,\ldots,x_5)=2x_1e_2\w e_3+2x_2e_1\w e_3-2x_3e_1\w e_2\\
&+x_4(e_1\w e_5+e_2\w e_5-e_3\w e_4)+x_5(-e_1\w e_4+e_2\w e_4+e_3\w e_5)\\
&+e_4\w e_5,
\end{align*}
then $\chi=\psi_*\bigl(\frac{\pd}{\pd y_1}\w\frac{\pd}{\pd y_2}\bigr)$ on 
$P$. This implies that $(P,\chi)$ is a set of affine evolution data
with $m=3$ and $n=5$, which does not arise from the construction 
of \cite[\S 4.1]{Joyc3}. So applying Theorem \ref{ev3thm1} will give a
family of special Lagrangian 3-folds in $\C^3$. These will be
studied at length in~\S\ref{ev5}--\S\ref{ev11}.
\label{ev4ex1}
\end{ex}

Now if a Lie group $G$ acts transitively on a symplectic 2-manifold
$\Si$ then often a Lie subgroup $G'$ of $G$ will also act transitively
on $\Si$, or on some open subset $\Si'$ of $\Si$. In our next result
we consider the relation between the families of SL 3-folds constructed
using $\Si,G$ and~$\Si',G'$.

\begin{prop} Let\/ $(\Si,\om)$ be a symplectic $2$-manifold, and\/
$G$ a connected Lie group with a faithful, transitive, Hamiltonian
action on $\Si$. Suppose $G'$ is a connected Lie subgroup of\/ $G$,
and $\Si'$ an open orbit of\/ $G'$ in $\Si$. Then the special Lagrangian 
$3$-folds in $\C^3$ constructed using $\Si$ and\/ $G$ by combining 
Theorem \ref{ev4thm} and the method of\/ \S\ref{ev3} include all 
those constructed using $\Si'$  and\/~$G'$. 
\label{ev4prop}
\end{prop}

\begin{proof} The construction above gives sets of linear evolution 
data $(P,\chi)$ and $(P'\chi')$ from $\Si,G$ and $\Si',G'$, where $P,P'$
lie in vector spaces $V^*,(V')^*$. It is easy to see that $V'$ is
a vector subspace of $V$, since the Lie algebra $\g'$ of $G'$ is 
a vector subspace of $\g$. Let $U$ be the {\it annihilator} 
$(V')^\circ$ of $V'$ in $V^*$. Then~$(V')^*\cong V^*/U$. 

Now the construction of SL 3-folds using $\Si',G'$ in \S\ref{ev3} 
involves a 1-parameter family of linear maps $\phi_t':(V')^*\ra\C^3$
in ${\cal C}_{P'}$ satisfying the o.d.e.\ \eq{ev3eq1}. Let $\phi_t:
V^*\ra\C^3$ be the pull-back of $\phi_t'$ from $(V')^*=V^*/U$ to $V^*$. 
It is easy to show that this family $\phi_t$ also lie in ${\cal C}_P$ 
and satisfy \eq{ev3eq1}, and that the SL 3-fold $N'$ constructed using 
the $\phi_t'$ is a subset of the SL 3-fold $N$ constructed using 
the $\phi_t$. Thus the SL 3-folds constructed using $\Si',G'$ are 
included in those constructed using~$\Si,G$. 
\end{proof}

In particular, the family of SL 3-folds coming from Example \ref{ev4ex1},
corresponding to the action of $\SL(2,\R)\lt\R^2$, will include families
of SL 3-folds corresponding to subgroups of $\SL(2,\R)\lt\R^2$. For 
example, as in part (iii) above we can set $\Si'=\R^2\sm\{0\}$ and
$G'=\SL(2,\R)$, and as in part (iv) above we can set $\Si'=\R^2$ and 
$G'=\SO(2)\lt\R^2$. So the families of SL 3-folds corresponding to parts 
(iii) and (iv), which we have already considered in \cite{Joyc3}, will
occur as special cases in the family of SL 3-folds to be studied
in~\S\ref{ev5}--\S\ref{ev11}.

Above we set $\Si=G/H$, and used the fact that $H$ acts naturally
on $T_H\Si=\g/\h$ preserving $\om_H$. It is tempting to assume that
this action of $H$ on $T_H\Si$ is {\it faithful}, as would be the case
in Riemannian rather than symplectic geometry. If this held then $H$ would
be a subgroup of $\SL(2,\R)$, and would be largest in Example~\ref{ev4ex1}.

However, $H$ need not act faithfully on $T_H\Si$, and in fact $G,H$
can have arbitrarily large dimension. Here is a class of examples
in which this happens. Let $C$ be $\R$ or ${\cal S}^1$, let $K$ be
a connected Lie group acting smoothly, transitively and faithfully
on $C$, and let $U$ be a nonzero vector space of 1-forms on $C$
which is invariant under~$K$.

Define $\Si$ to be $T^*C$ with its canonical symplectic structure $\om$,
and $G$ to be the semidirect product $K\lt U$ acting on $\Si$ by
\e
(x,y\,\d x)\,\,{\buildrel(\ka,u)\over\longmapsto}\,\,
\bigl(\ka(x),{\ts\frac{\d\ka}{\d x}}(x)y\,\d x+u(x)\bigr).
\label{ev4eq1}
\e
Here we write a point in $\Si$ as $(x,y\,\d x)$, where $x$ is a coordinate
in $C$ with values in $\R$ or $\R/\Z$, and $y\,\d x$ lies in $T^*_xC$, so
that $y\in\R$. Elements of $K\lt U$ are written $(\ka,u)$ for $\ka\in K$
and $u\in U$, so that $\ka:C\ra C$ is a differentiable map. It is easy to
see that \eq{ev4eq1} defines the action of a Lie group $G=K\lt U$ on $\Si$,
which is faithful and transitive and preserves~$\om$.

The possibilities for $C$ and $K$ are
\begin{itemize}
\item[(i)] $C=\R$ and $K=\R$ acting by translations.
\item[(ii)] $C=\R$ and $K=\R_+\lt\R$, acting by $x\,{\buildrel(a,b)\over
\longmapsto}\,ax+b$ for $a>0$ and~$b\in\R$.
\item[(iii)] $C={\cal S}^1$, thought of as $\U(1)$, and $K=\U(1)$
acting by multiplication.
\item[(iv)] $C={\cal S}^1$, thought of as $\mathbb{RP}^1$, and
$K={\mathbin{\rm PSL}}(2,\R)$ acting by projective transformations.
\end{itemize}

In case (iv), one can show that there are no non-zero, finite-dimensional,
$K$-invariant vector spaces of 1-forms $U$ on $C$, so we rule this case
out. In cases (i) and (iii) there are many suitable spaces of 1-forms $U$,
and so we may construct many sets of evolution data $(P,\chi)$ in $\R^n$.
However, calculation shows that we may always choose coordinates
$(x_1,\ldots,x_n)$ on $\R^n$, and so split $\R^n=\R^{n-1}\t\R$, such that
$\chi=v\w\frac{\pd}{\pd x_n}$, where $v$ is a linear or affine vector
field in $\R^{n-1}$, and $P=\ga\t\R$, where $\ga$ is an integral curve
of $v$ in~$\R^{n-1}$.

Thus, cases (i) and (iii) yield evolution data resulting from combining
Examples 4.5 and 4.6 of \cite{Joyc3}. The corresponding SL 3-folds will
all split as products $\Si\t\R$ in $\C^2\t\C$, where $\Si$ is an SL
2-fold in $\C^2$. We are not interested in such examples, so we rule
these cases out too.

This leaves case (ii). Here the natural candidates for $U$ are
\begin{equation*}
U_k=\bigl\{p(x)\,\d x:p(x)\text{ is a real polynomial of degree }<k\bigr\},
\end{equation*}
for $k\ge 1$. In the following example we define the corresponding
set of affine evolution data in $\R^{k+2}$, yielded by the construction
of Theorem~\ref{ev4thm}.

\begin{ex} Choose $k\ge 1$, and let $(x_1,\ldots,x_k,y_1,y_2)$ be 
coordinates on $\R^{k+2}$. Define a map $\psi:\R^2\ra\R^{k+2}$ by
$\psi:(x,y)\longmapsto (x,x^2,\ldots,x^k,y,xy)$. The image of $\psi$ is
\begin{equation*}
P=\bigl\{(x_1,\ldots,x_k,y_1,y_2)\in\R^{k+2}:
\text{$x_j=(x_1)^j$ for $j=2,\ldots,k$, $y_2=y_1x_1$}\bigr\},
\end{equation*}
which is diffeomorphic to $\R^2$. Calculation shows that
\begin{equation*}
\psi_*\bigl({\ts\frac{\pd}{\pd x}}\bigr)=y{\ts\frac{\pd}{\pd y_2}}
\!+\!{\ts\frac{\pd}{\pd x_1}}+2x{\ts\frac{\pd}{\pd x_2}}
\!+\!\cdots\!+\!kx^{k-1}{\ts\frac{\pd}{\pd x_k}}
\quad\text{and}\quad
\psi_*\bigl({\ts\frac{\pd}{\pd y}}\bigr)=
{\ts\frac{\pd}{\pd y_1}}\!+\!x{\ts\frac{\pd}{\pd y_2}},
\end{equation*}
and therefore
\begin{align*}
\psi_*\bigl({\ts\frac{\pd}{\pd x}}\w{\ts\frac{\pd}{\pd y}}\bigr)
&=-y{\ts\frac{\pd}{\pd y_1}}\w{\ts\frac{\pd}{\pd y_2}}
\!+\!{\ts\frac{\pd}{\pd x_1}}\w{\ts\frac{\pd}{\pd y_1}}
\!+\!2x{\ts\frac{\pd}{\pd x_2}}\w{\ts\frac{\pd}{\pd y_1}}\!+\!\cdots\!+\!
kx^{k-1}{\ts\frac{\pd}{\pd x_k}}\w{\ts\frac{\pd}{\pd y_1}}\\
&\quad +\!x{\ts\frac{\pd}{\pd x_1}}\w{\ts\frac{\pd}{\pd y_2}}
\!+\!2x^2{\ts\frac{\pd}{\pd x_2}}\w{\ts\frac{\pd}{\pd y_2}}\!+\!\cdots
\!+\!kx^k{\ts\frac{\pd}{\pd x_k}}\w{\ts\frac{\pd}{\pd y_2}}\\
&=-y_1{\ts\frac{\pd}{\pd y_1}}\!\w\!{\ts\frac{\pd}{\pd y_2}}
\!+\!{\ts\frac{\pd}{\pd x_1}}\!\w\!{\ts\frac{\pd}{\pd y_1}}
\!+\!2x_1{\ts\frac{\pd}{\pd x_2}}\!\w\!{\ts\frac{\pd}{\pd y_1}}\!+\!\cdots
\!+\!kx_{k-1}{\ts\frac{\pd}{\pd x_k}}\!\w\!{\ts\frac{\pd}{\pd y_1}}\\
&\quad +\!x_1{\ts\frac{\pd}{\pd x_1}}\w{\ts\frac{\pd}{\pd y_2}}
\!+\!2x_2{\ts\frac{\pd}{\pd x_2}}\w{\ts\frac{\pd}{\pd y_2}}\!+\!\cdots
\!+\!kx_k{\ts\frac{\pd}{\pd x_k}}\w{\ts\frac{\pd}{\pd y_2}}.
\end{align*}
Thus if we define an affine map $\chi:\R^{k+2}\ra\La^2\R^{k+2}$ by
\begin{align*}
\chi(x_1,&\ldots,x_k,y_1,y_2)\\
&=-2y_1{\ts\frac{\pd}{\pd y_1}}\w{\ts\frac{\pd}{\pd y_2}}
+2{\ts\frac{\pd}{\pd x_1}}\w{\ts\frac{\pd}{\pd y_1}}
+4x_1{\ts\frac{\pd}{\pd x_2}}\w{\ts\frac{\pd}{\pd y_1}}+\cdots
+2kx_{k-1}{\ts\frac{\pd}{\pd x_k}}\w{\ts\frac{\pd}{\pd y_1}}\\
&\quad+2x_1{\ts\frac{\pd}{\pd x_1}}\w{\ts\frac{\pd}{\pd y_2}}
+4x_2{\ts\frac{\pd}{\pd x_2}}\w{\ts\frac{\pd}{\pd y_2}}+\cdots
+2kx_k{\ts\frac{\pd}{\pd x_k}}\w{\ts\frac{\pd}{\pd y_2}},
\end{align*}
then $\chi=2\psi_*\bigl(\frac{\pd}{\pd x}\w\frac{\pd}{\pd y}\bigr)$ on 
$P$. This implies that $(P,\chi)$ is a set of affine evolution data with
$m=3$ and $n=k+2$. So applying Theorem \ref{ev3thm1} gives a family of 
special Lagrangian 3-folds in $\C^3$, which will be studied in~\S\ref{ev12}.
\label{ev4ex2}
\end{ex}

\section{A construction of SL 3-folds in $\C^3$}
\label{ev5}

We now apply the `evolution equation' construction of \S\ref{ev3} to
the set of affine evolution data defined in Example \ref{ev4ex1}. As
in Example \ref{ev4ex1}, let $P$ be the image in $\R^5$ of the map
$\psi:\R^2\ra\R^5$ given by
\e
\psi:(y_1,y_2)\mapsto \bigl(\ha(y_1^2+y_2^2),\ha(y_1^2-y_2^2),y_1y_2,
y_1,y_2\bigr),
\label{ev5eq1}
\e
and define $\chi:\R^5\ra\La^2\R^5$ by
\e
\begin{split}
&\chi(x_1,\ldots,x_5)=2x_1e_2\w e_3+2x_2e_1\w e_3-2x_3e_1\w e_2\\
&+x_4(e_1\!\w\!e_5\!+\!e_2\!\w\!e_5\!-\!e_3\!\w\!e_4)
+x_5(-e_1\!\w\!e_4\!+\!e_2\!\w\!e_4\!+\!e_3\!\w\!e_5)\\
&+e_4\w e_5,
\end{split}
\label{ev5eq2}
\e
where $e_j=\frac{\pd}{\pd x_j}$. Then $(P,\chi)$ is a set of affine
evolution data.

Let ${\bf z}_1,\ldots,{\bf z}_6$ be vectors in $\C^3$, and define an
affine map $\phi:\R^5\ra\C^3$ by
\e
\phi:(x_1,\ldots,x_5)\mapsto x_1{\bf z}_1+\cdots+x_5{\bf z}_5+{\bf z}_6.
\label{ev5eq3}
\e
Then as $\phi_*(e_k)={\bf z}_k$ for $k\le 5$, from \eq{ev5eq2} we see that 
\begin{align*}
\phi^*(\om)\cdot\chi&=
2x_1\om({\bf z}_2,{\bf z}_3)+2x_2\om({\bf z}_1,{\bf z}_3)
-2x_3\om({\bf z}_1,{\bf z}_2)\\
&+x_4\bigl(\om({\bf z}_1,{\bf z}_5)+\om({\bf z}_2,{\bf z}_5)
-\om({\bf z}_3,{\bf z}_4)\bigr)\\
&+x_5\bigl(-\om({\bf z}_1,{\bf z}_4)+\om({\bf z}_2,{\bf z}_4)
+\om({\bf z}_3,{\bf z}_5)\bigr)+\om({\bf z}_4,{\bf z}_5).
\end{align*}

Thus $\phi^*(\om)\vert_P\equiv 0$ if and only if
\ea
\om({\bf z}_2,{\bf z}_3)=\om({\bf z}_1,{\bf z}_3)=
\om({\bf z}_1,{\bf z}_2)&=0,
\label{ev5eq4}\\
\om({\bf z}_1,{\bf z}_5)+\om({\bf z}_2,{\bf z}_5)
-\om({\bf z}_3,{\bf z}_4)&=0,
\label{ev5eq5}\\
-\om({\bf z}_1,{\bf z}_4)+\om({\bf z}_2,{\bf z}_4)
+\om({\bf z}_3,{\bf z}_5)&=0,
\label{ev5eq6}\\
\text{and}\qquad \om({\bf z}_4,{\bf z}_5)&=0.
\label{ev5eq7}
\ea
Now $\phi$ lies in the set ${\cal C}_P$ of Definition \ref{ev3def}
if and only if equations \eq{ev5eq4}--\eq{ev5eq7} hold and $\phi(P)$
is 2-dimensional, which is an open condition on $\phi$. Hence the
${\bf z}_j$ have 36 real parameters satisfying 6 real equations, 
so that ${\cal C}_P$ has dimension~30.

Motivated by \eq{ev3eq1}, define a `cross product' 
$\t:\C^3\t\C^3\ra\C^3$ by
\e
({\bf r}\t{\bf s})^b={\bf r}^{a_1}{\bf s}^{a_2}(\Re\Om)_{a_1a_2a_3}
g^{a_3b},
\label{ev5eq8}
\e
regarding $\C^3$ as a real vector space, and using the index notation 
for tensors on $\C^3$. Calculation shows that in complex coordinates, 
we have
\e
(r_1,r_2,r_3)\t(s_1,s_2,s_3)=\ha(\bar r_2\bar s_3-\bar r_3\bar s_2,
\bar r_3\bar s_1-\bar r_1\bar s_3,\bar r_1\bar s_2-\bar r_2\bar s_1),
\label{ev5eq9}
\e
so that `$\t$' is complex anti-bilinear. Note that this cross product
is equivariant under the action of $\SU(3)$ on~$\C^3$.

In \S\ref{ev3} we explained how to construct special Lagrangian $m$-folds
using an evolution equation \eq{ev3eq1} for $\phi\in{\cal C}_P$. We shall
write this equation out explicitly for $\phi$ of the form \eq{ev5eq3}.
Let ${\bf z}_1(t),\ldots,{\bf z}_6(t)$ be differentiable functions 
$\R\ra\C^3$, and define $\phi_t$ by \eq{ev5eq3} for $t\in\R$. Then, 
comparing equations \eq{ev3eq1}, \eq{ev5eq2} and \eq{ev5eq8}, we see 
that \eq{ev3eq1} holds for the family $\bigl\{\phi_t:t\in\R\bigr\}$ 
if and only if
\begin{align*}
\frac{\d\phi_t}{\d t}&(x_1,\ldots,x_5)
=2x_1{\bf z}_2\!\t {\bf z}_3+2x_2{\bf z}_1\!\t {\bf z}_3-2x_3
{\bf z}_1\!\t {\bf z}_2\\
&+x_4({\bf z}_1\!\t\!{\bf z}_5\!+\!{\bf z}_2\!\t\!{\bf z}_5\!
-\!{\bf z}_3\!\t\!{\bf z}_4)
+x_5(-{\bf z}_1\!\t\!{\bf z}_4\!+\!{\bf z}_2\!\t\!{\bf z}_4\!
+\!{\bf z}_3\!\t\!{\bf z}_5)+{\bf z}_4\t {\bf z}_5.
\end{align*}
Using \eq{ev5eq3} we get expressions for $\d{\bf z}_j/\d t$ for
$j=1,\ldots,6$. So applying Theorem \ref{ev3thm1}, we prove:

\begin{thm} Suppose ${\bf z}_1,\ldots,{\bf z}_6:\R\ra\C^3$ are 
differentiable functions satisfying equations \eq{ev5eq4}--\eq{ev5eq7}
at\/ $t=0$ and
\begin{gather}
\frac{\d{\bf z}_1}{\d t}=2{\bf z}_2\!\t {\bf z}_3,\qquad
\frac{\d{\bf z}_2}{\d t}=2{\bf z}_1\!\t {\bf z}_3,\qquad
\frac{\d{\bf z}_3}{\d t}=-2{\bf z}_1\!\t {\bf z}_2,
\label{ev5eq10}\\
\frac{\d{\bf z}_4}{\d t}=
{\bf z}_1\!\t\!{\bf z}_5\!+\!{\bf z}_2\!\t\!{\bf z}_5\!
-\!{\bf z}_3\!\t\!{\bf z}_4,\quad
\frac{\d{\bf z}_5}{\d t}=-{\bf z}_1\!\t\!{\bf z}_4\!+\!
{\bf z}_2\!\t\!{\bf z}_4\!+\!{\bf z}_3\!\t\!{\bf z}_5,
\label{ev5eq11}\\
\text{and}\qquad
\frac{\d{\bf z}_6}{\d t}={\bf z}_4\t {\bf z}_5
\label{ev5eq12}
\end{gather}
for all\/ $t\in\R$, where `$\t$\!' is as in \eq{ev5eq9}. Define
a subset\/ $N$ of\/ $\C^3$ by
\e
\begin{split}
N=\bigl\{&\ha(y_1^2+y_2^2)\,{\bf z}_1(t)+
\ha(y_1^2-y_2^2)\,{\bf z}_2(t)+y_1y_2\,{\bf z}_3(t)\\
&+y_1\,{\bf z}_4(t)+y_2\,{\bf z}_5(t)+{\bf z}_6(t):
y_1,y_2,t\in\R\bigr\}.
\end{split}
\label{ev5eq13}
\e
Then $N$ is a special Lagrangian $3$-fold in $\C^3$ wherever
it is nonsingular.
\label{ev5thm}
\end{thm}

The results of \cite[\S 3]{Joyc3} also show that if \eq{ev5eq4}--\eq{ev5eq7}
hold at $t=0$ then they hold for all $t\in\R$, and that given initial 
values ${\bf z}_1(0),\ldots,{\bf z}_6(0)$, there exist unique solutions 
${\bf z}_1(t),\ldots,{\bf z}_6(t)$ to \eq{ev5eq10}--\eq{ev5eq12} for $t$ 
in $(-\ep,\ep)$ and some small $\ep>0$. In fact it will follow from later 
results that solutions always exist for all $t\in\R$, and this is why we
have used $t\in\R$ rather than $t\in(-\ep,\ep)$ above.

\subsection{Transformation of ${\bf z}_1,\ldots,{\bf z}_6$ under 
$\GL(2,\R)\lt\R^2$}
\label{ev51}

The evolution data $(P,\chi)$ we used above was derived in 
\S\ref{ev4} from the action of $\SL(2,\R)\lt\R^2$ on $\R^2$ by
symplectic affine transformations. We shall now show that the
construction of Theorem \ref{ev5thm} is invariant under the 
action not just of $\SL(2,\R)\lt\R^2$, but under the full group 
$\GL(2,\R)\lt\R^2$ of affine transformations of $\R^2$. That is, 
we shall define an action of $\GL(2,\R)\lt\R^2$ on the set of 
solutions ${\bf z}_1,\ldots,{\bf z}_6$ of \eq{ev5eq10}--\eq{ev5eq12}
which fixes the corresponding SL 3-folds $N$ of~\eq{ev5eq13}.

Consider the affine transformation of $\R^2$ given by
\e
(y_1,y_2)\mapsto(ay_1+by_2+e,cy_1+dy_2+f),
\label{ev5eq14}
\e
where $a,b,c,d,e,f\in\R$, and the determinant $\de=ad-bc$ is
nonzero. Suppose that ${\bf z}_1,\ldots,{\bf z}_6:\R\ra\C^3$  
satisfy \eq{ev5eq10}--\eq{ev5eq12}. The natural way to make the 
transformation \eq{ev5eq14} act upon ${\bf z}_1,\ldots,{\bf z}_6$
is to define ${\bf z}_1',\ldots,{\bf z}_6'$ by equating coefficients
of polynomials in $y_1',y_2'$ in the equation
\begin{align*}
&\ha(y_1^2+y_2^2)\,{\bf z}_1+\ha(y_1^2-y_2^2)\,{\bf z}_2+
y_1y_2\,{\bf z}_3+y_1\,{\bf z}_4+y_2\,{\bf z}_5+{\bf z}_6=\\
&\ha\bigl((y_1')^2+(y_2')^2\bigr){\bf z}_1'
+\ha((y_1')^2-(y_2')^2\bigr){\bf z}_2'+
y_1'y_2'{\bf z}_3'+y_1'{\bf z}_4'+y_2'{\bf z}_5'+{\bf z}_6',
\end{align*}
where $y_1=ay_1'+by_2'+e$ and~$y_2=cy_1'+dy_2'+f$.

Here each side is a polynomial in $y_1',y_2'$ with values in $\C^3$,
taken from \eq{ev5eq13}. A straightforward calculation gives expressions
for ${\bf z}_1',\ldots,{\bf z}_6'$ in terms of ${\bf z}_1,\ldots,{\bf z}_6$
and $a,\ldots,f$, so that for example
\begin{equation*}
{\bf z}_1'(t)=\ha(a^2+b^2+c^2+d^2){\bf z}_1(t)
+\ha(a^2+b^2-c^2-d^2){\bf z}_2(t)+(ac+bd){\bf z}_3(t).
\end{equation*}
If the transformation \eq{ev5eq14} lies in $\SL(2,\R)\lt\R^2$ then it
is easy to see from the construction of the evolution data out of the
action of $\SL(2,\R)\lt\R^2$ that ${\bf z}_1',\ldots,{\bf z}_6'$ must
satisfy \eq{ev5eq10}--\eq{ev5eq12}, and yield exactly the same
SL 3-fold $N$ in \eq{ev5eq13} as ${\bf z}_1,\ldots,{\bf z}_6$ do.

However, if \eq{ev5eq14} lies in $\GL(2,\R)\lt\R^2$ rather than
$\SL(2,\R)\lt\R^2$ then the ${\bf z}_j'$ will not in general 
satisfy \eq{ev5eq10}--\eq{ev5eq12}. This is because the data 
$\chi$ of \eq{ev5eq2} is essentially the same as $\frac{\pd}{\pd y_1}
\w\frac{\pd}{\pd y_2}$, but \eq{ev5eq14} multiplies $\frac{\pd}{\pd y_1}
\w\frac{\pd}{\pd y_2}$ by $\de=ad-bc$. Thus, in \eq{ev5eq10}--\eq{ev5eq12}
the $\d{\bf z}_j/\d t$ are also multiplied by $\de$. We deal with this
by replacing $t$ by $t'=\de^{-1}t$, and then the ${\bf z}_j'$ satisfy
\eq{ev5eq10}--\eq{ev5eq12} with respect to the new time variable $t'$.
Hence we prove:

\begin{prop} Suppose that\/ ${\bf z}_1,\ldots,{\bf z}_6:\R\ra\C^3$
satisfy \eq{ev5eq10}--\eq{ev5eq12}. Let\/ $a,b,c,d,e,f\in\R$ with\/ 
$\de=ad-bc\ne 0$, and define ${\bf z}_1',\ldots,{\bf z}_6':\R\ra\C^3$ by
\begin{align}
\begin{split}
{\bf z}_1'(t)&=\ha(a^2\!+\!b^2\!+\!c^2\!+\!d^2){\bf z}_1(\de t)
+\ha(a^2\!+\!b^2\!-\!c^2\!-\!d^2){\bf z}_2(\de t) \\
&\qquad\qquad +(ac+bd){\bf z}_3(\de t),
\end{split}
\label{ev5eq15}\\
\begin{split}
{\bf z}_2'(t)&=\ha(a^2\!-\!b^2\!+\!c^2\!-\!d^2){\bf z}_1(\de t)
+\ha(a^2\!-\!b^2\!-\!c^2\!+\!d^2){\bf z}_2(\de t) \\
&\qquad\qquad +(ac-bd){\bf z}_3(\de t),
\end{split}
\label{ev5eq16}\\
{\bf z}_3'(t)&=(ab+cd){\bf z}_1(\de t)+(ab-cd){\bf z}_2(\de t)
+(ad+bc){\bf z}_3(\de t),
\label{ev5eq17}\\
\begin{split}
{\bf z}_4'(t)&=(ae+cf){\bf z}_1(\de t)
+(ae-cf){\bf z}_2(\de t)+(af+ce){\bf z}_3(\de t)\\
&\qquad\qquad +a\,{\bf z}_4(\de t)+c\,{\bf z}_5(\de t),
\end{split}
\label{ev5eq18}\\
\begin{split}
{\bf z}_5'(t)&=(be+df){\bf z}_1(\de t)
+(be-df){\bf z}_2(\de t)+(bf+de){\bf z}_3(\de t)\\
&\qquad\qquad +b\,{\bf z}_4(\de t)+d\,{\bf z}_5(\de t),
\end{split}
\label{ev5eq19}\\
\begin{split}
{\bf z}_6'(t)&=\ha(e^2+f^2){\bf z}_1(\de t)
+\ha(e^2-f^2){\bf z}_2(\de t)+ef{\bf z}_3(\de t)\\
&\qquad\qquad +e\,{\bf z}_4(\de t)+f\,{\bf z}_5(\de t)+{\bf z}_6(\de t).
\end{split}
\label{ev5eq20}
\end{align}
Then ${\bf z}_1',\ldots,{\bf z}_6'$ satisfy \eq{ev5eq10}--\eq{ev5eq12}. 
Furthermore, the ${\bf z}_j'$ satisfy \eq{ev5eq4}--\eq{ev5eq7} if and
only if the ${\bf z}_j$ do, and in this case the special Lagrangian
$3$-folds $N,N'$ constructed in \eq{ev5eq13} from the ${\bf z}_j$ 
and\/ ${\bf z}_j'$ are the same.
\label{ev5prop}
\end{prop}

Suppose we are given solutions ${\bf z}_1,{\bf z}_2,{\bf z}_3:\R\ra\C^3$ 
to \eq{ev5eq10}, and we wish to solve \eq{ev5eq11} for ${\bf z}_4$ and 
${\bf z}_5$. Now \eq{ev5eq11} is {\it linear} in ${\bf z}_4,{\bf z}_5$,
so one solution is ${\bf z}_4={\bf z}_5=0$. Apply the proposition with 
$a=d=1$ and $b=c=0$ and arbitrary values of $e,f$. It gives new solutions
${\bf z}_1',\ldots,{\bf z}_5'$ to \eq{ev5eq10} and \eq{ev5eq11}, where
\begin{equation*}
{\bf z}_1'={\bf z}_1,\;\>
{\bf z}_2'={\bf z}_2,\;\>
{\bf z}_3'={\bf z}_3,\;\>
{\bf z}_4'=e\,{\bf z}_1\!+\!e\,{\bf z}_2\!+\!f\,{\bf z}_3\;\>\text{and}\;\>
{\bf z}_5'=f\,{\bf z}_1\!-\!f\,{\bf z}_2\!+\!e\,{\bf z}_3.
\end{equation*}
This gives:

\begin{cor} Suppose ${\bf z}_1,{\bf z}_2,{\bf z}_3:\R\ra\C^3$ satisfy
\eq{ev5eq10}. Define
\begin{equation*}
{\bf z}_4=e\,{\bf z}_1+e\,{\bf z}_2+f\,{\bf z}_3\quad\text{and}\quad
{\bf z}_5=f\,{\bf z}_1-f\,{\bf z}_2+e\,{\bf z}_3,
\end{equation*}
for $e,f\in\R$. Then ${\bf z}_4,{\bf z}_5$ satisfy~\eq{ev5eq11}.
\label{ev5cor}
\end{cor}

This will be helpful later in solving \eq{ev5eq11}, given solutions
to~\eq{ev5eq10}.

\subsection{Discussion of the construction}
\label{ev52}

Here is a parameter count for family of the special Lagrangian 
3-folds $N$ in $\C^3$ constructed by the theorem. The initial data 
${\bf z}_1(0),\ldots,{\bf z}_6(0)$ has 36 real parameters, as each 
${\bf z}_j(0)$ lies in $\C^3$. These are subject to 6 real conditions 
\eq{ev5eq4}--\eq{ev5eq7}, reducing them to 30 real parameters. That 
is, $\dim{\cal C}_P=30$ in the notation of Definition \ref{ev3def}, 
so the family of curves in ${\cal C}_P$ has dimension~29.

However, we saw in \S\ref{ev51} that $\GL(2,\R)\lt\R^2$ acts
on this family of curves in ${\cal C}_P$, and two curves related
by the group action give the same 3-fold. As $\GL(2,\R)\lt\R^2$
has dimension 6, this means that the family of distinct SL 3-folds 
in $\C^3$ constructed above has dimension~$29-6=23$.

If we identify SL 3-folds isomorphic under automorphisms of $\C^3$, 
the dimension reduces still further. The appropriate automorphism 
group is $\SU(3)\lt\C^3$, with dimension 14. Thus the family of
distinct SL 3-folds in $\C^3$ up to automorphisms of $\C^3$ has 
dimension~$23-14=9$. 

So the number of interesting real parameters in the construction 
of Theorem \ref{ev5thm} is 9. For comparison, the number of 
interesting parameters in the construction of \cite[\S 6]{Joyc3} is 3,
and the number in the construction of \cite[\S 7]{Joyc3} with $m=3$ is 2. 
Hence the construction above is quite a lot more general than those
of \cite[\S 6]{Joyc3}, and \cite[\S 7]{Joyc3} with $m=3$. In fact
the $m=3$ cases of \cite[\S 7]{Joyc3}, discussed in
\cite[Ex.~7.4 \& Ex.~7.5]{Joyc3}, occur as special cases of the
construction above.

In \S\ref{ev6}--\S\ref{ev11} we will study the solutions of the 
o.d.e.s \eq{ev5eq10}--\eq{ev5eq12}, and so construct special 
Lagrangian 3-folds in $\C^3$. The way we have divided the equations 
up suggests a three-stage solution process. For \eq{ev5eq10} shows 
that $\d{\bf z}_1/\d t$, $\d{\bf z}_2/\d t$ and $\d{\bf z}_3/\d t$ 
depend only on ${\bf z}_1,{\bf z}_2,{\bf z}_3$, and not on 
${\bf z}_4,{\bf z}_5$ or ${\bf z}_6$. Thus in the first stage we 
solve the nonlinear equations \eq{ev5eq10} for ${\bf z}_1,{\bf z}_2,
{\bf z}_3$, ignoring ${\bf z}_4,{\bf z}_5$ and~${\bf z}_6$. 

Then in the second stage we regard ${\bf z}_1,{\bf z}_2,{\bf z}_3$ as 
fixed, and solve equations \eq{ev5eq11} for ${\bf z}_4$ and ${\bf z}_5$.
Notice that \eq{ev5eq11} are {\it linear} in ${\bf z}_4,{\bf z}_5$,
which makes them much easier to solve. Also, Corollary \ref{ev5cor}
gives us two of the six solutions automatically. Finally, in the 
third stage we regard ${\bf z}_1,\ldots,{\bf z}_5$ as fixed and solve 
\eq{ev5eq12} for ${\bf z}_6$, which is just a matter of integration.

Now the first stage reduces to a problem we have already studied in 
\cite{Joyc3}. By ignoring ${\bf z}_4,{\bf z}_5$ and ${\bf z}_6$, 
we are effectively considering maps $\phi:\R^3\ra\C^3$ given by 
$\phi:(x_1,x_2,x_3)\mapsto x_1{\bf z}_1+x_2{\bf z}_2+x_3{\bf z}_3$. 
Then $P$ in $\R^3$ is the set of $(x_1,x_2,x_3)$ of the form 
$\bigl(\ha(y_1^2+y_2^2),\ha(y_1^2-y_2^2),y_1y_2\bigr)$ for 
$y_1,y_2\in\R$. This satisfies~$x_1^2=x_2^2+x_3^2$. 

Thus, we are evolving the image of a quadric cone in $\R^3$ under 
a linear map $\R^3\ra\C^3$. This is exactly the what we did in 
\cite[\S 6]{Joyc3}. The equations \eq{ev5eq10} are in fact equivalent 
to the problem considered in \cite[\S 6]{Joyc3}, so we shall use the 
material of \cite[\S 6]{Joyc3} to understand their solutions.

Here is some notation that will be useful. Define a map
$\Phi:\R^3\ra\C^3$ by
\e
\begin{split}
\Phi(y_1,y_2,t)&=\ha(y_1^2+y_2^2)\,{\bf z}_1(t)+
\ha(y_1^2-y_2^2)\,{\bf z}_2(t)+y_1y_2\,{\bf z}_3(t)\\
&\qquad\qquad +y_1\,{\bf z}_4(t)+y_2\,{\bf z}_5(t)+{\bf z}_6(t).
\end{split}
\label{ev5eq21}
\e
Then the SL 3-fold $N$ of \eq{ev5eq13} is the image of $\Phi$, that
is,~$N=\bigl\{\Phi(y_1,y_2,t):y_1,y_2,t\in\R\bigr\}$.

If $\Phi$ is an immersion then $N$ is a nonsingular immersed
3-submanifold. The points where $\Phi$ is not an immersion generally
lead to {\it singularities} of $N$. We will study the points where
$\Phi$ is not an immersion in \S\ref{ev6}. In particular, we will
show that $\Phi$ is an immersion outside a set of real codimension
one in the family of all $\Phi$ generated in Theorem \ref{ev5thm}.
Thus, {\it generic} SL 3-folds $N$ from Theorem \ref{ev5thm} are
nonsingular as immersed 3-submanifolds.

Another question we shall be interested in is whether the maps
$\Phi$ are {\it periodic} in $t$. That is, we wish to know whether
there exists $T>0$ such that $\Phi(y_1,y_2,t+T)=\Phi(y_1,y_2,t)$
for all $y_1,y_2,t\in\R$. If this holds then we can regard $\Phi$
as mapping $\R^2\t{\cal S}^1\ra\C^3$ rather than $\R^3\ra\C^3$,
where ${\cal S}^1=\R/T\Z$, so that if $\Phi$ is an immersion then
$N$ is an immersed copy of $\R^2\t{\cal S}^1$ rather than~$\R^3$.

Periodic solutions are interesting they give us examples of SL
3-folds in $\C^3$ with different topologies, and because they are
often suitable local models for singularities of SL 3-folds in
Calabi--Yau 3-folds, whereas the non-periodic solutions usually
are not suitable because they are not closed in $\C^3$, or for
other reasons.

\section{Singularities of these SL 3-folds}
\label{ev6}

We shall now study the singularities of the special Lagrangian
3-folds constructed in Theorem \ref{ev5thm}. A good way to do this
is to use the map $\Phi:\R^3\ra\C^3$ defined in \eq{ev5eq21}. 
Clearly $\Phi$ is smooth. If at each $(y_1,y_2,t)\in\R^3$ its
derivative $\d\Phi\vert_{(y_1,y_2,t)}:\R^3\ra\C^3$ is injective
then $\Phi$ is an {\it immersion}, and $N=\Image\Phi$ is nonsingular
as an immersed 3-submanifold.

Thus the singularities of $N$ come from points $(y_1,y_2,t)$ for
which $\d\Phi\vert_{(y_1,y_2,t)}$ is not injective. Generically,
if $\d\Phi\vert_{(y_1,y_2,t)}$ is not injective then $N$ is
singular at $\Phi(y_1,y_2,t)$. But we will see in cases (a) and (b)
of \S\ref{ev83} that it can happen that $\Phi$ is not an immersion,
but $N$ is a subset of a nonsingular 3-fold, so that the apparent
singularity is due to badly chosen coordinates.

We begin with a couple of lemmas about $\Phi$. The first is true by
construction.

\begin{lem} The map $\Phi$ of\/ \eq{ev5eq21} satisfies
\e
\begin{gathered}
\om\Bigl(\frac{\pd\Phi}{\pd y_1},\frac{\pd\Phi}{\pd y_2}\Bigr)=
\om\Bigl(\frac{\pd\Phi}{\pd y_1},\frac{\pd\Phi}{\pd t}\Bigr)=
\om\Bigl(\frac{\pd\Phi}{\pd y_2},\frac{\pd\Phi}{\pd t}\Bigr)=0\\
\text{and}\qquad
\frac{\pd\Phi}{\pd y_1}\t\frac{\pd\Phi}{\pd y_2}=\frac{\pd\Phi}{\pd t},
\end{gathered}
\label{ev6eq1}
\e
where `$\t$\!' is defined in~\eq{ev5eq9}.
\label{ev6lem1}
\end{lem}

The second gives a simple criterion to decide whether $\Phi$ is
an immersion.

\begin{lem} The map $\Phi$ of\/ \eq{ev5eq21} is an immersion near
$(y_1,y_2,t)\in\R^3$ if and only if\/ $\frac{\pd\Phi}{\pd t}(y_1,y_2,t)\ne
0$, or equivalently if and only if\/ $\frac{\pd\Phi}{\pd y_1}(y_1,y_2,t)$
and\/ $\frac{\pd\Phi}{\pd y_2}(y_1,y_2,t)$ are linearly independent.
\label{ev6lem2}
\end{lem}

\begin{proof} Since $\om\bigl(\frac{\pd\Phi}{\pd y_1},
\frac{\pd\Phi}{\pd y_2}\bigr)=0$ by \eq{ev6eq1}, one can show from
$\frac{\pd\Phi}{\pd y_1}\t\frac{\pd\Phi}{\pd y_2}=\frac{\pd\Phi}{\pd t}$ that
$\frac{\pd\Phi}{\pd y_1}$, $\frac{\pd\Phi}{\pd y_2}$ and $\frac{\pd\Phi}{\pd t}$
are linearly independent if and only if $\frac{\pd\Phi}{\pd t}\ne 0$, or
equivalently if and only if $\frac{\pd\Phi}{\pd y_1}$ and
$\frac{\pd\Phi}{\pd y_2}$ are linearly independent. But $\Phi$ is an
immersion near $(y_1,y_2,t)$ if and only if $\frac{\pd\Phi}{\pd y_1}$,
$\frac{\pd\Phi}{\pd y_2}$ and $\frac{\pd\Phi}{\pd t}$ are linearly independent.
\end{proof}

As in \S\ref{ev52}, we can do a parameter count on the family of
{\it singular} SL 3-folds in $\C^3$ constructed above. Consider
first the condition that $\Phi$ should not be an immersion at
(0,0,0). By Lemma \ref{ev6lem2}, this happens if and only if
$\frac{\pd\Phi}{\pd y_1}(0,0,0)$ and $\frac{\pd\Phi}{\pd y_2}(0,0,0)$
are linearly dependent, that is, if ${\bf z}_4(0)$ and ${\bf z}_5(0)$
are linearly dependent.

The set of linearly dependent pairs ${\bf z}_4(0),{\bf z}_5(0)$ in
$\C^3$ has dimension 7. Thus the set of initial data ${\bf z}_1(0),
\ldots,{\bf z}_6(0)$ with ${\bf z}_4(0)$ and ${\bf z}_5(0)$ linearly
dependent has $24+7=31$ real parameters. These are subject to 6 real
conditions \eq{ev5eq4}--\eq{ev5eq7}. But one of these,
$\om({\bf z}_4,{\bf z}_5)=0$, holds automatically as ${\bf z}_4(0)$
and ${\bf z}_5(0)$ are linearly dependent.

So the set of initial data ${\bf z}_1(0),\ldots,{\bf z}_6(0)$
satisfying \eq{ev5eq4}--\eq{ev5eq7} with ${\bf z}_4(0)$ and ${\bf z}_5(0)$
linearly dependent has $31-5=26$ real parameters. For comparison, we saw in
\S\ref{ev52} that the set of initial data ${\bf z}_1(0),\ldots,{\bf z}_6(0)$
satisfying \eq{ev5eq4}--\eq{ev5eq7} has 30 real parameters. Therefore the
condition that $\Phi$ should not be an immersion at $(0,0,0)$ is of
real codimension 4. By symmetry, the condition for $\Phi$ not to be an
immersion at any given point in $\R^3$ is also of codimension~4.

Thus we expect the family of singular 3-folds $N$ arising from Theorem
\ref{ev5thm} to be of codimension $4-3=1$ in the family of all such
3-folds. So the family of distinct singular SL 3-folds in $\C^3$ from
Theorem \ref{ev5thm}, up to automorphisms of $\C^3$, should have dimension
8. In particular, for generic 3-folds $N$ arising from Theorem \ref{ev5thm},
$\Phi$ is an immersion, and $N$ is a nonsingular immersed 3-submanifold.

Next we will describe the singularities of the 3-folds $N$ of Theorem
\ref{ev5thm} fairly explicitly, by modelling $N$ near a singular point.
A good way to do this is to expand $\Phi$ as a power series about the
singular point, to low order. For simplicity, we take the singular
point to be at (0,0,0) in $\R^3$ and~$\C^3$.

So let ${\bf z}_1,\ldots,{\bf z}_6$ and $\Phi$ satisfy all the
conditions of \S\ref{ev5}, and suppose that $\Phi(0,0,0)={\bf z}_6(0)=0$
and that $\d\Phi\vert_{(0,0,0)}$ is not injective. As above, this holds
if and only if ${\bf z}_4(0)$ and ${\bf z}_5(0)$ are linearly dependent.
We will expand $\Phi$ as a power series about 0 up to second order, and
use this to describe $N$ near its singular point~0.

Now in \S\ref{ev51} we described an action of $\GL(2,\R)\lt\R^2$
on the set of maps $\Phi$ satisfying the conditions of \S\ref{ev5},
that acts trivially on the corresponding SL 3-folds $N$. Since we
are really interested in $N$ rather than $\Phi$, we shall use this
action to put $\Phi$ in a more convenient form. Under the rotation
$\bigl(\begin{smallmatrix}a & b\\ c & d\end{smallmatrix}\bigr)
=\bigl(\begin{smallmatrix}\pha{-}\cos\th & \sin\th\\ -\sin\th &\cos\th
\end{smallmatrix}\bigr)$, with $e=f=0$, we see from \eq{ev5eq18} and
\eq{ev5eq19} that ${\bf z}_4(0)$ and ${\bf z}_5(0)$ are transformed~to
\begin{equation*}
{\bf z}_4'(0)=\cos\th\,{\bf z}_4(0)-\sin\th\,{\bf z}_5(0)
\quad\text{and}\quad
{\bf z}_5'(0)=\sin\th\,{\bf z}_4(0)+\cos\th\,{\bf z}_5(0).
\end{equation*}
As ${\bf z}_4(0)$ and ${\bf z}_5(0)$ are linearly dependent,
we may choose $\th$ so that~${\bf z}_5'(0)=0$.

So suppose that ${\bf z}_5(0)=0$. Take the initial data to be
\begin{equation*}
{\bf z}_1(0)={\bf v}+{\bf w},\;\>
{\bf z}_2(0)={\bf v}-{\bf w},\;\>
{\bf z}_3(0)={\bf x},\;\>
{\bf z}_4(0)={\bf u},\;\>
{\bf z}_5(0)={\bf z}_6(0)=0,
\end{equation*}
for vectors ${\bf u},{\bf v},{\bf w},{\bf x}$ in $\C^3$. Equations
\eq{ev5eq4}--\eq{ev5eq7} then reduce to
\e
\om({\bf u},{\bf w})=
\om({\bf u},{\bf x})=
\om({\bf v},{\bf w})=
\om({\bf v},{\bf x})=
\om({\bf w},{\bf x})=0.
\label{ev6eq2}
\e
Expanding ${\bf z}_1,\ldots,{\bf z}_6$ to low order in $t$ using
equations \eq{ev5eq10}--\eq{ev5eq12}, we find
\begin{alignat*}{2}
{\bf z}_1(t)&={\bf v}+{\bf w}+O(t),&\quad
{\bf z}_2(t)&={\bf v}-{\bf w}+O(t),\\
{\bf z}_3(t)&={\bf x}+O(t),&\quad
{\bf z}_4(t)&={\bf u}+t\,{\bf u}\t{\bf x}+O(t^2),\\
{\bf z}_5(t)&=2t\,{\bf u}\t{\bf w}+O(t^2),&\quad
{\bf z}_6(t)&=t^2\,{\bf u}\t({\bf u}\t{\bf w})+O(t^3),
\end{alignat*}
for small~$t$.

Now calculating using \eq{ev5eq9} shows that
\e
\begin{split}
{\bf u}\t({\bf u}\t{\bf w})&={\ts\frac{1}{4}}\bigl(g({\bf u},{\bf w})
+i\om({\bf u},{\bf w})\bigr){\bf u}-{\ts\frac{1}{4}}{\bf w}\\
&={\ts\frac{1}{4}}g({\bf u},{\bf w}){\bf u}-{\ts\frac{1}{4}}{\bf w},
\end{split}
\label{ev6eq3}
\e
as $\om({\bf u},{\bf w})=0$ by \eq{ev6eq2}. Substituting the above
expressions for ${\bf z}_j(t)$ into \eq{ev5eq21} and using \eq{ev6eq3}
to rewrite the ${\bf u}\t({\bf u}\t{\bf w})$ term in ${\bf z}_6(t)$,
we find that
\e
\begin{split}
\Phi(y_1,y_2,t)=&\bigl(y_1+{\ts\frac{1}{4}}
g({\bf u},{\bf w})t^2\bigr)\,{\bf u}+y_1^2\,{\bf v}
+\bigl(y_2^2-{\ts\frac{1}{4}}\ms{{\bf u}}t^2\bigr)\,{\bf w}\\
&+y_1y_2\,{\bf x}+y_1t\,{\bf u}\t{\bf x}+2y_2t\,{\bf u}\t{\bf w}\\
&+\text{third-order terms and above in $y_1,y_2,t$.}
\end{split}
\label{ev6eq4}
\e
This is the expansion of $\Phi$ up to second order in $y_1,y_2,t$.
Notice that the only first-order term is $y_1\,{\bf u}$.

Stratifying terms by their order in $y_1,y_2,t$ is probably not 
the most helpful view. Instead, it may be better to regard $y_1$ as
having twice the order of $y_2$ and $t$. To see this, observe from
\eq{ev6eq4} that
\e
\begin{split}
\Phi(\ep^2y_1,\ep y_2,\ep t)=
\ep^2\Bigl(&\bigl(
y_1+{\ts\frac{1}{4}}g({\bf u},{\bf w})t^2\bigr)\,{\bf u}
+\bigl(y_2^2-{\ts\frac{1}{4}}\ms{{\bf u}}t^2\bigr)\,{\bf w}\\
&+2y_2t\,{\bf u}\t{\bf w}\Bigr)+O(\ep^3),
\end{split}
\label{ev6eq5}
\e
for small $\ep$. That is, the dominant terms are those in
$y_1,y_2^2,ty_2$ and~$t^2$.

Suppose $\bf u$ and $\bf w$ are linearly independent, which is true
in the generic case. Then as $\om({\bf u},{\bf w})=0$ we see that
${\bf u},{\bf w}$ and ${\bf u}\t{\bf w}$ are linearly independent
and span a special Lagrangian $\R^3$ in $\C^3$. So \eq{ev6eq5}
shows that near 0 to lowest order, $N$ is the image of a map
$\R^3\ra\R^3$ given by
\begin{equation*}
(y_1,y_2,t)\longmapsto\bigl(y_1+{\ts\frac{1}{4}}g({\bf u},{\bf w})t^2,
y_2^2-{\ts\frac{1}{4}}\ms{{\bf u}}t^2,2y_2t\bigr),
\end{equation*}
with image in a special Lagrangian 3-plane. This map is a
{\it double cover} of $\R^3$, branched over the $x$-axis.

Therefore, if ${\bf u},{\bf w}$ are linearly independent then
near zero to lowest order, $N$ looks like a branched double cover
of a special Lagrangian 3-plane, branched over a real line. To
understand how $N$ deviates from this SL 3-plane to leading order,
we would need to know the terms in $\ep^3$ in \eq{ev6eq5}, which
come from the terms in $y_1y_2,y_1t,y_2^2t,y_2t^2$ and $t^3$ in
the expansion of $\Phi$. These can be calculated by the method
above, but we will not write them down as they are rather complex.

Most of the singularities of SL $m$-folds that we have met so far in
\cite{Joyc1}, \cite{Joyc2} and \cite{Joyc3} have been {\it conical}
singularities; that is, to lowest order the submanifold resembles a
cone with an {\it isolated} singularity at 0. These singularities
follow a different pattern. To lowest order they do resemble a
(degenerate) cone, the double SL 3-plane, but this cone should be 
regarded as {\it singular} along the branch locus, so that the 
singularity is not isolated. The interesting information about the 
singularity of $N$ is not just the cone itself, but also how $N$ 
deviates from the cone to leading order.

We assumed above that $\bf u$ and $\bf w$ are linearly independent.
Here is what happens if they are not.
\begin{itemize}
\item If ${\bf u}=0$ then ${\bf z}_4\equiv{\bf z}_5\equiv{\bf z}_6\equiv 0$,
and $\Phi$ reduces to
\begin{equation*}
\Phi(y_1,y_2,t)=\ha(y_1^2+y_2^2)\,{\bf z}_1(t)+
\ha(y_1^2-y_2^2)\,{\bf z}_2(t)+y_1y_2\,{\bf z}_3(t).
\end{equation*}
The image $N=\Image\Phi$ is a {\it cone} with singularity at 0, and
may be written
\begin{equation*}
N=\bigl\{x_1\,{\bf z}_1(t)+x_2\,{\bf z}_2(t)+x_3\,{\bf z}_3(t):
x_1,x_2,x_3,t\in\R,\quad x_1^2=x_2^2+x_3^2\bigr\}.
\end{equation*}
Cones of this kind were studied in \cite{Joyc3}, in particular
in~\cite[\S 6]{Joyc3}.

\item If ${\bf u}\ne 0$ and ${\bf w}=0$, we find that
${\bf z}_1\equiv{\bf z}_2$ for all $t$, that ${\bf z}_3$ is constant,
and that ${\bf z}_5\equiv{\bf z}_6\equiv 0$. This situation will be
studied in case (a) of \S\ref{ev83}, and an explicit expression for
$N$ is given in \eq{ev8eq9}. It turns out that $N$ is a subset of a
nonsingular product SL 3-fold $\Si\t\R$ in $\C^2\t\C$, and the
`singularity' is due to a poor choice of coordinates.

\item If ${\bf u},{\bf w}$ are nonzero but linearly dependent, we
need to expand $\Phi$ up to third order in $y_1,y_2,t$ and perform
a similar analysis to the above. It turns out that if ${\bf u},{\bf x}$
are linearly independent, then to lowest order $N$ flattens itself onto
the special Lagrangian 3-plane $\an{{\bf u},{\bf x},{\bf u}\t
{\bf x}}_{\sst\mathbb R}$ near zero, and is triply branched over the
real line $\an{{\bf u}}_{\sst\mathbb R}$. The details are complicated.
\end{itemize}

We conclude with two further remarks. Firstly, when ${\bf z}_5(0)=0$ as
above we have $\Phi(0,y_2,0)=\Phi(0,-y_2,0)$ for all $y_2\in\R$. This
means that away from the singularity $\Phi(0,0,0)$, the SL 3-fold $N$
intersects itself in a half-line, and so is not embedded. This applies
to every singular $N$ coming from Theorem \ref{ev5thm}, as any such $N$
can be transformed under the $\GL(2,\R)\lt\R^2$-action and translation
in $t$ to have~${\bf z}_5(0)=0$.

Secondly, one could extend the ideas above to study a certain class of
singularities of SL 3-folds in a very explicit way. Consider SL 3-folds
$N$ that are the images of real analytic maps $\Phi:\R^3\ra\C^3$, defined
near 0 in $\R^3$, and such that $\d\Phi\vert_{(0,0,0)}$ is not injective.
Generically $N$ is singular at $\Phi(0,0,0)$. By expanding $\Phi$ as a
power series about 0 and considering only low order terms, one could
hopefully construct local models for singularities of SL 3-folds.

\section{Division into cases using ${\bf z}_1,{\bf z}_2,{\bf z}_3$}
\label{ev7}

In order to study and understand the family of SL 3-folds $N$ 
constructed in Theorem \ref{ev5thm}, we shall find it helpful to
divide the solutions into cases according to the behaviour of the 
solutions ${\bf z}_1,{\bf z}_2,{\bf z}_3:\R\ra\C^3$ of \eq{ev5eq4} 
and \eq{ev5eq10}. For $t\in\R$, consider the real vector subspace 
$\an{{\bf z}_1(t),{\bf z}_2(t),{\bf z}_3(t)}_{\sst\mathbb R}$ in 
$\C^3$, where $\an{\ldots}_{\sst\mathbb R}$ is the span over~$\R$. 

It turns out that one of the most important things determining the 
behaviour of solutions ${\bf z}_1,{\bf z}_2,{\bf z}_3$ of \eq{ev5eq10} 
is the dimension of this vector space for generic $t$. Clearly the 
dimension lies between 0 and 3, so we divide into four cases:
\begin{itemize}
\item[(i)] $\dim\an{{\bf z}_1(t),{\bf z}_2(t),
{\bf z}_3(t)}_{\sst\mathbb R}=0$ for generic~$t\in\R$,
\item[(ii)] $\dim\an{{\bf z}_1(t),{\bf z}_2(t),
{\bf z}_3(t)}_{\sst\mathbb R}=1$ for generic~$t\in\R$,
\item[(iii)] $\dim\an{{\bf z}_1(t),{\bf z}_2(t),
{\bf z}_3(t)}_{\sst\mathbb R}=2$ for generic $t\in\R$, and
\item[(iv)] $\dim\an{{\bf z}_1(t),{\bf z}_2(t),
{\bf z}_3(t)}_{\sst\mathbb R}=3$ for generic~$t\in\R$.
\end{itemize}

It can be shown that in cases (i)--(iii) the dimension is constant
for all $t$. But in case (iv), the dimension can drop to 2 for
isolated $t\in\R$. The complexity of the solutions to 
\eq{ev5eq10}--\eq{ev5eq12} increases with the dimension. Thus cases 
(i) and (ii) are very straightforward, and we shall discuss them 
in the rest of the section. Case (iii) will be the subject of 
\S\ref{ev8}. Case (iv) is the most complicated, and will be divided 
into subcases and discussed in~\S\ref{ev9}--\S\ref{ev11}.
\medskip

\noindent{\bf Case (i).} 
\smallskip

\noindent In this case we have ${\bf z}_1\equiv{\bf z}_2\equiv
{\bf z}_3\equiv 0$ for all $t$. Thus ${\bf z}_4$ and ${\bf z}_5$ 
are constant by \eq{ev5eq11}, and \eq{ev5eq12} integrates to 
${\bf z}_6(t)=t\,{\bf z}_4\t {\bf z}_5+{\bf z}_6(0)$. So solutions 
exist for all $t\in\R$, and the SL 3-fold $N$ of \eq{ev5eq13} is
\begin{equation*}
N=\bigl\{y_1\,{\bf z}_4+y_2\,{\bf z}_5+
t\,{\bf z}_4\t {\bf z}_5+{\bf z}_6(0):y_1,y_2,t\in\R\bigr\}.
\end{equation*}
If ${\bf z}_4,{\bf z}_5$ are linearly independent and 
$\om({\bf z}_4,{\bf z}_5)=0$, as in \eq{ev5eq7}, this is an affine 
special Lagrangian $\R^3$ in $\C^3$. If ${\bf z}_4,{\bf z}_5$ are 
not linearly independent, it is a line or a point.
\medskip

\noindent{\bf Case (ii).} 
\smallskip

\noindent In this case ${\bf z}_1,{\bf z}_2,{\bf z}_3$ are proportional, 
so ${\bf z}_2\t{\bf z}_3\equiv{\bf z}_1\t{\bf z}_3\equiv{\bf z}_1\t
{\bf z}_2\equiv 0$, and thus ${\bf z}_1,{\bf z}_2,{\bf z}_3$ are 
constant by \eq{ev5eq10}. Let ${\bf z}_j=c_j{\bf v}$ for $j=1,2,3$, 
where $c_j\in\R$ and ${\bf v}\in\C^3$ is a unit vector. Define a 
quadratic polynomial $Q(y_1,y_2)$ by
\begin{equation*}
Q(y_1,y_2)=\ha c_1(y_1^2+y_2^2)+\ha c_2(y_1^2-y_2^2)+c_3y_1y_2.
\end{equation*}
Then from \eq{ev5eq13} we have 
\begin{equation*}
N=\bigl\{y_1\,{\bf z}_4(t)+y_2\,{\bf z}_5(t)+Q(y_1,y_2){\bf v}+
{\bf z}_6(t):y_1,y_2,t\in\R\bigr\}.
\end{equation*}

Define a quadric $P'$ in $\R^3$ to be $\bigl\{(y_1,y_2,y_3)\in\R^3:
y_3=Q(y_1,y_2)\bigr\}$, and for each $t\in\R$ define an affine map 
$\phi_t':\R^3\ra\C^3$ by
\begin{equation*}
\phi_t':(y_1,y_2,y_3)\mapsto y_1\,{\bf z}_4(t)+y_2\,{\bf z}_5(t)+
y_3\,{\bf v}+{\bf z}_6(t).
\end{equation*}
Then $N$ may be written $\bigl\{\phi_t'(p'):t\in\R$, 
$p'\in P'\bigr\}$. That is, $N$ is the total space of a family of 
quadrics $\phi_t'(P')$, which live in affine Lagrangian 3-planes
$\R^3$ in~$\C^3$.

Now in \cite{Joyc3} we constructed SL $m$-folds with this property, and
it is easy to show that the SL 3-folds of case (ii) above were studied
in \cite[\S 7]{Joyc3}, in particular in \cite[Ex.~7.4 \& Ex.~7.5]{Joyc3}.
By the classification of quadratic forms, as $Q$ is nonzero it is
equivalent under a linear transformation of $\R^2$ to one of the
standard forms
\begin{equation*}
y_1^2+y_2^2,\quad -y_1^2-y_2^2,\quad 
y_1^2-y_2^2\quad\text{and}\quad y_1^2.
\end{equation*}
In the first two cases the 3-fold $N$ of \eq{ev5eq13} is one of
those constructed in \cite[Ex.~7.4]{Joyc3}, in the third it is one 
of those constructed in \cite[Ex.~7.5]{Joyc3}, and in the fourth
$N$ is reducible and splits as a product $N'\t\R$ in $\C^2\t\C$, 
where $N'$ is an SL 2-fold in~$\C^2$.

\section{Case (iii) of \S\ref{ev7}}
\label{ev8}

We shall now study the SL 3-folds in $\C^3$ constructed in
Theorem \ref{ev5thm} corresponding to case (iii) of \S\ref{ev7}.
In \S\ref{ev81} we show that any ${\bf z}_1,{\bf z}_2,{\bf z}_3$ 
satisfying \eq{ev5eq4}, \eq{ev5eq10} and case (iii) of \S\ref{ev7} 
are equivalent under the natural symmetry group of the construction 
to ${\bf z}_1'={\bf z}_2'=({\rm e}^{it},-i{\rm e}^{-it},0)$
and~${\bf z}_3'=(0,0,1)$.

Then in \S\ref{ev82} we solve the remaining equations
\eq{ev5eq5}--\eq{ev5eq7} and \eq{ev5eq11}--\eq{ev5eq12} for
${\bf z}_4,{\bf z}_5$ and ${\bf z}_6$, with these values of
${\bf z}_1,{\bf z}_2,{\bf z}_3$. This leads to the main result
of this section, Theorem \ref{ev8thm2}, which gives an explicit
family of SL 3-folds in $\C^3$. Section \ref{ev83} then discusses
these 3-folds, doing a parameter count, showing that they are all
ruled by straight lines, and studying the periodic solutions.

\subsection{Solving the equations for ${\bf z}_1,{\bf z}_2$ and 
${\bf z}_3$}
\label{ev81}

As in case (iii) of \S\ref{ev7}, suppose that ${\bf z}_1,{\bf z}_2,
{\bf z}_3:\R\ra\C^3$ are solutions of \eq{ev5eq4} and \eq{ev5eq10} 
such that $\an{{\bf z}_1(t),{\bf z}_2(t),{\bf z}_3(t)}_{\sst\mathbb R}$
has dimension 2 for generic $t\in\R$. This dimension must be less than 
or equal to 2 for all $t\in\R$, by upper semicontinuity of dimension. 
But if the dimension is 0 or 1 for any $t$ then it is so for all, since 
${\bf z}_1,{\bf z}_2,{\bf z}_3$ are constant as in parts (i) and (ii) 
of \S\ref{ev7}. Thus $\dim\an{{\bf z}_1(t),{\bf z}_2(t),
{\bf z}_3(t)}_{\sst\mathbb R}=2$ for {\it all}\/~$t\in\R$.

Now in \S\ref{ev51} we defined an action of $\GL(2,\R)\lt\R^2$ 
on the set of solutions ${\bf z}_1,\ldots,{\bf z}_6$ to 
\eq{ev5eq10}--\eq{ev5eq12}. If we restrict our attention to
${\bf z}_1,{\bf z}_2,{\bf z}_3$ then $e$ and $f$ play no 
r\^ole, and the group acting is $\GL(2,\R)$. We shall use 
this $\GL(2,\R)$-action to write ${\bf z}_1,{\bf z}_2,{\bf z}_3$
in a convenient form.

\begin{prop} Let\/ ${\bf z}_1,{\bf z}_2,{\bf z}_3$ be as above. 
Then they are equivalent under the $\GL(2,\R)$-action described 
in Proposition \ref{ev5prop} to ${\bf z}_1',{\bf z}_2',{\bf z}_3':
\R\ra\C^3$ with\/ ${\bf z}_1'(0)={\bf z}_2'(0)$ 
and\/~$\bmd{{\bf z}_3'(0)}=1$.
\label{ev8prop1}
\end{prop}

\begin{proof} As ${\bf z}_1(0),{\bf z}_2(0),{\bf z}_3(0)$ span
a vector space of dimension 2, they satisfy $a_1{\bf z}_1(0)+
a_2{\bf z}_2(0)+a_3{\bf z}_3(0)=0$ for some $a_1,a_2,a_3\in\R$ not 
all zero. We may think of the ${\bf z}_j(0)$ as giving a linear map 
$S^2\R^2\ra\C^3$ with kernel $\R$, and so $(a_1,a_2,a_3)$ defines 
a point in the kernel in $S^2\R^2$. Under the action of $\GL(2,\R)$ 
on ${\bf z}_1,{\bf z}_2,{\bf z}_3$ defined in \eq{ev5eq15}--\eq{ev5eq17}, 
the vector $(a_1,a_2,a_3)$ transforms under the natural action of 
$\GL(2,\R)$ on~$\R^3=S^2\R^2$.

It can be shown that $\GL(2,\R)$ acts on $\R^3$ as the group 
preserving the {\it Lorentzian conformal structure} with null cone 
$a_1^2-a_2^2-a_3^2=0$. Thus there are three $\GL(2,\R)$-orbits of
nonzero vectors in $\R^3$: the `time-like' vectors with
$a_1^2-a_2^2-a_3^2>0$, the `space-like' vectors with 
$a_1^2-a_2^2-a_3^2<0$, and the `null' vectors 
with~$a_1^2-a_2^2-a_3^2=0$.

Therefore every nonzero vector $(a_1,a_2,a_3)$ is equivalent
under the $\GL(2,\R)$-action to $(1,0,0)$ or $(0,1,0)$ or
$(1,-1,0)$. Hence we can transform ${\bf z}_1,{\bf z}_2,{\bf z}_3$
under the $\GL(2,\R)$ action to ${\bf z}_1',{\bf z}_2',{\bf z}_3'$
satisfying either (a) ${\bf z}_1'(0)=0$, (b) ${\bf z}_2'(0)=0$, or
(c) ${\bf z}_1'(0)-{\bf z}_2'(0)=0$. We shall show that (a) and (b) 
contradict case (iii) of \S\ref{ev7}, which leaves case~(c).

In part (a), as ${\bf z}_2'(0),{\bf z}_3'(0)$ are linearly 
independent and $\om\bigl({\bf z}_2'(0),{\bf z}_3'(0)\bigr)=0$ 
by \eq{ev5eq4}, we see that ${\bf z}_2'(0),{\bf z}_3'(0)$ and
${\bf z}_2'(0)\t{\bf z}_3'(0)$ are linearly independent. But
\eq{ev5eq10} gives
\begin{equation*}
{\bf z}_1'(t)=2t\,{\bf z}_2'(0)\t{\bf z}_3'(0)\!+\!O(t^3),\;\>
{\bf z}_2'(t)={\bf z}_2'(0)\!+\!O(t^2),\;\>
{\bf z}_3'(t)={\bf z}_3'(0)\!+\!O(t^2)
\end{equation*}
for small $t$. Therefore ${\bf z}_1'(t),{\bf z}_2'(t)$ and 
${\bf z}_3'(t)$ are linearly independent for small, nonzero $t$.
This contradicts (iii). Similarly, part (b) leads to a
contradiction.

Thus part (c) holds, and we can transform the ${\bf z}_j$ 
to ${\bf z}_j'$ with ${\bf z}_1'(0)={\bf z}_2'(0)$. Clearly 
${\bf z}_3'(0)$ must be nonzero for case (iii) to hold.
We can then use a dilation in $\GL(2,\R)$ to rescale the
${\bf z}_j'$ to get $\bmd{{\bf z}_3'(0)}=1$. This completes
the proof.
\end{proof}

Apart from the $\GL(2,\R)$-action considered above, the
construction is also invariant under the action of $\SU(3)$
on $\C^3$. Given that $\bmd{{\bf z}_3'(0)}=1$, we can apply
an $\SU(3)$ transformation to fix ${\bf z}_3'(0)=(0,0,1)$.
Thus we see that any solution ${\bf z}_1,{\bf z}_2,{\bf z}_3$
of case (iii) can be transformed under the natural symmetry
groups $\GL(2,\R)$, $\SU(3)$ of the construction to a new
solution ${\bf z}_1',{\bf z}_2',{\bf z}_3'$ with 
${\bf z}_1'(0)={\bf z}_2'(0)$ and~${\bf z}_3'(0)=(0,0,1)$.

We can now solve for ${\bf z}_1',{\bf z}_2'$ and ${\bf z}_3'$ 
explicitly. By symmetry between ${\bf z}_1$ and ${\bf z}_2$ in the 
first two equations of \eq{ev5eq10} we see that ${\bf z}_1'(t)=
{\bf z}_2'(t)$ for all $t$. Also, $\d{\bf z}_3'/\d t=0$ by the third 
equation of \eq{ev5eq10}, so that ${\bf z}_3'(t)=(0,0,1)$ for all~$t$.

Set ${\bf z}_1'={\bf z}_2'=(w_1,w_2,w_3)$, for differentiable functions 
$w_1,w_2,w_3:\R\ra\C$. Then equations \eq{ev5eq10} become
\begin{equation*}
\frac{\d}{\d t}(w_1,w_2,w_3)=2(w_1,w_2,w_3)\t(0,0,1)=(\bar w_2,-\bar w_1,0),
\end{equation*}
so that
\begin{equation*}
\frac{\d w_1}{\d t}=\bar w_2,\quad
\frac{\d w_2}{\d t}=-\bar w_1\quad\text{and}\quad
\frac{\d w_3}{\d t}=0.
\end{equation*}
The first two equations give $\frac{\d^2w_1}{\d t^2}=-w_1$, which has
solutions $w_1(t)=X'{\rm e}^{it}+Y'{\rm e}^{-it}$ for $X',Y'\in\C$,
and so $w_2(t)=i\bar Y'{\rm e}^{it}-i\bar X'{\rm e}^{-it}$. The
third equation has solutions $w_3(t)=Z'$ for $Z'\in\C$. But we need
the solutions to satisfy \eq{ev5eq4}, which reduces to~$\Im Z'=0$.

Thus ${\bf z}_1',{\bf z}_2'$ and ${\bf z}_3'$ are given by
\begin{equation*}
{\bf z}_1'(t)={\bf z}_2'(t)=(X'{\rm e}^{it}+Y'{\rm e}^{-it},
i\bar Y'{\rm e}^{it}-i\bar X'{\rm e}^{-it},Z')
\quad\text{and}\quad {\bf z}_3'(t)=(0,0,1),
\end{equation*}
with $X',Y'\in\C$ and $Z'\in\R$. The condition that 
$\dim\an{{\bf z}_1(t),{\bf z}_2(t),{\bf z}_3(t)}_{\sst\mathbb R}=2$
is satisfied, for all $t$, if $X'$ and $Y'$ are not both zero.

Although we have used both the $\GL(2,\R)$ and $\SU(3)$ actions to
bring the ${\bf z}_j'$ to this special form, we have not used all
the freedom in these two group actions, and we can use what remains to
choose the values of $X',Y'$ and $Z'$. From \eq{ev5eq15}--\eq{ev5eq17},
the subgroup of $\bigl(\begin{smallmatrix}a & b\\ c &d\end{smallmatrix}
\bigr)$ in $\GL(2,\R)$ preserving the condition ${\bf z}_1={\bf z}_2$
and fixing ${\bf z}_3$ are those with $b=0$ and $ad=1$, so that~$\de=1$.

Thus, putting $b=0$ and $d=a^{-1}$, by \eq{ev5eq15}--\eq{ev5eq17}
we transform ${\bf z}_1',{\bf z}_2',{\bf z}_3'$ to
\begin{equation*}
{\bf z}_1''(t)={\bf z}_2''(t)=a^2{\bf z}_1'(t)+ac(0,0,1)
\quad\text{and}\quad {\bf z}_3''(t)=(0,0,1).
\end{equation*}
Choose $a=(\ms{X'}+\ms{Y'})^{-1/2}$ and $c=-a^{-1}Z'$. Then we get
\begin{equation*}
{\bf z}_1''(t)={\bf z}_2''(t)=(X''{\rm e}^{it}+Y''{\rm e}^{-it},
i\bar Y''{\rm e}^{it}-i\bar X''{\rm e}^{-it},0)
\quad\text{and}\quad {\bf z}_3''(t)=(0,0,1),
\end{equation*}
where $\ms{X''}+\ms{Y''}=1$. Finally, applying the $\SU(3)$ matrix
\begin{equation*}
\begin{pmatrix}\bar X'' & -iY'' & 0 \\ -i\bar Y'' & X''& 0 \\ 
0 & 0 & 1 \end{pmatrix}
\end{equation*}
we transform the ${\bf z}_j''$ to ${\bf z}_j'''$, where
\begin{equation*}
{\bf z}_1'''(t)={\bf z}_2'''(t)=({\rm e}^{it},-i{\rm e}^{-it},0)
\quad\text{and}\quad {\bf z}_3'''(t)=(0,0,1).
\end{equation*}
Hence we have proved:

\begin{thm} Any ${\bf z}_1,{\bf z}_2,{\bf z}_3:\R\ra\C^3$ satisfying
\eq{ev5eq4}, \eq{ev5eq10} and part\/ {\rm(iii)} of\/ \S\ref{ev7} are 
equivalent under the natural actions of\/ $\GL(2,\R)$ and\/ $\SU(3)$~to
\e
\hat{\bf z}_1(t)=\hat{\bf z}_2(t)=({\rm e}^{it},-i{\rm e}^{-it},0)
\quad\text{and}\quad \hat{\bf z}_3(t)=(0,0,1).
\label{ev8eq1}
\e
\label{ev8thm1}
\end{thm}

\subsection{Solving the equations for ${\bf z}_4,{\bf z}_5$ and ${\bf z}_6$}
\label{ev82}

Now let us fix ${\bf z}_1(t)={\bf z}_2(t)=({\rm e}^{it},-i{\rm e}^{-it},0)$
and ${\bf z}_3(t)=(0,0,1)$, as in Theorem \ref{ev8thm1}, and solve
\eq{ev5eq11} and \eq{ev5eq12} for ${\bf z}_4,{\bf z}_5$ and ${\bf z}_6$.
Write ${\bf z}_4=(p_1,p_2,p_3)$ and ${\bf z}_5=(q_1,q_2,q_3)$ for
differentiable functions $p_j,q_j:\R\ra\C$. Then the first and second
equations of \eq{ev5eq11} become
\begin{gather}
\begin{gathered}
\frac{\d p_1}{\d t}=i{\rm e}^{it}\bar q_3+\ha\bar p_2,\qquad
\frac{\d p_2}{\d t}=-{\rm e}^{-it}\bar q_3-\ha\bar p_1\\
\text{and}\qquad
\frac{\d p_3}{\d t}=-i{\rm e}^{it}\bar q_1+{\rm e}^{-it}\bar q_2,
\end{gathered}
\label{ev8eq2}\\
\frac{\d q_1}{\d t}=-\ha\bar q_2,\quad
\frac{\d q_2}{\d t}=\ha\bar q_1
\quad\text{and}\quad \frac{\d q_3}{\d t}=0.
\label{ev8eq3}
\end{gather}
The solutions of \eq{ev8eq3} are easily shown to be
\e
q_1(t)=A{\rm e}^{it/2}+B{\rm e}^{-it/2},\;\>
q_2(t)=i\bar A{\rm e}^{-it/2}-i\bar B{\rm e}^{it/2},\;\>
q_3(t)=C,
\label{ev8eq4}
\e
for $A,B,C\in\C$. Substituting these into \eq{ev8eq2} gives
\begin{align*}
\frac{\d p_1}{\d t}&=i\bar C{\rm e}^{it}+\ha\bar p_2,\qquad
\frac{\d p_2}{\d t}=-\bar C{\rm e}^{-it}-\ha\bar p_1\\
\text{and}\quad \frac{\d p_3}{\d t}&=-iA{\rm e}^{-it/2}-i\bar A{\rm e}^{it/2}
+iB{\rm e}^{-3it/2}-i\bar B{\rm e}^{3it/2}.
\end{align*}

The first two reduce to $\frac{\d^2 p_1}{\d t^2}+\frac{1}{4}p_1
=-(\ha C+\bar C){\rm e}^{it}$, and are then easily solved, and
the third gives $p_3$ by integration. Thus the solutions are
\begin{align*}
p_1(t)&=\bigl({\ts\frac{2}{3}}C+{\ts\frac{4}{3}}\bar C\bigr){\rm e}^{it}
+D{\rm e}^{it/2}+E{\rm e}^{-it/2},\\
p_2(t)&=-i\bigl({\ts\frac{2}{3}}C+{\ts\frac{4}{3}}\bar C\bigr){\rm e}^{-it}
-i\bar D{\rm e}^{-it/2}+i\bar E{\rm e}^{it/2},\\
p_3(t)&=2A{\rm e}^{-it/2}-2\bar A{\rm e}^{it/2}
-{\ts\frac{2}{3}}B{\rm e}^{-3it/2}-{\ts\frac{2}{3}}\bar B{\rm e}^{3it/2}+F,
\end{align*}
for $D,E,F\in\C$. This gives us the full solutions ${\bf z}_4,{\bf z}_5$
to~\eq{ev5eq11}.

Now ${\bf z}_4$ and ${\bf z}_5$ are required to satisfy equations
\eq{ev5eq5}--\eq{ev5eq7}. Calculation with the formulae above shows that
\begin{align*}
\om({\bf z}_1,{\bf z}_5)+\om({\bf z}_2,{\bf z}_5)
-\om({\bf z}_3,{\bf z}_4)&=-\Im F,\\
-\om({\bf z}_1,{\bf z}_4)+\om({\bf z}_2,{\bf z}_4)
+\om({\bf z}_3,{\bf z}_5)&=\Im C,\\
\text{and}\quad\om({\bf z}_4,{\bf z}_5)&=\Im(2A\bar D+2B\bar E+C\bar F).
\end{align*}
Thus \eq{ev5eq5}--\eq{ev5eq7} hold if and only if~$\Im C=\Im F=
\Im(2A\bar D+2B\bar E+C\bar F)=0$.

The two automatic solutions to \eq{ev5eq11} given by Corollary
\ref{ev5cor} correspond to $\Re C$ and $\Re F$, with $e=\Re C$ and
$f=\Re F$. Therefore, as in \S\ref{ev51}, we can use translational
symmetry in $\R^2$ to set $\Re C=\Re F=0$ without reducing the
set of SL 3-folds we construct. Thus, without loss of generality
we can set $C=F=0$, and then \eq{ev5eq5}--\eq{ev5eq7} hold if and
only if~$\Im(A\bar D+B\bar E)=0$.

We summarize our progress so far in the following proposition.

\begin{prop} Define ${\bf z}_1,\ldots,{\bf z}_5:\R\ra\C^3$ by
\ea
{\bf z}_1(t)&={\bf z}_2(t)=({\rm e}^{it},-i{\rm e}^{-it},0),\quad
{\bf z}_3(t)=(0,0,1),
\label{ev8eq5}\\
\begin{split}
{\bf z}_4(t)&=\bigl(D{\rm e}^{it/2}+E{\rm e}^{-it/2},
-i\bar D{\rm e}^{-it/2}+i\bar E{\rm e}^{it/2},\\
&\qquad 2A{\rm e}^{-it/2}-2\bar A{\rm e}^{it/2}-{\ts\frac{2}{3}}
B{\rm e}^{-3it/2}-{\ts\frac{2}{3}}\bar B{\rm e}^{3it/2}\bigr),
\end{split}
\label{ev8eq6}\\
\text{and}\quad
{\bf z}_5(t)&=\bigl(A{\rm e}^{it/2}+B{\rm e}^{-it/2},
i\bar A{\rm e}^{-it/2}-i\bar B{\rm e}^{it/2},0\bigr),
\label{ev8eq7}
\ea
where $A,B,D,E$ are complex numbers with\/ $\Im(A\bar D+B\bar E)=0$.
Then ${\bf z}_1,\ldots,{\bf z}_5$ satisfy \eq{ev5eq4}--\eq{ev5eq7}
and\/~\eq{ev5eq10}--\eq{ev5eq11}. 
\label{ev8prop2}
\end{prop}

Next we solve \eq{ev5eq11} for ${\bf z}_6$. Substituting \eq{ev8eq6}
and \eq{ev8eq7} into \eq{ev5eq11} gives a rather complicated formula
for $\d{\bf z}_6/\d t$, which may be integrated in the normal way.
We find that ${\bf z}_6=(r_1,r_2,r_3)$, where
\begin{align*}
r_1(t)&=\!-\!{\ts\frac{1}{6}}A\bar B{\rm e}^{2it}
\!+\!(A\bar A\!+\!{\ts\frac{1}{3}}B\bar B){\rm e}^{it}
\!-\!i(A^2\!+\!\bar AB)t\!-\!{\ts\frac{2}{3}}AB{\rm e}^{-it}
\!-\!{\ts\frac{1}{6}}B^2{\rm e}^{-2it}\!+\!G,\\
r_2(t)&={\ts\frac{i}{6}}\bar B^2{\rm e}^{2it}
\!-\!{\ts\frac{2i}{3}}\bar A\bar B{\rm e}^{it}\!+\!(\bar A^2\!-\!A\bar B)t
\!-\!i(A\bar A\!+\!{\ts\frac{1}{3}}B\bar B){\rm e}^{-it}
\!-\!{\ts\frac{i}{6}}\bar A B{\rm e}^{-2it}\!+\!H,\\
r_3(t)&=-\ha(A\bar E+\bar B D){\rm e}^{it}
-i\Re(A\bar D-B\bar E)t-\ha(\bar AE+B\bar D){\rm e}^{-it}+I,
\end{align*}
for $G,H,I\in\C$. Collecting all the above material together,
and setting $G=H=I=0$ for simplicity, we have proved:

\begin{thm} Let\/ $A,B,D,E\in\C$ with\/ $\Im(A\bar D\!+\!B\bar E)\!=\!0$,
and define $N$~to~be
\e
\begin{split}
\biggl\{&\Bigl(y_1^2{\rm e}^{it}
\!+\!y_1\bigl(D{\rm e}^{it/2}\!+\!E{\rm e}^{-it/2}\bigr)
\!+\!y_2\bigl(A{\rm e}^{it/2}\!+\!B{\rm e}^{-it/2}\bigr)
\!-\!{\ts\frac{1}{6}}A\bar B{\rm e}^{2it}\\
&\quad +(A\bar A\!+\!{\ts\frac{1}{3}}B\bar B){\rm e}^{it}
\!-\!i(A^2+\bar AB)t\!-\!{\ts\frac{2}{3}}AB{\rm e}^{-it}
\!-\!{\ts\frac{1}{6}}B^2{\rm e}^{-2it},\\
&-\!iy_1^2{\rm e}^{-it}\!+\!iy_1\bigl(-\bar D{\rm e}^{-it/2}
\!+\!\bar E{\rm e}^{it/2}\bigr)\!+\!iy_2\bigl(\bar A{\rm e}^{-it/2}
\!-\!\bar B{\rm e}^{it/2}\bigr)\!+\!{\ts\frac{i}{6}}\bar B^2{\rm e}^{2it}\\
&\quad
-\!{\ts\frac{2i}{3}}\bar A\bar B{\rm e}^{it}\!+\!(\bar A^2-A\bar B)t
\!-\!i(A\bar A\!+\!{\ts\frac{1}{3}}B\bar B){\rm e}^{-it}
\!-\!{\ts\frac{i}{6}}\bar A B{\rm e}^{-2it},\\
&y_1y_2\!+\!y_1\bigl(2A{\rm e}^{\!-it/2}\!-\!2\bar A{\rm e}^{it/2}
\!-\!{\ts\frac{2}{3}}B{\rm e}^{\!-3it/2}\!-\!{\ts\frac{2}{3}}
\bar B{\rm e}^{3it/2}\bigr)\!-\!\ha(A\bar E\!+\!\bar B D){\rm e}^{it}\\
&\quad -i\Re(A\bar D-B\bar E)t
-\ha(\bar AE+B\bar D){\rm e}^{-it}\Bigr):y_1,y_2,t\in\R\biggr\}.
\end{split}
\label{ev8eq8}
\e
Then $N$ is a special Lagrangian $3$-fold in $\C^3$. Furthermore,
any special Lagrangian $3$-fold in $\C^3$ constructed using Theorem
\ref{ev5thm} and satisfying case {\rm(iii)} of\/ \S\ref{ev7} is
isomorphic to one of this family under an automorphism of\/~$\C^3$.
\label{ev8thm2}
\end{thm}

\subsection{Discussion}
\label{ev83}

As in \S\ref{ev52}, here is a parameter count for the set of SL 3-folds
constructed in Theorem \ref{ev8thm2}. The 3-folds depend on 4 complex
numbers $A,B,D,E$, which is 8 real parameters, and they satisfy one
real equation $\Im(A\bar D+B\bar E)=0$, reducing the number of
parameters to 7. Should we reduce this further to allow for
isomorphisms between members of the family?

Well, in \S\ref{ev81} we used up the $\GL(2,\R)$ and $\SU(3)$ symmetries
to fix ${\bf z}_1,{\bf z}_2$ and ${\bf z}_3$, and in \S\ref{ev82} we used
the $\R^2$ translational symmetry to set $\Re C=\Re F=0$, and the $\C^3$
translational symmetry to fix $G=H=I=0$. So the $\GL(2,\R)\lt\R^2$ and
$\SU(3)\lt\C^3$ symmetry groups are both already fully accounted for.
The only remaining symmetry in the problem is translation in
time,~$t\mapsto t+c$.

Let $c\in\R$, and replace $t,A,B,D,E$ by $t',A',B',D',E'$, where
\begin{equation*}
t'=t+c,\;\>
A'={\rm e}^{ic/2}A,\;\>
B'={\rm e}^{3ic/2}B,\;\>
D'={\rm e}^{ic/2}D,\;\>
E'={\rm e}^{3ic/2}E.
\end{equation*}
Then a point $(z_1,z_2,z_3)$ in $N$ is replaced by $(z_1',z_2',z_3')$,
where
\begin{gather*}
z_1'={\rm e}^{ic}\bigl(z_1-i(A^2+\bar AB)c\bigr),\qquad
z_2'={\rm e}^{-ic}\bigl(z_2+(\bar A^2-A\bar B)c\bigr)\\
\text{and}\qquad z_3'=z_3-i\Re(A\bar D-B\bar E)c.
\end{gather*}
This map $(z_1,z_2,z_3)\mapsto(z_1',z_2',z_3')$ lies in $\SU(3)\lt\C^3$.

Therefore if $N'$ is defined as in \eq{ev8eq8}, but using constants
$A'={\rm e}^{ic/2}A$, $B'={\rm e}^{3ic/2}B$, $D'={\rm e}^{ic/2}D$ and
$E'={\rm e}^{3ic/2}E$, then $N'$ is isomorphic to $N$ under an
$\SU(3)\lt\C^3$ transformation. So we should reduce the number of real
parameters by one, and the family of distinct SL 3-folds from Theorem
\ref{ev8thm2} up to automorphisms of $\C^3$ has 6 dimensions. For
comparison, the family from Theorem \ref{ev5thm}, of which this is a
special case, has 9 dimensions.

Next we consider the conditions for the evolution to be periodic in $t$.
From \eq{ev8eq8} it is clear that the evolution will be periodic, with
period $4\pi$, if and only if the coefficients of $t$ vanish. That is,
the evolution is periodic if $A^2+\bar AB=0$, $\bar A^2-A\bar B=0$ and
$\Re(B\bar E-A\bar D)=0$. The first two equations hold if and only if
$A=0$, and then the second equation becomes $\Re(B\bar E)=0$. But when
$A=0$ the conditions on $A,B,D$ and $E$ in Theorem \ref{ev8thm2} become
$\Im(B\bar E)=0$. Thus~$B\bar E=0$.

Hence the evolution is periodic when $A=B\bar E=0$, that is, if either
\begin{itemize}
\item[(a)] $A=B=0$, or
\item[(b)] $A=E=0$.
\end{itemize}
In case (a), the 3-fold $N$ of \eq{ev8eq8} is given by
\e
\begin{split}
\Bigl\{\bigl(&y_1^2{\rm e}^{it}
+y_1(D{\rm e}^{it/2}+E{\rm e}^{-it/2}),\\
-&iy_1^2{\rm e}^{-it}+iy_1(-\bar D{\rm e}^{-it/2}
+\bar E{\rm e}^{it/2}),y_1y_2\bigr):y_1,y_2,t\in\R\Bigr\}.
\end{split}
\label{ev8eq9}
\e
Notice that for each $(z_1,z_2,z_3)$ in $N$, the third coordinate
$z_3=y_1y_2$ is {\it real}. This implies that $N$ is a subset of
$\Si\t\R$ in $\C^2\t\C$, where $\Si$ is an SL 2-fold in $\C^2$, which
is in fact a {\it complex parabola} with respect to an alternative
complex structure on $\C^2$. However, the parametrization $(y_1,y_2,t)$
does not respect the product structure.

In case (b), the 3-fold $N$ of \eq{ev8eq8} is given by
\begin{align*}
\Bigl\{\bigl(&y_1^2{\rm e}^{it}+y_1D{\rm e}^{it/2}+y_2B{\rm e}^{-it/2}
+{\ts\frac{1}{3}}B\bar B{\rm e}^{it}-{\ts\frac{1}{6}}B^2{\rm e}^{-2it},\\
-&iy_1^2{\rm e}^{-it}-iy_1\bar D{\rm e}^{-it/2}-iy_2\bar B{\rm e}^{it/2}
-{\ts\frac{i}{3}}B\bar B{\rm e}^{-it}+{\ts\frac{i}{6}}\bar B^2{\rm e}^{2it},\\
&y_1y_2-y_1({\ts\frac{2}{3}}B{\rm e}^{-3it/2}+{\ts\frac{2}{3}}
\bar B{\rm e}^{3it/2})-\ha\bar B D{\rm e}^{it}
-\ha B\bar D{\rm e}^{-it}\bigr):y_1,y_2,t\in\R\Bigr\}.
\end{align*}
Here for each $(z_1,z_2,z_3)$ in $N$ we have $z_2=-i\bar z_1$ and
$z_3$ is real. That is, $N$ is a subset of the special Lagrangian
3-plane
\begin{equation*}
\bigl\{(z_1,z_2,z_3)\in\C^3:z_2=-i\bar z_1,\quad z_3\in\R\bigr\}.
\end{equation*}
So case (b) is just an unusual parametrization of an $\R^3$ in $\C^3$.
Thus, the periodic solutions in Theorem \ref{ev8thm2} are not very
interesting.

Finally we note that the SL 3-folds of Theorem \ref{ev8thm2} are
ruled by straight lines. Define $\Phi:\R^3\ra\C^3$ as in
\eq{ev5eq21}, so that $\Phi(y_1,y_2,t)$ is the vector in \eq{ev8eq8}.
Then, as \eq{ev8eq8} contains no terms in $y_2^2$, for each fixed
$y_1,t\in\R$ the set $\bigl\{\Phi(y_1,y_2,t):y_2\in\R\bigr\}$ is
a real straight line in $\C^3$. So $N$ is fibred by straight lines.
Such submanifolds are called {\it ruled}. Ruled special Lagrangian
3-folds are the subject of~\cite{Joyc5}.

\section{Case (iv) of \S\ref{ev7}}
\label{ev9}

We now move on to study case (iv) of \S\ref{ev7}, in which
${\bf z}_1,{\bf z}_2$ and ${\bf z}_3$ are linearly independent
for generic $t$. In \S\ref{ev91} we use the symmetry of the 
construction to reduce ${\bf z}_1,{\bf z}_2,{\bf z}_3$ to a
convenient form in coordinates, and then in Theorem \ref{ev9thm}
we solve equation \eq{ev5eq10} explicitly using material 
from~\cite[\S 6]{Joyc3}. 

Section \ref{ev92} rewrites equations \eq{ev5eq11} and 
\eq{ev5eq12} for ${\bf z}_4,{\bf z}_5$ and ${\bf z}_6$
and derives some properties of their solutions, and 
\S\ref{ev93} discusses the difficulties of solving them.
In two cases of Theorem \ref{ev9thm} we can solve \eq{ev5eq11}
and \eq{ev5eq12} explicitly, and we will do this in sections
\ref{ev10} and~\ref{ev11}.

\subsection{Solving the equations for ${\bf z}_1,{\bf z}_2$ and 
${\bf z}_3$}
\label{ev91}

We begin by choosing a convenient form for ${\bf z}_1,{\bf z}_2$ 
and ${\bf z}_3$ in coordinates, using the symmetry groups 
$\GL(2,\R)$ and $\SU(3)$ of the construction, as we did in 
Theorem~\ref{ev8thm1}.

\begin{prop} Any ${\bf z}_1,{\bf z}_2,{\bf z}_3:\R\ra\C^3$
satisfying \eq{ev5eq4}, \eq{ev5eq10} and part\/ {\rm(iv)} of\/
\S\ref{ev7} are equivalent under the natural actions of\/ 
$\SL(2,\R)$ and\/ $\SU(3)$~to
\e
\hat{\bf z}_1=(w_1,0,0),\quad
\hat{\bf z}_2=(0,w_2,0)\quad\text{and}\quad
\hat{\bf z}_3=(0,0,w_3)
\label{ev9eq1}
\e
for differentiable functions $w_1,w_2,w_3:\R\ra\C$.
\label{ev9prop1}
\end{prop}

\begin{proof} Suppose for simplicity that ${\bf z}_1(0),{\bf z}_2(0)$
and ${\bf z}_3(0)$ are linearly independent in $\C^3$. Consider the
action of $\SL(2,\R)$ in $\GL(2,\R)$ upon 
${\bf z}_1,{\bf z}_2$ and 
${\bf z}_3$ given in \eq{ev5eq15}--\eq{ev5eq17}, where $\de=ad-bc=1$.
We can think of this as an action of $\SL(2,\R)$ on $\R^3$, where
$(a_1,a_2,a_3)$ in $\R^3$ is mapped to $a_1{\bf z}_1+a_2{\bf z}_2+
a_3{\bf z}_3$. This action preserves the quadratic form 
$a_1^2-a_2^2-a_3^2$ upon~$\R^3$.

Now there is a second, positive definite quadratic form on $\R^3$
given by
\begin{equation*}
(a_1,a_2,a_3)\mapsto\bms{a_1{\bf z}_1(0)+a_2{\bf z}_2(0)+a_3{\bf z}_3(0)}.
\end{equation*}
By standard results on simultaneous diagonalization of quadratic
forms, there is a basis of $\R^3$ in which both quadratic forms
are diagonal. Choose an element of $\SL(2,\R)$ which takes this
basis to vectors proportional to $(1,0,0)$, $(0,1,0)$ and $(0,0,1)$.

Let ${\bf z}_1',{\bf z}_2',{\bf z}_3'$ be the transforms of
${\bf z}_1,{\bf z}_2,{\bf z}_3$ under this element of $\SL(2,\R)$.
Then the quadratic form
\begin{equation*}
(a_1,a_2,a_3)\mapsto\bms{a_1{\bf z}_1'(0)+a_2{\bf z}_2'(0)+a_3{\bf z}_3'(0)}
\end{equation*}
is diagonal with respect to the standard basis of $\R^3$. That is,
${\bf z}_1'(0),{\bf z}_2'(0)$ and ${\bf z}_3'(0)$ are {\it orthogonal}.
Now from \eq{ev5eq4} we see that ${\bf z}_1'(0),{\bf z}_2'(0),
{\bf z}_3'(0)$ span a Lagrangian plane in $\C^3$. But any three 
nonzero orthogonal elements in a Lagrangian plane in $\C^3$ are 
conjugate under a matrix in $\SU(3)$ to $(w_1,0,0)$, $(0,w_2,0)$ 
and $(0,0,w_3)$ for some~$w_j\in\C\sm\{0\}$. 

Applying this $\SU(3)$ matrix to ${\bf z}_1',{\bf z}_2',{\bf z}_3'$ 
gives $\hat{\bf z}_1,\hat{\bf z}_2,\hat{\bf z}_3$, where \eq{ev9eq1}
holds for $t=0$. But it is easy to see from \eq{ev5eq10} that if
\eq{ev9eq1} holds at $t=0$ then it holds for all $t$. This completes
the proof.
\end{proof}

So suppose that ${\bf z}_1,{\bf z}_2,{\bf z}_3$ are given by
\e
{\bf z}_1=(w_1,0,0),\quad
{\bf z}_2=(0,w_2,0)\quad\text{and}\quad
{\bf z}_3=(0,0,w_3)
\label{ev9eq2}
\e
for differentiable functions $w_1,w_2,w_3:\R\ra\C$ which are all
nonzero for generic $t$. Then \eq{ev5eq4} automatically holds,
and the evolution equations \eq{ev5eq10} are equivalent to the 
o.d.e.s
\e
\frac{\d w_1}{\d t}=\ov{w_2w_3},\quad
\frac{\d w_2}{\d t}=-\,\ov{w_3w_1},\quad
\frac{\d w_3}{\d t}=-\,\ov{w_1w_2}.
\label{ev9eq3}
\e
These are the same as equation (31) of \cite[Th.~6.1]{Joyc3}. Thus we
can use the material of \cite[\S 6]{Joyc3} to understand their solutions
very explicitly. In particular, \cite[Prop.~6.2]{Joyc3} and the discussion
after it gives

\begin{prop} For any given initial data $w_1(0),w_2(0),w_3(0)$, 
solutions $w_j(t)$ of\/ \eq{ev9eq3} exist for all\/ $t\in\R$. 
Wherever the $w_j(t)$ are nonzero, these functions may be written 
\begin{equation*}
w_1={\rm e}^{i\th_1}\sqrt{\al_1+u},\quad
w_2={\rm e}^{i\th_2}\sqrt{\al_2-u}\quad\text{and}\quad
w_3={\rm e}^{i\th_3}\sqrt{\al_3-u},
\end{equation*}
where $\al_j\in\R$ and\/ $u,\th_1,\th_2,\th_3:\R\ra\R$ are 
differentiable functions. Define
\begin{equation*}
Q(u)=(\al_1+u)(\al_2-u)(\al_3-u)\quad\text{and}\quad
\th=\th_1+\th_2+\th_3.
\end{equation*}
Then $Q(u)^{1/2}\sin\th\equiv A$ for some $A\in\R$, and\/
$u$ and\/ $\th_j$ satisfy
\begin{alignat*}{2}
\Bigl(\frac{\d u}{\d t}\Bigr)^2&=4\bigl(Q(u)-A^2\bigr),\quad&
\frac{\d\th_1}{\d t}&=-\,\frac{A}{\al_1+u},\\
\frac{\d\th_2}{\d t}&=\frac{A}{\al_2-u}\quad&\text{and}\qquad 
\frac{\d\th_3}{\d t}&=\frac{A}{\al_3-u}.
\end{alignat*}
\label{ev9prop2}
\end{prop}

As in \cite[\S 6]{Joyc3}, the solutions to these equations 
involve the {\it Jacobi elliptic functions}, to which we now give 
a brief introduction. The following material can be found in
Chandrasekharan~\cite[Ch.~VII]{Chan}.

For each $k\in[0,1]$, the Jacobi elliptic functions $\sn(t,k)$, $\cn(t,k)$,
$\dn(t,k)$ with modulus $k$ are the unique solutions to the o.d.e.s
\begin{align*}
\bigl({\ts\frac{\d}{\d t}}\sn(t,k)\bigr)^2&=\bigl(1-\sn^2(t,k)\bigr)
\bigl(1-k^2\sn^2(t,k)\bigr),\\
\bigl({\ts\frac{\d}{\d t}}\cn(t,k)\bigr)^2&=\bigl(1-\cn^2(t,k)\bigr)
\bigl(1-k^2+k^2\cn^2(t,k)\bigr),\\
\bigl({\ts\frac{\d}{\d t}}\dn(t,k)\bigr)^2&=-\bigl(1-\dn^2(t,k)\bigr)
\bigl(1-k^2-\dn^2(t,k)\bigr),
\end{align*}
with initial conditions
\begin{alignat*}{3}
\sn(0,k)&=0,\quad & \cn(0,k)&=1,\quad & \dn(0,k)&=1,\\
{\ts\frac{\d}{\d t}}\sn(0,k)&=1,\quad&
{\ts\frac{\d}{\d t}}\cn(0,k)&=0,\quad&
{\ts\frac{\d}{\d t}}\dn(0,k)&=0.
\end{alignat*}
They satisfy the identities
\begin{equation*}
\sn^2(t,k)+\cn^2(t,k)=1 \quad\text{and}\quad k^2\sn^2(t,k)+\dn^2(t,k)=1,
\end{equation*}
and the differential equations
\begin{gather*}
{\ts\frac{\d}{\d t}}\sn(t,k)=\cn(t,k)\dn(t,k),\qquad
{\ts\frac{\d}{\d t}}\cn(t,k)=-\sn(t,k)\dn(t,k)\\
\text{and}\qquad {\ts\frac{\d}{\d t}}\dn(t,k)=-k^2\sn(t,k)\cn(t,k).
\end{gather*}
When $k=0$ or 1 they reduce to trigonometric functions:
\begin{alignat}{2}
\sn(t,0)&=\sin t,&\quad \cn(t,0)&=\cos t,\quad \dn(t,0)=1,
\label{ev9eq4}\\
\sn(t,1)&=\tanh t&,\quad \cn(t,1)&=\dn(t,1)=\sech t.
\label{ev9eq5}
\end{alignat}
For $k\in[0,1)$ the Jacobi elliptic functions are periodic in~$t$.

As in \cite[\S 5.4]{Joyc3}, we choose $\al_1,\al_2,\al_3$ uniquely such that
$\al_j>0$ and $\frac{1}{\al_1}=\frac{1}{\al_2}+\frac{1}{\al_3}$. This then
implies that $Q$ has a maximum at 0, and $0\le A^2\le \al_1\al_2\al_3$.
From the results of \cite[\S 5.4]{Joyc3} and \cite[\S 6]{Joyc3} we can
then deduce the following theorem.

\begin{thm} In the situation above, suppose that\/ $u$ has a minimum 
at\/ $t=0$ and that\/ $\th_2(0)=\th_3(0)=0$, $A\ge 0$ and\/ 
$\al_2\le\al_3$. Then exactly one of the following four cases holds.
\begin{itemize}
\item[{\rm(a)}] $A=0$ and\/ $\al_2=\al_3$, and\/ $w_1,w_2,w_3$ are
given by
\e
\begin{gathered}
w_1(t)=\sqrt{3\al_1}\tanh\,\bigl(t\sqrt{3\al_1}\,\bigr)\qquad\text{and}\\
w_2(t)=w_3(t)=\sqrt{3\al_1}\,\sech\bigl(t\sqrt{3\al_1}\,\bigr).
\end{gathered}
\label{ev9eq6}
\e
\item[{\rm(b)}] $A=0$ and\/ $\al_2<\al_3$, and\/ $w_1,w_2,w_3$ are
given by
\begin{gather*}
w_1=(\al_1+\al_2)^{1/2}\sn(\si t,\tau),\quad
w_2=(\al_1+\al_2)^{1/2}\cn(\si t,\tau)\\
\text{and}\qquad w_3=(\al_1+\al_3)^{1/2}\dn(\si t,\tau),
\end{gather*}
where $\si=(\al_1+\al_3)^{1/2}$ and\/~$\tau=(\al_1+\al_2)^{1/2}
(\al_1+\al_3)^{-1/2}$. Note that\/ $w_1,w_2,w_3$ are periodic in~$t$.
\item[{\rm(c)}] $0<A<(\al_1\al_2\al_3)^{1/2}$. Let the roots of\/ 
$Q(u)-A^2$ be $\ga_1,\ga_2,\ga_3$, ordered such that\/ 
$\ga_1\le 0\le\ga_2\le\ga_3$. Then $u$ and\/ $\th_1,\th_2,\th_3$
are given by
\begin{align*}
u(t)&=\ga_1+(\ga_2-\ga_1)\sn^2(\si t,\tau),\\
\th_1(t)&=\th_1(0)-A\int_0^t\frac{\d s}{\al_1+\ga_1+(\ga_2-\ga_1)
\sn^2(\si s,\tau)},\\
\th_2(t)&=A\int_0^t\frac{\d s}{\al_2-\ga_1-(\ga_2-\ga_1)\sn^2(\si s,\tau)}\\
\text{and}\quad
\th_3(t)&=A\int_0^t\frac{\d s}{\al_3-\ga_1-(\ga_2-\ga_1)\sn^2(\si s,\tau)},
\end{align*}
where $\si=(\ga_3-\ga_1)^{1/2}$ and\/~$\tau=(\ga_2-\ga_1)^{1/2}
(\ga_3-\ga_1)^{-1/2}$.
\item[{\rm(d)}] $A=(\al_1\al_2\al_3)^{1/2}$. Define $a_1,a_2,a_3\in\R$ by
\begin{equation*}
a_1=-\al_1^{-1/2}(\al_2\al_3)^{1/2},\;\>
a_2=\al_2^{-1/2}(\al_3\al_1)^{1/2},\;\>
a_3=\al_3^{-1/2}(\al_1\al_2)^{1/2}.
\end{equation*}
Then $a_1+a_2+a_3=0$, as $\frac{1}{\al_1}=\frac{1}{\al_2}+\frac{1}{\al_3}$,
and\/ $w_1,w_2,w_3$ are given by
\begin{equation*}
w_1(t)=i\sqrt{\al_1}\,{\rm e}^{ia_1t},\quad
w_2(t)=\sqrt{\al_2}\,{\rm e}^{ia_2t}
\quad\text{and}\quad w_3(t)=\sqrt{\al_3}\,{\rm e}^{ia_3t}.
\end{equation*}
\end{itemize}
\label{ev9thm}
\end{thm}

Here we have made a number of assumptions to simplify the 
formulae, none of which really reduces the generality of the result. 
Supposing $u$ has a minimum at $t=0$ means that in part (c) we get 
$u(t)=\ga_1+(\ga_2-\ga_1)\sn^2(\si t,\tau)$ rather than $u(t)=\ga_1+
(\ga_2-\ga_1)\sn^2(\si t+\up,\tau)$ for some $\up\in\R$, and similarly
for parts (a) and (b). Thus the assumption can be removed by adding a
constant to~$t$.

Setting $\th_2(0)=\th_3(0)=0$ has the effect of tidying up the 
constant phase factors in $w_1,w_2,w_3$, so that the $w_j$ are real 
in cases (a) and (b). We can change the sign of $A$ by replacing 
$w_1$ by $-w_1$ and $t$ by $-t$, and we can swap $\al_2$ and 
$\al_3$ by swapping $w_2$ and $w_3$. Thus the assumptions that 
$A\ge 0$ and $\al_2\le\al_3$ are no real restriction.

Note that case (a) follows from case (b), as $\tau=1$ when $\al_2=\al_3$, 
so the $\sn$, $\cn$ and $\dn$ expressions reduce to $\tanh$ and $\sech$
by \eq{ev9eq5}. Also, in case (a) the equations $\al_2=\al_3$ and 
$\frac{1}{\al_1}=\frac{1}{\al_2}+\frac{1}{\al_3}$ imply that $\al_2=\al_3=2\al_1$, 
so we have written the solutions solely in terms of~$\al_1$.

\subsection{Rewriting the equations for ${\bf z}_4,{\bf z}_5$ and
${\bf z}_6$}
\label{ev92}

Next we shall we rewrite equations \eq{ev5eq11} and \eq{ev5eq12}. 
Choose the slightly unusual form
\e
{\bf z}_4=(p_1,p_2,q_3)\quad\text{and}\quad {\bf z}_5=(q_1,-q_2,p_3),
\label{ev9eq7}
\e
where $p_j,q_j:(-\ep,\ep)\ra\C$ are differentiable functions.
Then \eq{ev5eq11} becomes
\begin{align*}
\frac{\d}{\d t}(p_1,p_2,q_3)&=(w_1,w_2,0)\t(q_1,-q_2,p_3)-
(0,0,w_3)\t(p_1,p_2,q_3),\\
\frac{\d}{\d t}(q_1,-q_2,p_3)&=(-w_1,w_2,0)\t(p_1,p_2,q_3)+
(0,0,w_3)\t(q_1,-q_2,p_3).
\end{align*}
Expanding using \eq{ev5eq9}, this yields
\begin{gather}
\begin{gathered}
\frac{\d p_1}{\d t}=\ha(\bar w_2\bar p_3+\bar w_3\bar p_2),\qquad
\frac{\d p_2}{\d t}=-\ha(\bar w_3\bar p_1+\bar w_1\bar p_3)\\
\text{and}\qquad \frac{\d p_3}{\d t}=-\ha(\bar w_1\bar p_2+\bar w_2\bar p_1),
\end{gathered}
\label{ev9eq8}\\
\begin{gathered}
\frac{\d q_1}{\d t}=\ha(\bar w_2\bar q_3+\bar w_3\bar q_2),\qquad
\frac{\d q_2}{\d t}=-\ha(\bar w_3\bar q_1+\bar w_1\bar q_3)\\
\text{and}\qquad \frac{\d q_3}{\d t}=-\ha(\bar w_1\bar q_2+\bar w_2\bar q_1).
\end{gathered}
\label{ev9eq9}
\end{gather}
These are real linear in the $p_j$ and $q_j$, and are the same
equations, so that any solution to \eq{ev9eq8} in $p_1,p_2,p_3$
is also a solution to \eq{ev9eq9} in $q_1,q_2,q_3$. Notice also 
that from \eq{ev9eq3}, $p_j=w_j$ solves \eq{ev9eq8}, and $q_j=w_j$ 
solves \eq{ev9eq9}. These are the two automatic solutions of 
\eq{ev5eq11} we get from Corollary~\ref{ev5cor}.

Writing ${\bf z}_6=(r_1,r_2,r_3)$ for $r_j:(-\ep,\ep)\ra\C$, 
equation \eq{ev5eq12} becomes
\e
\begin{gathered}
\frac{\d r_1}{\d t}=\ha(\bar p_2\bar p_3+\bar q_3\bar q_2),\qquad
\frac{\d r_2}{\d t}=\ha(\bar q_3\bar q_1-\bar p_1\bar p_3)\\
\text{and}\qquad 
\frac{\d r_3}{\d t}=-\ha(\bar p_1\bar q_2+\bar p_2\bar q_1),
\end{gathered}
\label{ev9eq10}
\e
and integrating these equations gives $r_1,r_2$ and~$r_3$.

In our next result we work out equations \eq{ev5eq4}--\eq{ev5eq7} 
for ${\bf z}_1,\ldots,{\bf z}_6$ of this form. The proof is very
easy, and we omit it.

\begin{lem} When ${\bf z}_1,\ldots,{\bf z}_6$ are defined
as above, equation \eq{ev5eq4} holds automatically, and
equations \eq{ev5eq5}--\eq{ev5eq7} are equivalent to
\e
\begin{split}
\Im(w_1\bar p_1-w_2\bar p_2-w_3\bar p_3)&=
\Im(w_1\bar q_1-w_2\bar q_2-w_3\bar q_3)\\
&=\Im(p_1\bar q_1-p_2\bar q_2-p_3\bar q_3)=0.
\end{split}
\label{ev9eq11}
\e
Furthermore, if the $w_j,p_j$ and\/ $q_j$ satisfy \eq{ev9eq3},
\eq{ev9eq8} and\/ \eq{ev9eq9}, but not necessarily \eq{ev9eq11},
then $\Im(w_1\bar p_1-w_2\bar p_2-w_3\bar p_3)$,
$\Im(w_1\bar q_1-w_2\bar q_2-w_3\bar q_3)$ and\/
$\Im(p_1\bar q_1-p_2\bar q_2-p_3\bar q_3)$ are constant.
\label{ev9lem1}
\end{lem}

The last part can be useful in solving equations \eq{ev9eq8}
and \eq{ev9eq9}, once $w_1,w_2$ and $w_3$ are chosen, because
each solution $p_1,p_2,p_3$ of \eq{ev9eq8} gives a conserved 
quantity in \eq{ev9eq9}, and so reduces the number of real
variables by one. Of course, \eq{ev9eq8} and \eq{ev9eq9} are
really the same, and have the same six-dimensional space of
solutions. Thus, given three independent solutions of 
\eq{ev9eq8}, we may be able to find the other three 
solutions by a kind of matrix inversion.

The equations simplify further in cases (a) and (b) of Theorem 
\ref{ev9thm}, when $w_1,w_2,w_3$ are real. Then the real parts
of \eq{ev9eq8} involve only $\Re(p_j)$, and the imaginary
parts of \eq{ev9eq8} involve only $\Im(p_j)$. So using
Lemma \ref{ev9lem1} we deduce:

\begin{lem} In cases {\rm(a)} and\/ {\rm(b)} of Theorem 
\ref{ev9thm}, the solutions $p_1,p_2,p_3$ of\/ \eq{ev9eq8} 
are of the form $p_j=u_j+iv_j$, where $u_j,v_j:\R\ra\R$ are 
differentiable functions satisfying
\begin{gather}
\begin{gathered}
\frac{\d u_1}{\d t}=\ha(w_2u_3+w_3u_2),\qquad
\frac{\d u_2}{\d t}=-\ha(w_3u_1+w_1 u_3)\\
\text{and}\qquad \frac{\d u_3}{\d t}=-\ha(w_1u_2+w_2u_1),
\end{gathered}
\label{ev9eq12}\\
\begin{gathered}
\frac{\d v_1}{\d t}=-\ha(w_2 v_3+w_3 v_2),\qquad
\frac{\d v_2}{\d t}=\ha(w_3v_1+w_1 v_3)\\
\text{and}\qquad \frac{\d v_3}{\d t}=\ha(w_1v_2+w_2v_1).
\end{gathered}
\label{ev9eq13}
\end{gather}
Furthermore, if\/ \eq{ev9eq12} and\/ \eq{ev9eq13} hold, then 
$u_1v_1-u_2v_2-u_3v_3$ is constant.
\label{ev9lem2}
\end{lem}

If we know the full solutions to \eq{ev9eq12} we can determine
the solutions to \eq{ev9eq13}, and vice versa. Here is how.
Suppose that $u_1^j,u_2^j,u_3^j$ are solutions to \eq{ev9eq12}
for $j=1,2,3$, and $v_1^j,v_2^j,v_3^j$ solutions to \eq{ev9eq13}
for $j=1,2,3$. Then the last part of the lemma shows that
\begin{equation*}
\begin{pmatrix}u_1^1&u_2^1&u_3^1\\u_1^2&u_2^2&u_3^2\\u_1^3&u_2^3&u_3^3
\end{pmatrix}\begin{pmatrix}1&0&0\\0&-1&0\\0&0&-1\end{pmatrix}
\begin{pmatrix}v^1_1&v^2_1&v^3_1\\v^1_2&v^2_2&v^3_2\\v^1_3&v^2_3&v^3_3
\end{pmatrix}
\end{equation*}
is a constant matrix. In particular, if for $j=1,2,3$ the 
$(u_1^j,u_2^j,u_3^j)$ are linearly independent for some $t$,
then they are linearly independent for all $t$, and we may
{\it define} the $v_i^j$ by
\begin{equation*}
\begin{pmatrix}v^1_1&v^2_1&v^3_1\\v^1_2&v^2_2&v^3_2\\v^1_3&v^2_3&v^3_3
\end{pmatrix}=\begin{pmatrix}u_1^1&-u_2^1&-u_3^1\\
u_1^2&-u_2^2&-u_3^2\\u_1^3&-u_2^3&-u_3^3\end{pmatrix}^{-1},
\end{equation*}
where the matrix inverse exists. The $v_1^j,v_2^j,v_3^j$ will
then be solutions to \eq{ev9eq13} for $j=1,2,3$, and span the
full set of solutions.

\subsection{Discussion}
\label{ev93}

So far we have solved equation \eq{ev5eq10} for ${\bf z}_1,
{\bf z}_2,{\bf z}_3$ explicitly in Theorem \ref{ev9thm}, and 
we have rewritten equations \eq{ev5eq11} and \eq{ev5eq12} for 
${\bf z}_4,{\bf z}_5$ and ${\bf z}_6$ in a convenient form, 
and said something about the properties of their solutions.
However, we have not yet solved \eq{ev5eq11} and \eq{ev5eq12}
explicitly.

In cases (a) and (d) of Theorem \ref{ev9thm} the author has been able 
to solve equations \eq{ev5eq11} and \eq{ev5eq12} fairly explicitly.
These cases will be treated at length in sections \ref{ev10} and 
\ref{ev11} respectively. However, in cases (b) and (c) of Theorem 
\ref{ev9thm} the author has made little progress in solving 
\eq{ev5eq11} and~\eq{ev5eq12}.

In case (c), it may help to split into two cases $\al_2=\al_3$ 
and $\al_2<\al_3$. The case $\al_2=\al_3$ is simpler, because 
$w_2\equiv w_3$, and may be more amenable to explicit solution.
In particular, \eq{ev9eq9} is unchanged by swapping round
$p_2$ and $p_3$. Thus we can separately consider symmetric 
solutions with $p_2=p_3$, and antisymmetric solutions with
$p_1=0$ and~$p_2=-p_3$.

It would be interesting to know about the {\it periodicity} of
solutions to equations \eq{ev9eq8}--\eq{ev9eq10} in cases (b) and
(c) of Theorem \ref{ev9thm}. In particular, in case (b) the $w_j$ 
are periodic with period $T>0$. If equations \eq{ev9eq8}--\eq{ev9eq10} 
admitted interesting solutions with period $nT$ for some integer 
$n\ge 1$, then the corresponding SL 3-fold $N$ would be a closed, 
immersed copy of $\R^2\t{\cal S}^1$ rather than~$\R^3$.

In fact a weaker form of periodicity would be sufficient, where the
map $\Phi$ of \eq{ev5eq21} satisfies $\Phi(y_1+c_1,y_2+c_2,t+nT)=
\Phi(y_1,y_2,t)$ for $c_1,c_2\in\R$. We also saw in \cite[\S 6]{Joyc3} that
many of the solutions $w_1,w_2,w_3$ in part (c) are periodic, so we
can ask the same question.

\section{Case (a) of Theorem \ref{ev9thm}}
\label{ev10}

We now study the SL 3-folds $N$ in $\C^3$ which come from the
construction of Theorem \ref{ev5thm} when ${\bf z}_1,{\bf z}_2$
and ${\bf z}_3$ are as in case (a) of Theorem \ref{ev9thm}. For
simplicity we set $\al_1=\frac{1}{3}$, so that the functions 
$w_j$ in case (a) of Theorem \ref{ev9thm} are $w_1(t)=\tanh t$ 
and $w_2(t)=w_3(t)=\sech t$. Since we are free to rescale 
$\al_1,\al_2,\al_3$ using dilations in $\GL(2,\R)$, under the 
$\GL(2,\R)$-action discussed in \S\ref{ev51}, fixing $\al_1$
in this way does not change the resulting set of SL 3-folds.

\begin{thm} In the situation of\/ \S\ref{ev9}, set\/ 
$w_1(t)=\tanh t$ and\/ $w_2(t)=w_3(t)=\sech t$, as in case {\rm(a)} 
of Theorem \ref{ev9thm}. Equation \eq{ev9eq8} then becomes
\e
\begin{gathered}
\frac{\d p_1}{\d t}=\ha(\bar p_2+\bar p_3)\sech t,\qquad
\frac{\d p_2}{\d t}=-\ha(\bar p_1\sech t+\bar p_3\tanh t)\\
\text{and}\qquad \frac{\d p_3}{\d t}=
-\ha(\bar p_1\sech t+\bar p_2\tanh t).
\end{gathered}
\label{ev10eq1}
\e
The full solutions to these equations are
\ea
p_1&=A\,\tanh t\!+\!B\bigl(f(t)\tanh t\!-\!2\sqrt{\cosh t}\bigr)
\!-\!\frac{iD\,f(t)}{\cosh t}\!+\!\frac{iE}{\sqrt{\cosh t}},
\label{ev10eq2}\\
\begin{split}
p_2&=A\,\sech t+B\,f(t)\sech t+C\,\sqrt{\cosh t}\\
&\quad+iD\Bigl(\cosh t-\frac{f(t)\sinh t}{2\sqrt{\cosh t}}\Bigr)
+\frac{iE\,\sinh t}{2\sqrt{\cosh t}}+\frac{iF}{\sqrt{\cosh t}},
\end{split}
\label{ev10eq3}\\
\begin{split}
p_3&=A\,\sech t+B\,f(t)\sech t-C\,\sqrt{\cosh t}\\
&\quad+iD\Bigl(\cosh t-\frac{f(t)\sinh t}{2\sqrt{\cosh t}}\Bigr)
+\frac{iE\,\sinh t}{2\sqrt{\cosh t}}-\frac{iF}{\sqrt{\cosh t}},
\end{split}
\label{ev10eq4}
\ea
where $A,B,C,D,E,F\in\R$ and\/~$f(t)=\int_0^t\sqrt{\cosh s}\,\d s$.
\label{ev10thm1}
\end{thm}                                                            

\begin{proof} Equation \eq{ev10eq1} follows immediately. One can 
verify that \eq{ev10eq2}--\eq{ev10eq4} are solutions to \eq{ev10eq1} 
by subtituting them in, and using $\frac{\d f}{\d t}=\sqrt{\cosh t}$.
The six solutions of \eq{ev10eq1} this gives, with coefficients 
$A,\ldots,F$, are easily seen to be linearly independent. But as 
\eq{ev10eq1} is a well-behaved first-order o.d.e., its solutions
are determined by the initial data $p_1(0),p_2(0),p_3(0)$, which
comprise 3 complex or 6 real numbers. Thus \eq{ev10eq1} can have
no more than 6 linearly independent solutions, so 
\eq{ev10eq2}--\eq{ev10eq4} are the full solutions to~\eq{ev10eq1}.
\end{proof}

To derive equations \eq{ev10eq2}--\eq{ev10eq4}, the author used 
Lemma \ref{ev9lem2}, the fact that $p_j=w_j$ is automatically a 
solution, and the symmetry between $p_2$ and $p_3$, which means
that we can consider separately solutions with $p_2=p_3$, and 
those with $p_1=0$ and $p_2=-p_3$. Next we work out the conditions
for \eq{ev5eq4}--\eq{ev5eq7} to hold for these solutions.

\begin{lem} Let\/ $w_1,w_2,w_3$ be as above, define $p_1,p_2,p_3$ 
by \eq{ev10eq2}--\eq{ev10eq4} for some $A,\ldots,F\in\R$, and
define $q_1,q_2,q_3$ by \eq{ev10eq2}--\eq{ev10eq4}, but replacing
$p_j$ with\/ $q_j$ and\/ $A,\ldots,F$ with\/ $A',\ldots,F'\in\R$. This 
defines ${\bf z}_1,\ldots,{\bf z}_5$, as in \S\ref{ev9}. Equations
\eq{ev5eq4}--\eq{ev5eq7} hold for these ${\bf z}_1,\ldots,{\bf z}_5$
if and only if
\e
D=D'=0\quad\text{and}\quad AD'+BE'+CF'=A'D+B'E+C'F.
\label{ev10eq5}
\e
\label{ev10lem}
\end{lem}

\begin{proof} When $t=0$, we have $w_1(0)=1$, $w_2(0)=w_3(0)=1$,
$p_1(0)=-2B+iE$, $p_2(0)=A+C+iD+iF$, $p_3(0)=A-C+iD-iF$, 
$q_1(0)=-2B'+iE'$, $q_2(0)=A'+C'+iD'+iF'$ and $q_3(0)=A'-C'+iD'-iF'$.
Thus \eq{ev9eq11} holds at $t=0$ if and only if \eq{ev10eq5} holds,
so by Lemma \ref{ev9lem1} equations \eq{ev5eq4}--\eq{ev5eq7} hold when 
$t=0$ if and only if \eq{ev10eq5} holds. But \eq{ev5eq4}--\eq{ev5eq7} 
hold at $t=0$ if and only if they hold for all~$t$.
\end{proof}

Thus we should set $D$ and $D'$ to zero. But we can show by changing 
coordinates in $\R^2$ from $(y_1,y_2)$ to $(y_1+A,y_2+A')$ as in
\S\ref{ev51} that we are also free to set $A$ and $A'$ to zero, 
without restricting the SL 3-folds constructed in Theorem 
\ref{ev5thm}. So the remaining interesting parameters are $B,C,E,F$ 
and $B',C',E',F'$, which must satisfy the restriction~$BE'+CF'=B'E+C'F$.

Drawing together much of the work above, we get the following
result, which is the explicit working out of those special
Lagrangian 3-folds of Theorem \ref{ev5thm} coming out of
part (a) of Theorem~\ref{ev9thm}.

\begin{thm} Define functions $w_j,p_j,q_j:\R\ra\C$ for $j=1,2,3$ by
\begin{align*}
w_1(t)&=\tanh t,\qquad w_2(t)=w_3(t)=\sech t,\\
p_1(t)&=B\bigl(f(t)\tanh t-2\sqrt{\cosh t}\,\bigr)
+\frac{iE}{\sqrt{\cosh t}},\\
p_2(t)&=B\,f(t)\sech t+C\,\sqrt{\cosh t}
+\frac{iE\,\sinh t}{2\sqrt{\cosh t}}+\frac{iF}{\sqrt{\cosh t}},\\
p_3(t)&=B\,f(t)\sech t-C\,\sqrt{\cosh t}
+\frac{iE\,\sinh t}{2\sqrt{\cosh t}}-\frac{iF}{\sqrt{\cosh t}},\\
q_1(t)&=B'\bigl(f(t)\tanh t-2\sqrt{\cosh t}\,\bigr)
+\frac{iE'}{\sqrt{\cosh t}},\\
q_2(t)&=B'\,f(t)\sech t+C'\,\sqrt{\cosh t}
+\frac{iE'\,\sinh t}{2\sqrt{\cosh t}}+\frac{iF'}{\sqrt{\cosh t}},\\
q_3(t)&=B'\,f(t)\sech t-C'\,\sqrt{\cosh t}
+\frac{iE'\,\sinh t}{2\sqrt{\cosh t}}-\frac{iF'}{\sqrt{\cosh t}},
\end{align*}
where $B,\ldots,F'\in\R$ satisfy $BE'+CF'=B'E+C'F$, and\/
$f(t)=\int_0^t\sqrt{\cosh s}\,\d s$. Let\/ $r_1,r_2,r_3:\R\ra\C$ 
be the unique solutions of 
\begin{equation*}
\frac{r_1}{\d t}=\ha(\bar p_2\bar p_3+\bar q_3\bar q_2),\quad
\frac{r_2}{\d t}=\ha(\bar q_3\bar q_1-\bar p_1\bar p_3),\quad
\frac{r_3}{\d t}=-\ha(\bar p_1\bar q_2+\bar p_2\bar q_1)
\end{equation*}
with\/ $r_j(0)=0$. Define a subset $N$ of\/ $\C^3$ by 
\begin{align*}
N=\Bigl\{\bigl(&\ha(y_1^2+y_2^2)w_1(t)+y_1p_1(t)+y_2q_1(t)+r_1(t),\\
&\ha(y_1^2-y_2^2)w_2(t)+y_1p_2(t)-y_2q_2(t)+r_2(t),\\
&y_1y_2w_3(t)+y_1q_3(t)+y_2p_3(t)+r_3(t)\bigr):t,y_1,y_2\in\R\Bigr\}\,.
\end{align*}
Then $N$ is a special Lagrangian $3$-fold.
\label{ev10thm2}
\end{thm}

Observe that if $E=E'=CF=C'F'=0$ then $w_1,p_1,q_1$ and $r_1$ are 
real, so that $z_1\in\R$ for each $(z_1,z_2,z_3)\in N$. This implies 
that writing $\C^3=\C\t\C^2$, we have $N\subseteq\R\t\Si$, where
$\Si$ is a special Lagrangian 2-fold in $\C^2$. This is not obvious 
from the explicit expression for $N$, because the coordinates
$t,y_1,y_2$ are not compatible with the product structure on $N$. But
SL 2-folds in $\C^2$ are holomorphic with respect to an alternative
complex structure $J$ on $\C^2$, and in fact $\Si$ is a {\it complex
parabola} with respect to~$J$.

\section{Case (d) of Theorem \ref{ev9thm}}
\label{ev11}

Next we study the SL 3-folds $N$ which come from the construction
of Theorem \ref{ev5thm} when ${\bf z}_1,{\bf z}_2,{\bf z}_3$ are as
in case (d) of Theorem \ref{ev9thm}. Throughout this section we use
the notation of \S\ref{ev9}. Let $\al_1,\al_2,\al_3>0$ satisfy
$\frac{1}{\al_1}=\frac{1}{\al_2}+\frac{1}{\al_3}$, and as in case (d) of
Theorem \ref{ev9thm}, define $a_1,a_2,a_3\in\R$ by
\e
\begin{gathered}
a_1=-\al_1^{-1/2}(\al_2\al_3)^{1/2},\qquad
a_2=\al_2^{-1/2}(\al_3\al_1)^{1/2}\\
\text{and}\qquad a_3=\al_3^{-1/2}(\al_1\al_2)^{1/2},
\end{gathered}
\label{ev11eq1}
\e
so that $a_1+a_2+a_3=0$. Define $w_1,w_2,w_3:\R\ra\C$ by
\e
\begin{gathered}
w_1(t)=i\sqrt{\al_1}\,{\rm e}^{ia_1t},\qquad
w_2(t)=\sqrt{\al_2}\,{\rm e}^{ia_2t}\\
\text{and}\qquad w_3(t)=\sqrt{\al_3}\,{\rm e}^{ia_3t}.
\end{gathered}
\label{ev11eq2}
\e
Then \eq{ev9eq3} holds, and when we define ${\bf z}_1,{\bf z}_2,{\bf z}_3$
by \eq{ev9eq2} they satisfy~\eq{ev5eq10}.

In \S\ref{ev111} we solve equations \eq{ev5eq11} and \eq{ev5eq12} for
${\bf z}_4,{\bf z}_5$ and ${\bf z}_6$, and so write down in Theorem
\ref{ev11thm1} a large family of explicit SL 3-folds in $\C^3$. Section
\ref{ev112} then studies solutions periodic in $t$, which are surprisingly
abundant, interprets them geometrically, and gives a parameter count for
the construction.

\subsection{Solving the equations for ${\bf z}_4,{\bf z}_5$ and ${\bf z}_6$}
\label{ev111}

We shall explicitly solve equations \eq{ev5eq11} and \eq{ev5eq12}
for ${\bf z}_4,{\bf z}_5$ and ${\bf z}_6$ with this choice of
${\bf z}_1,{\bf z}_2,{\bf z}_3$. But \eq{ev5eq11} is equivalent
to equations \eq{ev9eq8} and \eq{ev9eq9} for $p_1,p_2,p_3$ and
$q_1,q_2,q_3$, as in \S\ref{ev92}. As \eq{ev9eq8} and \eq{ev9eq9}
are identical, we need only solve \eq{ev9eq8}. Here is a convenient
way of rewriting it.

\begin{prop} In the situation above, write
\e
p_1=i{\rm e}^{ia_1t}\be_1,\quad
p_2={\rm e}^{ia_2t}\be_2\quad\text{and}\quad
p_3={\rm e}^{ia_3t}\be_3,
\label{ev11eq3}
\e
where $\be_j:\R\ra\C$ is a differentiable function. Then
\eq{ev9eq8} is equivalent to
\begin{equation*}
\frac{\d}{\d t}\begin{pmatrix}\be_1\\ \be_2\\ \be_3\\
\bar\be_1\\ \bar\be_2\\ \bar\be_3\end{pmatrix}
=\ha i \begin{pmatrix}
-2a_1 & 0 & 0 & 0 & -\sqrt{\al_3} & -\sqrt{\al_2} \\
0 & -2a_2 & 0 & \sqrt{\al_3} & 0 & \pha{-}\sqrt{\al_1} \\
0 & 0 & -2a_3 & \sqrt{\al_2} & \pha{-}\sqrt{\al_1} & 0 \\
0 & \pha{-}\sqrt{\al_3} & \pha{-}\sqrt{\al_2} & 2a_1 & 0 & 0 \\
-\sqrt{\al_3} & 0 & -\sqrt{\al_1} & 0 & 2a_2 & 0 \\
-\sqrt{\al_2} & -\sqrt{\al_1} & 0 & 0 & 0 & 2a_3
\end{pmatrix}
\begin{pmatrix}\be_1\\ \be_2\\ \be_3\\
\bar\be_1\\ \bar\be_2\\ \bar\be_3\end{pmatrix}.
\end{equation*}
\label{ev11prop1}
\end{prop}

\begin{proof} Substitute \eq{ev11eq3} into the first equation of
\eq{ev9eq8}, and use the values for $w_1,w_2,w_3$ in \eq{ev11eq2}. We get
\begin{equation*}
\frac{\d}{\d t}\bigl(i{\rm e}^{ia_1t}\be_1\bigr)
=\ha{\rm e}^{-i(a_2+a_3)t}(\sqrt{\al_3}\,\bar\be_2+\sqrt{\al_2}\,\bar\be_3).
\end{equation*}
But ${\rm e}^{-i(a_2+a_3)t}={\rm e}^{ia_1t}$ as $a_1+a_2+a_3=0$. So 
dividing by $i{\rm e}^{ia_1t}$ gives
\begin{equation*}
\frac{\d}{\d t}\be_1+ia_1\be_1=-\ha i(\sqrt{\al_3}\,\bar\be_2
+\sqrt{\al_2}\,\bar\be_3),
\end{equation*}
which is equivalent to the first line of the equation we have to prove.
The other lines follow in the same way.
\end{proof}

To solve for $\be_1,\be_2,\be_3$ we need the eigenvalues and
eigenvectors of this matrix.

\begin{prop} Let\/ $M$ be the $6\t 6$ real matrix appearing in
Proposition \ref{ev11prop1}. Then the eigenvalues of\/ $M$ are $0$
(twice), $\la,-\la,3\la$ and\/ $-3\la$, where $\la>0$ satisfies
\e
\la^2=a_1^2-a_2a_3=a_2^2-a_3a_1=a_3^2-a_1a_2.
\label{ev11eq4}
\e
There exist in $\R^3$ nonzero real vectors $(b_1,b_2,b_3)$, which is 
unique, $(c_1,c_2,c_3)$, $(d_1,d_2,d_3)$, $(e_1,e_2,e_3)$ and\/ 
$(f_1,f_2,f_3)$, such that
\begin{alignat}{2}
M&\begin{pmatrix} \sqrt{\al_1} \\ \sqrt{\al_2} \\ \sqrt{\al_3} \\
\sqrt{\al_1} \\ \sqrt{\al_2} \\ \sqrt{\al_3}\end{pmatrix}=0,\qquad&
M&\begin{pmatrix} b_1 \\ b_2 \\ b_3 \\ -b_1 \\ -b_2 \\ -b_3 \end{pmatrix}
=\begin{pmatrix} \sqrt{\al_1} \\ \sqrt{\al_2} \\ \sqrt{\al_3} \\
\sqrt{\al_1} \\ \sqrt{\al_2} \\ \sqrt{\al_3}\end{pmatrix},
\label{ev11eq5}\\
M&\begin{pmatrix} c_1 \\ c_2 \\ c_3 \\ d_1 \\ d_2 \\ d_3 \end{pmatrix}=
\la\begin{pmatrix} c_1 \\ c_2 \\ c_3 \\ d_1 \\ d_2 \\ d_3 \end{pmatrix},
\qquad&
M&\begin{pmatrix} d_1 \\ d_2 \\ d_3 \\ c_1 \\ c_2 \\ c_3 \end{pmatrix}=
-\la\begin{pmatrix} d_1 \\ d_2 \\ d_3 \\ c_1 \\ c_2 \\ c_3 \end{pmatrix},
\label{ev11eq6}\\
M&\begin{pmatrix} e_1 \\ e_2 \\ e_3 \\ f_1 \\ f_2 \\ f_3 \end{pmatrix}=
3\la\begin{pmatrix} e_1 \\ e_2 \\ e_3 \\ f_1 \\ f_2 \\ f_3 \end{pmatrix}
\quad&\text{and}\quad
M&\begin{pmatrix} f_1 \\ f_2 \\ f_3 \\ e_1 \\ e_2 \\ e_3 \end{pmatrix}=
-3\la\begin{pmatrix} f_1 \\ f_2 \\ f_3 \\ e_1 \\ e_2 \\ e_3\end{pmatrix}.
\label{ev11eq7}
\end{alignat}
\label{ev11prop2}
\end{prop}

\begin{proof} From \eq{ev11eq1} we see that $M=LM'L^{-1}$, where
\begin{align*}
L&=\begin{pmatrix}
\sqrt{\al_1} & 0 & 0 & 0 & 0 & 0 \\
0 & \sqrt{\al_2} & 0 & 0 & 0 & 0 \\
0 & 0 & \sqrt{\al_3} & 0 & 0 & 0 \\
0 & 0 & 0 & \sqrt{\al_1} & 0 & 0 \\
0 & 0 & 0 & 0 & \sqrt{\al_2} & 0 \\
0 & 0 & 0 & 0 & 0 & \sqrt{\al_3} 
\end{pmatrix}\\
\text{and}\qquad
M'&=\begin{pmatrix}
-2a_1 & 0 & 0 & 0 & a_1 & a_1 \\
0 & -2a_2 & 0 & a_2 & 0 & a_2 \\
0 & 0 & -2a_3 & a_3 & a_3 & 0 \\
0 & -a_1 & -a_1 & 2a_1 & 0 & 0 \\
-a_2 & 0 & -a_2 & 0 & 2a_2 & 0 \\
-a_3 & -a_3 & 0 & 0 & 0 & 2a_3
\end{pmatrix}.
\end{align*}
Thus the eigenvalues of $M$ are the same as those of $M'$, which are
the roots of the polynomial $p(x)=\det(M'-xI)$. Direct calculation
shows that
\begin{align*}
p(x)=x^6+&\bigl(2(a_2a_3+a_3a_1+a_1a_2)-4(a_1^2+a_2^2+a_3^2)\bigr)x^4\\
+&\bigl(9(a_2^2a_3^2+a_3^2a_1^2+a_1^2a_2^2)-6(a_1^2a_2a_3+a_2^2a_3a_1
+a_3^2a_1a_2)\bigr)x^2.
\end{align*}
Now $a_1<0$ and $a_2,a_3>0$ with $a_1+a_2+a_3=0$. Thus we may
substitute $a_2=-\si a_1$, $a_3=(\si-1)a_1$, for $\si=-a_2/a_1$
in $(0,1)$. Then $p(x)$ becomes
\begin{equation*}
p(x)=x^6-10a_1^2(\si^2-\si-1)x^4+9a_1^4(\si^4-2\si^3+3\si^2-2\si+1)x^2,
\end{equation*}
which factorizes as 
\begin{equation*}
p(x)=x^2\bigl(x^2-a_1^2(\si^2-\si+1)\bigr)\bigl(x^2-9a_1^2(\si^2-\si+1)\bigr).
\end{equation*}
Thus the roots of $p$ are zero (twice), $\la,-\la,3\la$ and $-3\la$, where 
\begin{equation*}
\la^2=a_1^2(\si^2-\si+1)=a_1^2-(-\si a_1)((\si-1)a_1)=a_1^2-a_2a_3.
\end{equation*}
This proves the first equation of \eq{ev11eq4}. Note also that as
$\si\in(0,1)$ and $a_1<0$, we have $a_1^2(\si^2-\si+1)>0$, so we can take
$\la$ to be real and positive. The second two equations of \eq{ev11eq4}
follow from $a_1+a_2+a_3=0$, as for instance
\begin{equation*}
(a_1^2-a_2a_3)-(a_2^2-a_3a_1)=(a_1+a_2+a_3)(a_1-a_2)=0.
\end{equation*}

The first equation of \eq{ev11eq5} follows from \eq{ev11eq1}. The second is 
equivalent to
\begin{equation*}
-\begin{pmatrix}2a_1 & -\sqrt{\al_3} & -\sqrt{\al_2} \\
\sqrt{\al_3} & 2a_2 & \sqrt{\al_1} \\ \sqrt{\al_2} & \sqrt{\al_1} & 2a_3
\end{pmatrix}\begin{pmatrix} b_1 \\ b_2 \\ b_3 \end{pmatrix}=
\begin{pmatrix} \sqrt{\al_1} \\ \sqrt{\al_2} \\ \sqrt{\al_3} \end{pmatrix},
\end{equation*}
which has a unique solution $(b_1,b_2,b_3)$, as the $3\t 3$ matrix
appearing here has determinant $4a_1a_2a_3\ne 0$, and so is invertible.
To prove \eq{ev11eq6}, observe that as $\la$ is a real eigenvalue of
the real matrix $M$, there exists a real eigenvector 
$(c_1\, c_2\, c_3\, d_1\, d_2\, d_3)^T$ of $M$ with eigenvalue $\la$.
This gives the first equation of \eq{ev11eq6}, and the second then
holds because of the form of $M$. One can also easily show that 
$(c_1,c_2,c_3)$ and $(d_1,d_2,d_3)$ are both nonzero. Equation
\eq{ev11eq7} is proved in the same way, as $3\la$ is an
eigenvalue of~$M$.
\end{proof}

The following identities will be helpful later.

\begin{prop} The constants $b_j,\ldots,f_j$ of Proposition
\ref{ev11prop2} satisfy
\ea
\sqrt{\al_1}\,c_1-\sqrt{\al_2}\,c_2-\sqrt{\al_3}\,c_3
-\sqrt{\al_1}\,d_1+\sqrt{\al_2}\,d_2+\sqrt{\al_3}\,d_3&=0,
\label{ev11eq8}\\
\sqrt{\al_1}\,e_1-\sqrt{\al_2}\,e_2-\sqrt{\al_3}\,e_3
-\sqrt{\al_1}\,f_1+\sqrt{\al_2}\,f_2+\sqrt{\al_3}\,f_3&=0,
\label{ev11eq9}\\
b_1c_1-b_2c_2-b_3c_3+b_1d_1-b_2d_2-b_3d_3&=0,
\label{ev11eq10}\\
b_1e_1-b_2e_2-b_3e_3+b_1f_1-b_2f_2-b_3f_3&=0,
\label{ev11eq11}\\
c_1e_1-c_2e_2-c_3e_3-d_1f_1+d_2f_2+d_3f_3&=0,
\label{ev11eq12}\\
\text{and}\qquad
c_1f_1-c_2f_2-c_3f_3-d_1e_1+d_2e_2+d_3e_3&=0.
\label{ev11eq13}
\ea
Furthermore, none of\/ $\pm\la$ or $\pm 3\la$ is equal to $a_1,a_2$
or~$a_3$.
\label{ev11prop3}
\end{prop}

\begin{proof} The $6\t 6$ matrix $M$ studied above is not symmetric,
because the signs of the $\sqrt{\al_1}$ terms are not right. However,
if we define
\begin{equation*}
K=\begin{pmatrix}
1 & 0 & 0 & 0 & 0 & 0 \\
0 & i & 0 & 0 & 0 & 0 \\
0 & 0 & i & 0 & 0 & 0 \\
0 & 0 & 0 & i & 0 & 0 \\
0 & 0 & 0 & 0 & 1 & 0 \\
0 & 0 & 0 & 0 & 0 & 1 
\end{pmatrix}
\end{equation*}
then $KMK^{-1}$ is a {\it symmetric} complex matrix, and 
Proposition \ref{ev11prop2} shows that 
\begin{equation*}
\begin{pmatrix} \sqrt{\al_1} \\ i\sqrt{\al_2} \\ i\sqrt{\al_3} \\
i\sqrt{\al_1} \\ \sqrt{\al_2} \\ \sqrt{\al_3}\end{pmatrix},\quad
\begin{pmatrix} c_1 \\ ic_2 \\ ic_3 \\ id_1 \\ d_2 \\ d_3\end{pmatrix},\quad
\begin{pmatrix} d_1 \\ id_2 \\ id_3 \\ ic_1 \\ c_2 \\ c_3\end{pmatrix},\quad
\begin{pmatrix} e_1 \\ ie_2 \\ ie_3 \\ if_1 \\ f_2 \\ f_3\end{pmatrix}
\quad\text{and}\quad
\begin{pmatrix} f_1 \\ if_2 \\ if_3 \\ ie_1 \\ e_2 \\ e_3\end{pmatrix}
\end{equation*}
are eigenvectors of $KMK^{-1}$, with eigenvalues $0,\la,-\la,3\la$
and $-3\la$ respectively.

Now if ${\bf u},{\bf v}$ are eigenvectors of a symmetric matrix with
different eigenvalues, then ${\bf u}^T{\bf v}=0$. Taking inner products
between the five eigenvectors above proves \eq{ev11eq8}, \eq{ev11eq9},
\eq{ev11eq12} and \eq{ev11eq13}. Also, $(b_1,ib_2,ib_3,-ib_1,-b_2,-b_3)^T$ is
essentially an eigenvector with eigenvalue 0, modulo nilpotent behaviour,
and taking inner products of this with the $\la$ and $3\la$ eigenvectors
gives \eq{ev11eq10} and~\eq{ev11eq11}.

To prove the last part, as in the proof of Proposition \ref{ev11prop2}
put $a_2=-\si a_1$, $a_3=(\si-1)a_1$ and $\la^2=a_1^2(\si^2-\si+1)$
for some $\si\in(0,1)$. Then $a_1=\pm\la$ when $a_1^2=a_1^2(\si^2-\si+1)$,
that is, when $\si=0$ or 1, which contradicts $\si\in(0,1)$. Similarly,
$a_2=\pm\la$ when $a_1^2\si^2=a_1^2(\si^2-\si+1)$, that is, when $\si=1$,
contradicting $\si\in(0,1)$, and $a_3=\pm\la$ is ruled out in the same way.
To eliminate $a_j=\pm 3\la$ for $j=1,2,3$, note that $\md{a_j}\le -a_1$,
but $\la^2=a_1^2(\si^2-\si+1)$ for $\si\in\R$ implies
that~$3\la\ge -3\sqrt{3}\,a_1/2>-a_1$.
\end{proof}

From Proposition \ref{ev11prop2} we see that the o.d.e.\ in Proposition
\ref{ev11prop1} has solution
\begin{align*}
\begin{pmatrix}\be_1\\ \be_2\\ \be_3\\
\bar\be_1\\ \bar\be_2\\ \bar\be_3\end{pmatrix}=
(A-\ha Bt)\begin{pmatrix} \sqrt{\al_1} \\ \sqrt{\al_2} \\ \sqrt{\al_3} \\
\sqrt{\al_1} \\ \sqrt{\al_2} \\ \sqrt{\al_3} \end{pmatrix}
\!+\!iB\begin{pmatrix} b_1\\b_2\\b_3\\-b_1\\-b_2\\-b_3 \end{pmatrix}
&+\!C{\rm e}^{i\la t/2\pha{3}}
\begin{pmatrix} c_1\\c_2\\c_3\\d_1\\d_2\\d_3 \end{pmatrix}
\!+\!\hat C{\rm e}^{-i\la t/2\pha{3}}
\begin{pmatrix} d_1\\d_2\\d_3\\c_1\\c_2\\c_3 \end{pmatrix}\\
&+\!D{\rm e}^{3i\la t/2}
\begin{pmatrix} e_1\\e_2\\e_3\\f_1\\f_2\\f_3\end{pmatrix}
\!+\!\hat D{\rm e}^{-3i\la t/2}
\begin{pmatrix} f_1\\f_2\\f_3\\e_1\\e_2\\e_3 \end{pmatrix}
\end{align*}
for $A,B,C,\hat C,D,\hat D\in\C$. But the last three rows of this equation
should be the complex conjugates of the first three rows, which implies
that $A$ and $B$ are real and $\hat C=\bar C$, $\hat D=\bar D$. This gives
$\be_1,\be_2,\be_3$, and then \eq{ev11eq3} gives the general solution for
$p_j$. We have proved:

\begin{prop} In the situation above, the general solution of\/
\eq{ev9eq8} is
\ea
\begin{split}
p_1(t)&=i\bigl((A-\ha Bt)\sqrt{\al_1}\,{\rm e}^{ia_1t}
+iBb_1{\rm e}^{ia_1t}+Cc_1{\rm e}^{i(a_1+\la/2)t}\\
&\quad +\bar Cd_1{\rm e}^{i(a_1-\la/2)t}+De_1{\rm e}^{i(a_1+3\la/2)t}
+\bar Df_1{\rm e}^{i(a_1-3\la/2)t}\bigr),
\end{split}
\label{ev11eq14}\\
\begin{split}
p_2(t)&=(A-\ha Bt)\sqrt{\al_2}\,{\rm e}^{ia_2t}
+iBb_2{\rm e}^{ia_2t}+Cc_2{\rm e}^{i(a_2+\la/2)t}\\
&\quad +\bar Cd_2{\rm e}^{i(a_2-\la/2)t}+De_2{\rm e}^{i(a_2+3\la/2)t}
+\bar Df_2{\rm e}^{i(a_2-3\la/2)t},
\end{split}
\label{ev11eq15}\\
\begin{split}
p_3(t)&=(A-\ha Bt)\sqrt{\al_3}\,{\rm e}^{ia_3t}
+iBb_3{\rm e}^{ia_3t}+Cc_3{\rm e}^{i(a_3+\la/2)t}\\
&\quad +\bar Cd_3{\rm e}^{i(a_3-\la/2)t}+De_3{\rm e}^{i(a_3+3\la/2)t}
+\bar Df_3{\rm e}^{i(a_3-3\la/2)t},
\end{split}
\label{ev11eq16}
\ea
for $A,B\in\R$ and\/ $C,D\in\C$, where $b_j,\ldots,f_j$ and\/ $\la$
are as in Proposition~\ref{ev11prop2}.
\label{ev11prop4}
\end{prop}

Define $p_1,p_2,p_3$ by \eq{ev11eq14}--\eq{ev11eq16}, and similarly 
define $q_1,q_2,q_3$ by \eq{ev11eq14}--\eq{ev11eq16}, but using constants
$A',B'\in\R$ and $C',D'\in\C$ instead of $A,B,C,D$. Then $p_j$ and $q_j$
satisfy \eq{ev9eq8} and \eq{ev9eq9}. From \S\ref{ev92}, if we define
${\bf z}_4,{\bf z}_5$ by \eq{ev9eq7}, then they satisfy \eq{ev5eq11}.
We have found the general solution of the evolution equation 
\eq{ev5eq11} for ${\bf z}_4$ and ${\bf z}_5$, with the given values 
of~${\bf z}_1,{\bf z}_2,{\bf z}_3$.

Now from \S\ref{ev5}, we need ${\bf z}_4$ and ${\bf z}_5$ to satisfy
\eq{ev5eq5}--\eq{ev5eq7}. By Lemma \ref{ev9lem1}, these are equivalent
to equation \eq{ev9eq11}. Using the values of $w_j$ and $p_j$ above,
we find that
\begin{align*}
\Im(w_1\bar p_1-w_2\bar p_2-w_3\bar p_3)=&
-B(\sqrt{\al_1}\,b_1-\sqrt{\al_2}\,b_2-\sqrt{\al_3}\,b_3)\\
-\Im(C{\rm e}^{i\la t/2})
(\sqrt{\al_1}\,c_1&-\sqrt{\al_2}\,c_2-\sqrt{\al_3}\,c_3
-\sqrt{\al_1}\,d_1+\sqrt{\al_2}\,d_2+\sqrt{\al_3}\,d_3)\\
-\Im(D{\rm e}^{3i\la t/2})
(\sqrt{\al_1}\,e_1&-\sqrt{\al_2}\,e_2-\sqrt{\al_3}\,e_3
-\sqrt{\al_1}\,f_1+\sqrt{\al_2}\,f_2+\sqrt{\al_3}\,f_3)\\
=&-B(\sqrt{\al_1}\,b_1-\sqrt{\al_2}\,b_2-\sqrt{\al_3}\,b_3),
\end{align*}
where in the last line we have used \eq{ev11eq8} and \eq{ev11eq9}. Since 
in general $\sqrt{\al_1}\,b_1-\sqrt{\al_2}\,b_2-\sqrt{\al_3}\,b_3\ne 0$,
we see that $\Im(w_1\bar p_1-w_2\bar p_2-w_3\bar p_3)=0$ if and only 
if~$B=0$.

Similarly, using \eq{ev11eq8}--\eq{ev11eq13} we find that
$\Im(w_1\bar q_1-w_2\bar q_2-w_3\bar q_3)=0$ if and only if $B'=0$, 
and $\Im(p_1\bar q_1-p_2\bar q_2-p_3\bar q_3)=0$ if and only if
\begin{align*}
&(A'B-AB')(\sqrt{\al_1}\,b_1-\sqrt{\al_2}\,b_2-\sqrt{\al_3}\,b_3)+\\
&\Im(C\bar C')(c_1^2\!-\!c_2^2\!-\!c_3^2\!-\!d_1^2\!+\!d_2^2\!+\!d_3^2)
\!+\!\Im(D\bar D')(e_1^2\!-\!e_2^2\!-\!e_3^2\!-\!f_1^2\!+\!f_2^2\!+\!f_3^2)=0.
\end{align*}
Thus we should set $B=B'=0$. But we can show by changing coordinates 
in $\R^2$ from $(y_1,y_2)$ to $(y_1+A,y_2+A')$ as in \S\ref{ev51} that 
we are also free to set $A$ and $A'$ to be zero, without restricting
the SL 3-folds constructed in Theorem~\ref{ev5thm}. 

Next we solve equation \eq{ev5eq12} for ${\bf z}_6$. From \S\ref{ev92},
this is equivalent to the o.d.e.s \eq{ev9eq10} for $r_1,r_2,r_3$, where
${\bf z}_6=(r_1,r_2,r_3)$. Using the expressions above for $p_j$ and
$q_j$ and remembering that $A=A'=B=B'=0$, we get rather complicated
expressions for $\d r_j/\d t$ for $j=1,2,3$, which can then be
integrated to get $r_1,r_2$ and~$r_3$.

We sum up all the above work in the following result, which is the
explicit working out of those special Lagrangian 3-folds of Theorem
\ref{ev5thm} coming out of part (d) of Theorem~\ref{ev9thm}.

\begin{thm} Define functions $w_j,p_j,q_j,r_j:\R\ra\C$ for $j=1,2,3$ by
\begin{align*}
w_1(t)=&i\sqrt{\al_1}\,{\rm e}^{ia_1t},\quad
w_2(t)=\sqrt{\al_2}\,{\rm e}^{ia_2t},\quad
w_3(t)=\sqrt{\al_3}\,{\rm e}^{ia_3t},
\displaybreak[0]\\
p_1(t)=&i\bigl(Cc_1{\rm e}^{i(a_1+\la/2)t}
\!\!+\!\bar Cd_1{\rm e}^{i(a_1-\la/2)t}
\!\!+\!De_1{\rm e}^{i(a_1+3\la/2)t}
\!\!+\!\bar Df_1{\rm e}^{i(a_1-3\la/2)t}\bigr),\\
p_2(t)=&Cc_2{\rm e}^{i(a_2+\la/2)t}\!\!
+\!\bar Cd_2{\rm e}^{i(a_2-\la/2)t}
\!\!+\!De_2{\rm e}^{i(a_2+3\la/2)t}\!\!
+\!\bar Df_2{\rm e}^{i(a_2-3\la/2)t},\\
p_3(t)=&Cc_3{\rm e}^{i(a_3+\la/2)t}\!\!
+\!\bar Cd_3{\rm e}^{i(a_3-\la/2)t}
\!\!+\!De_3{\rm e}^{i(a_3+3\la/2)t}\!\!
+\!\bar Df_3{\rm e}^{i(a_3-3\la/2)t},
\displaybreak[0]\\
q_1(t)=&i\bigl(C'c_1{\rm e}^{i(a_1+\la/2)t}
\!\!+\!\bar C'd_1{\rm e}^{i(a_1-\la/2)t}
\!\!+\!D'e_1{\rm e}^{i(a_1+3\la/2)t}
\!\!+\!\bar D'f_1{\rm e}^{i(a_1-3\la/2)t}\bigr),\\
q_2(t)=&C'c_2{\rm e}^{i(a_2+\la/2)t}\!\!
+\!\bar C'd_2{\rm e}^{i(a_2-\la/2)t}
\!\!+\!D'e_2{\rm e}^{i(a_2+3\la/2)t}\!\!
+\!\bar D'f_2{\rm e}^{i(a_2-3\la/2)t},\\
q_3(t)=&C'c_3{\rm e}^{i(a_3+\la/2)t}\!\!
+\!\bar C'd_3{\rm e}^{i(a_3-\la/2)t}
\!\!+\!D'e_3{\rm e}^{i(a_3+3\la/2)t}\!\!
+\!\bar D'f_3{\rm e}^{i(a_3-3\la/2)t},
\displaybreak[0]\\
r_1(t)=&-\!{\ts\frac{i}{2(a_1+3\la)}}(D^2\!+\!D'{}^2)f_2f_3
{\rm e}^{i(a_1+3\la)t}\!-\!{\ts\frac{i}{2(a_1-3\la)}}
(\bar D^2\!+\!\bar D'{}^2)e_2e_3{\rm e}^{i(a_1-3\la)t}\\
&-{\ts\frac{i}{2(a_1+\la)}}\bigl((C^2+C'{}^2)d_2d_3+(\bar CD+\bar C'D')
(c_2f_3+f_2c_3)\bigr){\rm e}^{i(a_1+\la)t}\\
&-{\ts\frac{i}{2(a_1-\la)}}\bigl((\bar C^2\!+\!\bar C'{}^2)c_2c_3
\!+\!(C\bar D\!+\!C'\bar D')(d_2e_3\!+\!e_2d_3)\bigr)
{\rm e}^{i(a_1-\la)t}\!+\!E_1,
\displaybreak[0]\\
r_2(t)=&{\ts\frac{1}{2(a_2+3\la)}}(D^2-D'{}^2)f_1f_3{\rm e}^{i(a_2+3\la)t}
+{\ts\frac{1}{2(a_2-3\la)}}(\bar D^2-\bar D'{}^2)e_1e_3
{\rm e}^{i(a_2-3\la)t}\\
&+{\ts\frac{1}{2(a_2+\la)}}\bigl((C^2-C'{}^2)d_1d_3+(\bar CD-\bar C'D')
(c_1f_3+f_1c_3)\bigr){\rm e}^{i(a_2+\la)t}\\
&+{\ts\frac{1}{2(a_2-\la)}}\bigl((\bar C^2\!-\!\bar C'{}^2)c_1c_3
\!+\!(C\bar D\!-\!C'\bar D')(d_1e_3\!+\!e_1d_3)\bigr)
{\rm e}^{i(a_2-\la)t}\!+\!E_2,
\displaybreak[0]\\
r_3(t)=&{\ts\frac{1}{a_3+3\la}}DD'f_1f_2{\rm e}^{i(a_3+3\la)t}
+{\ts\frac{1}{a_3-3\la}}\bar D\bar D'{\rm e}^{i(a_3-3\la)t}\\
&+{\ts\frac{1}{2(a_3+\la)}}\bigl(2CC'd_1d_2+(\bar CD'+\bar C'D)
(c_1f_2+f_1c_2)\bigr){\rm e}^{i(a_3+\la)t}\\
&+{\ts\frac{1}{2(a_3-\la)}}\bigl(2\bar C\bar C'c_1c_2+(C\bar D'+C'\bar D)
(d_1e_2+e_1d_2)\bigr){\rm e}^{i(a_3-\la)t}+E_3,
\end{align*}
where $E_1,E_2,E_3\in\C$ and\/ $C,D,C',D'\in\C$ satisfy
\e
\begin{split}
&\Im(C\bar C')(c_1^2-c_2^2-c_3^2-d_1^2+d_2^2+d_3^2)\\
+&\Im(D\bar D')(e_1^2-e_2^2-e_3^2-f_1^2+f_2^2+f_3^2)=0.
\end{split}
\label{ev11eq17}
\e
Define a subset $N$ of\/ $\C^3$ by
\e
\begin{split}
N=\Bigl\{\bigl(&\ha(y_1^2+y_2^2)w_1(t)+y_1p_1(t)+y_2q_1(t)+r_1(t),\\
&\ha(y_1^2-y_2^2)w_2(t)+y_1p_2(t)-y_2q_2(t)+r_2(t),\\
&y_1y_2w_3(t)+y_1q_3(t)+y_2p_3(t)+r_3(t)\bigr):y_1,y_2,t\in\R\Bigr\}\,.
\end{split}
\label{ev11eq18}
\e
Then $N$ is a special Lagrangian $3$-fold.
\label{ev11thm1}
\end{thm}

Note that in the expressions for $r_1,r_2$ and $r_3$ we have divided by
factors $a_j\pm\la$ and $a_j\pm 3\la$ for $j=1,2,3$. This is legitimate
because none of $\pm\la$ or $\pm 3\la$ is equal to $a_1,a_2$ or $a_3$
by Proposition \ref{ev11prop3}, so none of these factors vanish. If one of
the factors had been zero, we would have had to replace the corresponding
term by a multiple of~$t$.

Observe that all the functions in the theorem are linear combinations
of exponentials ${\rm e}^{i\al t}$ for $\al\in\R$. It would seem a
reasonable guess that this is because the SL 3-fold $N$ is actually
symmetric under a subgroup $\U(1)$ or $\R$ in $\SU(3)$, which acts
by multiplication by such exponentials in suitable coordinates.
However, this is {\it not}\/ the case, and generically the SL 3-folds
of Theorem \ref{ev11thm1} have only discrete symmetry groups.

\subsection{Periodicity conditions}
\label{ev112}

We now discuss {\it periodicity} in $t$ of the SL 3-folds $N$ of
Theorem \ref{ev11thm1}. Let $\Phi:\R^3\ra\C^3$ be defined as in
\eq{ev5eq21}, so that $\Phi(y_1,y_2,t)$ is the vector in \eq{ev11eq18},
and $N=\Image\Phi$. We want to know when $\Phi(y_1,y_2,t+T)=\Phi(y_1,y_2,t)$
for some $T>0$ and all $y_1,y_2,t\in\R$. The corresponding immersed SL
3-folds $N$ will then have topology ${\cal S}^1\t\R^2$ rather than~$\R^3$.

As in the proof of Proposition \ref{ev11prop2}, write $a_2=-\si a_1$,
$a_3=(\si-1)a_1$ and $\la=-a_1\sqrt{\si^2-\si+1}$ for some $\si\in(0,1)$.
Also let $\tau=\sqrt{\si^2-\si+1}$, so that $\tau>0$ and $\la=-a_1\tau$.
The exponentials in the expressions for $w_j,p_j,q_j,r_j$ have the 27
periods
\e
\frac{2\pi}{a_j},\;\> \frac{2\pi}{a_j\pm\la/2},\;\>
\frac{2\pi}{a_j\pm\la},\;\> \frac{2\pi}{a_j\pm 3\la/2}\;\>\text{and}\;\>
\frac{2\pi}{a_j\pm 3\la}\;\>\text{for $j=1,2,3$.}
\label{ev11eq19}
\e
For generic values of $C,D,C',D'$, it is clear that $\Phi$ will be
periodic if and only if these periods have a common multiple. But this
holds exactly when $\si$ and $\tau$ lie in $\Q$. So we need to
understand the set of pairs of rational numbers $(\si,\tau)$ satisfying
$\tau^2=\si^2-\si+1$, with $\si\in(0,1)$ and~$\tau>0$.

Now finding the rational points on a conic is a well-known problem
in elementary number theory, and has a standard solution method,
which is to parametrize the conic in the usual way. A suitable
parametrization of the conic $\tau^2=\si^2-\si+1$ is
\begin{equation*}
\si=(1-2s)/(1-s^2)\qquad\text{and}\qquad \tau=(1-s+s^2)/(1-s^2).
\end{equation*}
This has the property that $\si$ and $\tau$ are both rational if
and only if $s$ is rational, and that $\si\in(0,1)$ and $\tau>0$
if and only if~$s\in(0,\ha)$.

For example, when $s=p/q$ for coprime $p,q\in\Z$ with $0<2p<q$,
we may~set
\e
\begin{split}
a_1=p^2-q^2,\quad a_2&=q^2-2pq,\quad a_3=2pq-p^2 \\
\text{and}\quad \la&=p^2-pq+q^2.
\end{split}
\label{ev11eq20}
\e
Then $a_1,a_2,a_3$ and $\la$ satisfy the equations above with this
value of $s$. Careful consideration shows that
\begin{equation*}
\hcf(a_1,a_2,a_3)=\hcf(a_1,a_2,a_3,\la)=
\begin{cases} 1 &\quad p+q\not\equiv 0\;\mod 3,\\
3 &\quad p+q\equiv 0\;\mod 3.\end{cases}
\end{equation*}
When $p+q\equiv 0\;\mod 3$, we replace the $a_j$ and $\la$ by
\e
\begin{split}
a_1={\ts\frac{1}{3}}(p^2-q^2),\quad a_2&={\ts\frac{1}{3}}(q^2-2pq),\quad
a_3={\ts\frac{1}{3}}(2pq-p^2)\\
\text{and}\quad \la&={\ts\frac{1}{3}}(p^2-pq+q^2),
\end{split}
\label{ev11eq21}
\e
so that in both cases we have $\hcf(a_1,a_2,a_3)=\hcf(a_1,a_2,a_3,\la)=1$.
Having chosen $a_1,a_2,a_3$, we determine $\al_1,\al_2,\al_3$ uniquely
by inverting~\eq{ev11eq1}.

As $p,q$ are coprime at least one of them is odd, so $\la$ is odd.
Thus $a_j,a_j\pm\la$ and $a_j\pm 3\la$ are integers, and $a_j\pm\la/2$
and $a_j\pm 3\la/2$ are half-integers but not integers. It follows that
\begin{equation*}
w_j(t\!+\!2\pi)\!=\!w_j(t),\;\> p_j(t\!+\!2\pi)\!=\!-p_j(t),\;\>
q_j(t\!+\!2\pi)\!=\!-q_j(t),\;\> r_j(t\!+\!2\pi)\!=\!r_j(t).
\end{equation*}
Hence $\Phi(y_1,y_2,t+2\pi)=\Phi(-y_1,-y_2,t)$ for all $y_1,y_2,t\in\R$.
This implies that $\Phi(y_1,y_2,t+4\pi)=\Phi(y_1,y_2,t)$. Moreover, as
$\hcf(a_1,a_2,a_3)=1$ one can show that $4\pi$ is the least $T>0$ with
$\Phi(y_1,y_2,t+T)=\Phi(y_1,y_2,t)$ for all $y_1,y_2,t\in\R$, so that
$\Phi$ is periodic in $t$ with period~$4\pi$.

Therefore we can regard $\Phi$ as a map $\R^3/\Z\ra\C^3$, where $\Z$
acts on $\R^3$ by
\e
(y_1,y_2,y_3){\buildrel n\over\longmapsto}
\bigl((-1)^ny_1,(-1)^ny_2,t+2\pi n\bigr)\qquad\text{for $n\in\Z$.}
\label{ev11eq22}
\e
As $\R^3/\Z$ is diffeomorphic to ${\cal S}^1\t\R^2$, though not with
the obvious $\R^2$ coordinates, we see that if $\Phi$ is an immersion
then $N$ is an immersed 3-submanifold diffeomorphic to~${\cal S}^1\t\R^2$.

Consider the asymptotic behaviour of these 3-folds $N$ at infinity in
$\C^3$. From \eq{ev11eq18} we see that $\Phi(y_1,y_2,t)\ra\iy$ in $\C^3$
as $(y_1,y_2)\ra\iy$ in $\R^2$. But when $(y_1,y_2)$ is large then the
dominant terms in \eq{ev11eq18} are the quadratic terms in $y_1$ and $y_2$,
and we can neglect the lower order terms. Thus we expect that $N$ should
be asymptotic to $N_0=\Phi_0(\R^3)$ at infinity in $\C^3$, to leading
order, where $\Phi_0:\R^3\ra\C^3$ is given by
\begin{equation*}
\Phi_0(y_1,y_2,t)=
\bigl({\ts\frac{i}{2}}(y_1^2+y_2^2)\sqrt{\al_1}\,
{\rm e}^{ia_1t},\ha(y_1^2-y_2^2)\sqrt{\al_2}\,{\rm e}^{ia_2t},
y_1y_2\sqrt{\al_3}\,{\rm e}^{ia_3t}\bigr).
\end{equation*}

Calculation shows that $N_0$ is an SL $T^2$-cone in $\C^3$, which may
be written
\begin{equation*}
\bigl\{(i{\rm e}^{ia_1t}x_1,{\rm e}^{ia_2t}x_2,{\rm e}^{ia_3t}x_3):
x_j,t\in\R,\quad x_1\ge 0,\quad a_1x_1^2+a_2x_2^2+a_3x_3^2=0\bigr\}.
\end{equation*}
We have already met this example in \cite{Joyc2}, in particular
in \cite[\S 7, case (a)]{Joyc2}, \cite[Th.~8.7]{Joyc2}
and~\cite[Ex.~9.5]{Joyc2}.

Regarding $\Phi_0$ as mapping $\R^3/\Z\ra N_0$, where $\Z$ acts on
$\R^3$ by \eq{ev11eq22}, we find that $\Phi_0$ is generically 2:1,
since $\Phi_0(y_1,y_2,t)=\Phi_0(-y_1,-y_2,t)$. So we should think
of $N$ as converging at infinity to a {\it double cover} of the
$T^2$-cone $N_0$. That is, towards infinity two points of $N$
converge to each point of $N_0$, and at infinity $N$ is a
$T^2$-cone which is wrapped twice round the $T^2$-cone~$N_0$.

Note that the convergence of $N$ to $N_0$ is of a rather weak kind.
Let $r$ be the radius function on $\C^3$, so that $y_1^2+y_2^2=O(r)$.
The largest terms in \eq{ev11eq18} we have neglected are linear in
$y_1,y_2$, and so they are $O(r^{1/2})$. Thus $N$ `converges' to $N_0$
to order $O(r^{1/2})$ for large $r$, so that $N$ actually gets further
away from $N_0$ towards infinity, rather than closer.

We summarize the material above in the following theorem.

\begin{thm} For each\/ $s\in(0,\ha)\cap\Q$ the construction of Theorem
\ref{ev11thm1} yields a family of closed special Lagrangian $3$-folds
in $\C^3$ depending on $E_1,E_2,E_3$ and\/ $C,D,C',D'\in\C$ satisfying
\eq{ev11eq17}. Generic members of the family are nonsingular immersed\/
$3$-submanifolds diffeomorphic to ${\cal S}^1\t\R^2$. Each\/ $3$-fold
is weakly asymptotic to order $O(r^{1/2})$ at infinity in $\C^3$ to a
double cover of the special Lagrangian $T^2$-cone
\begin{equation*}
\bigl\{(i{\rm e}^{ia_1t}x_1,{\rm e}^{ia_2t}x_2,{\rm e}^{ia_3t}x_3):
x_j,t\in\R,\quad x_1\ge 0,\quad a_1x_1^2+a_2x_2^2+a_3x_3^2=0\bigr\},
\end{equation*}
where $a_1,a_2,a_3\in\Z$ depend on $s=p/q$ as in \eq{ev11eq20}
or~\eq{ev11eq21}.
\label{ev11thm2}
\end{thm} 

It seems rather odd to the author that the periodicity conditions
in this problem turn out to have such a neat, and geometrically
interesting, answer. We saw in \eq{ev11eq19} that for a general
member of the family to be periodic we need 27 periods to be
relatively rational. As these periods depend only on $a_1,a_2$
and $a_3$, one would expect this to be a very overdetermined
problem, with no interesting solutions, but in fact there are
infinitely many.

It is also surprising that once the rationality conditions are solved,
when $a_1,a_2,a_3$ are integers it turns out that $\la$ is necessarily
an integer, rather than just a rational number. This has the effect
that at infinity $N$ is a double cover of $N_0$, rather than a
multiple cover of some high degree.

We finish by giving a parameter count for the SL 3-folds of Theorem
\ref{ev11thm1}. They depend upon parameters $\al_1,\al_2,\al_3\in\R$
and $C,D,C',D',E_1,E_2,E_3\in\C$, which is 17 real parameters. These
must satisfy $1/\al_1=1/\al_2+1/\al_3$ and \eq{ev11eq17}, reducing it
to 15 parameters. Of the symmetry groups $\GL(2,\R)\lt\R^2$ and
$\SU(3)\lt\C^3$ we have used all but the dilations in $\GL(2,\R)$
and the translations in $\C^3$. We must subtract 7 parameters for
these, leaving 8 parameters.

There is one other symmetry to take account of, which is translation
in time, $t\mapsto t+c$. This has the following effect: if we replace
$C,D,C',D'$ by ${\rm e}^{i\la c/2}C$, ${\rm e}^{3i\la c/2}D$,
${\rm e}^{i\la c/2}C'$ and ${\rm e}^{3i\la c/2}D'$ respectively,
then the corresponding SL 3-fold $N'$ is equivalent to $N$ under
an $\SU(3)\lt\C^3$ transformation. So subtracting one parameter, we
see that the family of SL 3-folds from Theorem \ref{ev11thm1} up to
automorphisms of $\C^3$ has dimension 7. For comparison, the whole
family from Theorem \ref{ev5thm} has dimension~9.

\section{The family of SL 3-folds from Example \ref{ev4ex2}}
\label{ev12}

We now apply the `evolution equation' construction of \S\ref{ev3} 
to the set of affine evolution data defined in Example \ref{ev4ex2}.
The material of this section runs parallel to sections \ref{ev5} and 
\ref{ev8}, and so we will leave out many of the details.

\subsection{Application of the method of \S\ref{ev3}}
\label{ev121}

As in Example \ref{ev4ex2}, let $k\ge 1$, and let $(x_1,\ldots,x_k,y_1,y_2)$ 
be coordinates on $\R^{k+2}$. Define $P$ to be the image in $\R^{k+2}$ of 
the map $\psi:\R^2\ra\R^{k+2}$ given by
\e
\psi:(x,y)\longmapsto (x,x^2,\ldots,x^k,y,xy),
\label{ev12eq1}
\e
and $\chi:\R^{k+2}\ra\La^2\R^{k+2}$ to be the affine map
\e
\begin{split}
\chi(&x_1,\ldots,x_k,y_1,y_2)=\\
&-2y_1{\ts\frac{\pd}{\pd y_1}}\!\w\!{\ts\frac{\pd}{\pd y_2}}
\!+\!2{\ts\frac{\pd}{\pd x_1}}\!\w\!{\ts\frac{\pd}{\pd y_1}}
\!+\!4x_1{\ts\frac{\pd}{\pd x_2}}\!\w\!{\ts\frac{\pd}{\pd y_1}}\!+\!\cdots
\!+\!2kx_{k-1}{\ts\frac{\pd}{\pd x_k}}\!\w\!{\ts\frac{\pd}{\pd y_1}}\\
&+\!2x_1{\ts\frac{\pd}{\pd x_1}}\!\w\!{\ts\frac{\pd}{\pd y_2}}
\!+\!4x_2{\ts\frac{\pd}{\pd x_2}}\!\w\!{\ts\frac{\pd}{\pd y_2}}\!+\!\cdots
\!+\!2kx_k{\ts\frac{\pd}{\pd x_k}}\!\w\!{\ts\frac{\pd}{\pd y_2}}.
\end{split}
\label{ev12eq2}
\e
Then $(P,\chi)$ is a set of affine evolution data with $m=3$ and~$n=k+2$.

Let ${\bf p}_0,\ldots,{\bf p}_k$ and ${\bf q}_1,{\bf q}_2$ be vectors 
in $\C^3$, and define an affine map $\phi:\R^{k+2}\ra\C^3$ by
\e
\phi:(x_1,\ldots,x_k,y_1,y_2)\mapsto 
{\bf p}_0+x_1{\bf p}_1+\cdots+x_k{\bf p}_k
+y_1{\bf q}_1+y_2{\bf q}_2.
\label{ev12eq3}
\e
Then from \eq{ev12eq2} we see that 
\begin{align*}
\phi^*(\om)\cdot\chi&=2\bigl(-y_1\om({\bf q}_1,{\bf q}_2)
\!+\!\om({\bf p}_1,{\bf q}_1)\!+\!2x_1\om({\bf p}_2,{\bf q}_1)
\!+\!\cdots\!+\!kx_{k-1}\om({\bf p}_k,{\bf q}_1)\\
&\quad+x_1\om({\bf p}_1,{\bf q}_2)\!+\!2x_2\om({\bf p}_2,{\bf q}_2)
\!+\!\cdots\!+\!kx_k\om({\bf p}_k,{\bf q}_2)\bigr).
\end{align*}

Thus $\phi^*(\om)\vert_P\equiv 0$ if and only if
\begin{gather}
\om({\bf q}_1,{\bf q}_2)=0,\qquad
\quad \om({\bf p}_k,{\bf q}_2)=0 \qquad \text{and}
\label{ev12eq4}\\
j\,\om({\bf p}_j,{\bf q}_1)+(j-1)\om({\bf p}_{j-1},{\bf q}_2)=0
\quad\text{for $1\le j\le k$.}
\label{ev12eq5}
\end{gather}
These are the analogues of equations~\eq{ev5eq4}--\eq{ev5eq7}.

Let ${\bf p}_0(t),\ldots,{\bf p}_k(t)$ and ${\bf q}_1(t),{\bf q}_2(t)$
be differentiable functions $\R\ra\C^3$, and define $\phi_t$ by 
\eq{ev12eq3} for $t\in\R$. Then as in \S\ref{ev5}, we use \eq{ev12eq2}
to show that equation \eq{ev3eq1} of \S\ref{ev3} holds for the family 
$\bigl\{\phi_t:t\in\R\bigr\}$ if and only if
\begin{align*}
\frac{\d\phi_t}{\d t}(x_1,\ldots,y_2)=&-2y_1\,{\bf q}_1\!\t\!{\bf q}_2
\!+\!2\,{\bf p}_1\!\t\!{\bf q}_1\!+\!4x_1\,{\bf p}_2\!\t\!{\bf q}_1
\!+\!\cdots\!+\!2kx_{k-1}\,{\bf p}_k\!\t\!{\bf q}_1\\
&+2x_1\,{\bf p}_1\t{\bf q}_2+4x_2\,{\bf p}_2\t{\bf q}_2
+\cdots+2k\,x_k{\bf p}_k\t{\bf q}_2,
\end{align*}
where the cross product `$\t$' is defined in \eq{ev5eq9}. Using 
\eq{ev12eq3} we get expressions for $\d{\bf p}_i/\d t$ and
$\d{\bf q}_j/\d t$. So applying Theorem \ref{ev3thm1}, we prove the
following analogue of Theorem \ref{ev5thm}:

\begin{thm} Let\/ $k\ge 1$, and suppose ${\bf p}_0,\ldots,{\bf p_k},
{\bf q}_1,{\bf q}_2:\R\ra\C^3$ are differentiable functions satisfying 
equations \eq{ev12eq4} and\/ \eq{ev12eq5} at\/ $t=0$ and
\begin{gather}
\frac{\d{\bf p}_0}{\d t}=2\,{\bf p}_1\!\t {\bf q}_1, \qquad
\frac{\d{\bf p}_k}{\d t}=2k\,{\bf p}_k\!\t {\bf q}_2,
\label{ev12eq6}\\
\frac{\d{\bf p}_j}{\d t}=2(j+1)\,{\bf p}_{j+1}\!\t {\bf q}_1
+2j\,{\bf p}_j\!\t {\bf q}_2\quad\text{for $1\le j\le k-1$,}
\label{ev12eq7}\\
\frac{\d{\bf q}_1}{\d t}=-2\,{\bf q}_1\!\t {\bf q}_2
\qquad\text{and}\qquad \frac{\d{\bf q}_2}{\d t}=0
\label{ev12eq8}
\end{gather}
for all\/ $t\in\R$, where `$\t$\!' is as in \eq{ev5eq9}. Define
a subset\/ $N$ of\/ $\C^3$ to be
\e
\bigl\{{\bf p}_0(t)\!+\!x\,{\bf p}_1(t)\!+\!\cdots\!+\!x^k\,{\bf p}_k(t)
\!+\!y\,{\bf q}_1(t)\!+\!xy\,{\bf q}_2(t):x,y,t\in\R\bigr\}.
\label{ev12eq9}
\e
Then $N$ is a special Lagrangian $3$-fold in $\C^3$ wherever
it is nonsingular.
\label{ev12thm1}
\end{thm}

As in \S\ref{ev5}, if \eq{ev12eq4} and \eq{ev12eq5} hold at $t=0$ 
then they hold for all $t\in\R$, and given initial values 
${\bf p}_j(0),{\bf q}_j(0)$, there exist unique solutions 
${\bf p}_j(t),{\bf q}_k(t)$ to \eq{ev12eq6}--\eq{ev12eq8} for $t$ 
in $(-\ep,\ep)$ and some small $\ep>0$. In fact solutions always 
exist for all $t\in\R$, and this is why we have used $t\in\R$ 
rather than $t\in(-\ep,\ep)$ above.

Following \eq{ev5eq21}, define $\Phi:\R^3\ra\C^3$ by
\e
\Phi(x,y,t)={\bf p}_0(t)+x\,{\bf p}_1(t)+\cdots+x^k\,{\bf p}_k(t)
+y\,{\bf q}_1(t)+xy\,{\bf q}_2(t),
\label{ev12eq10}
\e
so that $N=\Image\Phi$. As in \S\ref{ev6}, one can show that $\Phi$ fails
to be an immersion for a subset of real codimension one in the set of 
all initial data ${\bf p}_i(0),{\bf q}_j(0)$ satisfying \eq{ev12eq4} 
and \eq{ev12eq5}. In particular, for {\it generic} initial data $\Phi$ 
is an immersion, and $N$ a nonsingular immersed 3-submanifold 
diffeomorphic to~$\R^3$.

Observe also that the SL 3-folds of Theorem \ref{ev12thm1} are ruled 
by straight lines, as we showed in \S\ref{ev83} for the 3-folds of
Theorem \ref{ev8thm2}. As \eq{ev12eq10} contains no terms in $y^2$, 
for each fixed $x,t\in\R$ the set $\bigl\{\Phi(x,y,t):y\in\R\bigr\}$ 
is a real straight line in $\C^3$. So $N$ is fibred by straight lines, 
and is a {\it ruled submanifold}. Ruled SL 3-folds are studied 
in~\cite{Joyc5}.

We have already met the families of SL 3-folds constructed above
when $k=1$ and 2. When $k=1$ we have
\begin{equation*}
\Phi(x,y,t)={\bf p}_0(t)+x\,{\bf p}_1(t)+y\,{\bf q}_1(t)+xy\,{\bf q}_2(t).
\end{equation*}
Clearly, $N=\Image\Phi$ is fibred by the images of the quadric $x_3=x_1x_2$
in $\R^3$ under the affine maps $\phi_t:\R^3\ra\C^3$ given by
\begin{equation*}
\phi_t:(x_1,x_2,x_3)\mapsto{\bf p}_0(t)+x_1\,{\bf p}_1(t)
+x_2\,{\bf q}_1(t)+x_3\,{\bf q}_2(t).
\end{equation*}
It is easy to show that these SL 3-folds are isomorphic to those
constructed in \cite[Ex.~7.5]{Joyc3}, by evolving the image of an
equivalent quadric in~$\R^3$.

When $k=2$, the SL 3-folds above are equivalent to those of \S\ref{ev8},
under automorphisms of $\C^3$. Putting $k=2$ in \eq{ev12eq10} gives
\begin{equation*}
\Phi(x,y,t)={\bf p}_0(t)+x\,{\bf p}_1(t)+x^2\,{\bf p}_2(t)
+y\,{\bf q}_1(t)+xy\,{\bf q}_2(t),
\end{equation*}
whereas in the construction of \S\ref{ev8} we have
\begin{equation*}
\Phi(y_1,y_2,t)=y_1^2\,{\bf z}_1(t)+y_1y_2\,{\bf z}_3(t)+
y_1\,{\bf z}_4(t)+y_2{\bf z}_5(t)+{\bf z}_6(t),
\end{equation*}
remembering that ${\bf z}_1\equiv{\bf z}_2$. Comparing these two
equations, we see that they agree under the correspondence
\begin{equation*}
x\leftrightarrow y_1,\quad
y\leftrightarrow y_2,\quad
{\bf p}_0\leftrightarrow{\bf z}_6,\quad
{\bf p}_1\leftrightarrow{\bf z}_4,\quad
{\bf p}_2\leftrightarrow{\bf z}_1,\quad
{\bf q}_1\leftrightarrow{\bf z}_5,\quad
{\bf q}_2\leftrightarrow{\bf z}_3.
\end{equation*}
Also, the variable $t$ in \S\ref{ev8} corresponds to $2t$ in this 
section, which is due to the fact that in Example \ref{ev4ex1} we 
constructed $\chi$ from $\frac{\pd}{\pd y_1}\w\frac{\pd}{\pd y_2}$, but 
in Example \ref{ev4ex2} we constructed $\chi$ from~$2\frac{\pd}{\pd x}
\w\frac{\pd}{\pd y}$.

\subsection{Symmetries of the construction}
\label{ev122}

In \S\ref{ev51} we constructed an action of $\GL(2,\R)\lt\R^2$
on the set of solutions ${\bf z}_j$ to \eq{ev5eq10}--\eq{ev5eq12},
which acts trivially on the corresponding SL 3-folds $N$ of 
\eq{ev5eq13}. This group $\GL(2,\R)\lt\R^2$ acts on $\R^2$, and
consists of the original symmetry group $G=\SL(2,\R)\lt\R^2$ 
used to construct the evolution data in Example \ref{ev4ex1},
together with dilations of~$\R^2$.

We shall now do the same thing for the construction above. In
this case, the appropriate group of symmetries of $\R^2$ is
the set of transformations of the form
\e
(x,y)\mapsto(ax+b,cy+d_0+d_1x+\cdots+d_{k-1}x^{k-1}),
\label{ev12eq11}
\e
where $a,b,c$ and $d_0,\ldots,d_{k-1}\in\R$, and $\de=ac$ is nonzero.
The subgroup of transformations with $\de=1$ and $a,c>0$ are the group 
$G\lt U_k$ used in \S\ref{ev4} to construct Example~\ref{ev4ex2}. 

Here $G=\R_+\lt\R$ is the $a,b$ and $c=a^{-1}$ part of the action, and
$U_k\cong\R^k$ is the vector space of 1-forms $p(x)\d y$ for $p(x)$ a
polynomial of degree less than $k$, so that $p(x)=d_0+d_1x+\cdots+
d_{k-1}x^{k-1}$. By allowing $\de$ to be nonzero rather than fixing it
to be 1 we include dilations of $\R^2$ in the group. 

Following the proof of Proposition \ref{ev5prop}, we may show that
the transformation \eq{ev12eq11} of $\R^2$ corresponds to the following
transformation of the ${\bf p}_i$ and~${\bf q}_j$.

\begin{prop} Suppose ${\bf p}_0,\ldots,{\bf p}_k,{\bf q}_1,{\bf q}_2:
\R\ra\C^3$ satisfy \eq{ev12eq6}--\eq{ev12eq8}. Let\/ $a,b,c$ and\/
$d_0,\ldots,d_{k-1}\in\R$ with\/ $\de=ac\ne 0$, and define
${\bf p}_0',\ldots,{\bf p}_k',{\bf q}_1',{\bf q}_2':\R\ra\C^3$ by
\begin{align*}
{\bf p}_0'(t)&=\sum_{i=0}^kb^i\,{\bf p}_i(\de t)+
d_0\,{\bf q}_1(\de t)+bd_0\,{\bf q}_2(\de t),\\
{\bf p}_j'(t)&=\sum_{i=j}^k\binom{i}{j}a^jb^{i-j}\,{\bf p}_i(\de t)
\!+\!d_j\,{\bf q}_1(\de t)\!+\!(ad_{j-1}\!+\!bd_j)\,{\bf q}_2(\de t),
\;\> 1\le j<k,\\
{\bf p}_k'(t)&=a^k\,{\bf p}_k(\de t)+ad_{k-1}\,{\bf q}_2(\de t),\\
{\bf q}_1'(t)&=c\,{\bf q}_1(\de t)+bc\,{\bf q}_2(\de t)
\qquad\text{and}\qquad {\bf q}_2'(t)=ac\,{\bf q}_2(\de t).
\end{align*}
Then ${\bf p}_0',\ldots,{\bf p}_k',{\bf q}_1',{\bf q}_2'$ satisfy
\eq{ev12eq6}--\eq{ev12eq8}. Furthermore, the ${\bf p}_i',{\bf q}_j'$
satisfy \eq{ev12eq4} and\/ \eq{ev12eq5} if and only if the ${\bf p}_i,
{\bf q}_j$ do, and then the special Lagrangian $3$-folds $N,N'$ 
constructed in \eq{ev12eq9} from the ${\bf p}_j,{\bf q}_j$ and\/ 
${\bf p}_i',{\bf q}_j'$ are the same.
\label{ev12prop1}
\end{prop}

We can now give a parameter count for the construction, following the
method of \S\ref{ev52}. The initial data ${\bf p}_0(0),\ldots,{\bf p}_k(0),
{\bf q}_1(0),{\bf q}_2(0)$ has $6(k+3)$ parameters. These must satisfy 
\eq{ev12eq4} and \eq{ev12eq5}, which is $k+2$ equations. Thus the set 
${\cal C}_P$ of Definition \ref{ev3def} has dimension~$5k+16$.

From this we must subtract symmetries of three kinds. Firstly, the
internal symmetry group of Proposition \ref{ev12prop1} has dimension
$k+3$. Secondly, the automorphisms $\SU(3)\lt\C^3$ of $\C^3$ have
dimension 14. Thirdly, we subtract one dimension to allow for
translation in time, $t\mapsto t+c$. Taking all these into account,
we find that the family of distinct special Lagrangian 3-folds in
$\C^3$ constructed in Theorem \ref{ev12thm1}, up to automorphisms
of $\C^3$, has dimension~$4k-2$.

\subsection{Solving the equations}
\label{ev123}

We shall now solve equations \eq{ev12eq6}--\eq{ev12eq8} for the ${\bf p}_i$ 
and ${\bf q}_j$ fairly explicitly, under the restrictions \eq{ev12eq4} 
and \eq{ev12eq5}, and making use of the symmetries discussed above.
Our treatment follows \S\ref{ev8}, and indeed the case $k=2$ is
equivalent to the construction of \S\ref{ev8}. 

We begin by putting ${\bf q}_1,{\bf q}_2$ in a convenient form. 
Divide into the two cases
\begin{itemize}
\item[(a)] ${\bf q}_1(0)$ and ${\bf q}_2(0)$ are linearly dependent, and
\item[(b)] ${\bf q}_1(0)$ and ${\bf q}_2(0)$ are linearly independent.
\end{itemize}
It is easy to show from \eq{ev12eq8} that in case (a) ${\bf q}_1,{\bf q}_2$ 
are constant, and the SL $3$-fold $N$ of \eq{ev12eq9} splits as a product 
$\Si\t\R$ in $\C^2\t\C$, where $\Si$ is an SL $2$-fold in $\C^2$. Thus
case (a) is not very interesting, and we shall not consider it further.

In case (b) we may follow the proof of Theorem \ref{ev8thm1} to prove:

\begin{prop} Let\/ ${\bf p}_0,\ldots,{\bf p}_k$ and\/ 
${\bf q}_1,{\bf q}_2$ satisfy equations \eq{ev12eq4}--\eq{ev12eq8}, and 
suppose ${\bf q}_1(0),{\bf q}_2(0)$ are linearly independent. Then we 
may transform the ${\bf p}_i,{\bf q}_j$ under $\SU(3)$ and the symmetries 
of Proposition \ref{ev12prop1} to ${\bf p}_i',{\bf q}_j'$, where
\begin{equation*}
{\bf q}_1'(t)=({\rm e}^{it},i{\rm e}^{-it},0) \quad\text{and}\quad
{\bf q}_2'(t)=(0,0,1).
\end{equation*}
\label{ev12prop2}
\end{prop}

So let us fix ${\bf q}_1(t)=({\rm e}^{it},i{\rm e}^{-it},0)$ and 
${\bf q}_2(t)=(0,0,1)$ for the rest of the section. Next we rewrite 
equations \eq{ev12eq4}--\eq{ev12eq7} for ${\bf p}_0,\ldots,{\bf p}_k$.
The proof follows immediately from \eq{ev5eq9} and the definition
of~$\om$.

\begin{prop} Set\/ ${\bf q}_1(t)=({\rm e}^{it},i{\rm e}^{-it},0)$ 
and\/ ${\bf q}_2(t)=(0,0,1)$, and write ${\bf p}_j=(a_j,b_j,c_j)$ 
for $j=0,\ldots,k$, where $a_j,b_j,c_j:\R\ra\C$ are differentiable 
functions. Then \eq{ev12eq8} holds, and equations \eq{ev12eq6} and\/
\eq{ev12eq7} are equivalent to
\begin{align}
\frac{\d a_0}{\d t}&=i{\rm e}^{it}\bar c_1,\quad
\frac{\d b_0}{\d t}={\rm e}^{-it}\bar c_1,\quad
\frac{\d c_0}{\d t}=-i{\rm e}^{it}\bar a_1-{\rm e}^{-it}\bar b_1,
\label{ev12eq12}\\
\frac{\d a_j}{\d t}&=(j+1)i{\rm e}^{it}\bar c_{j+1}+j\bar b_j
\qquad\text{for $j=1,\ldots,k-1$,}
\label{ev12eq13}\\
\frac{\d b_j}{\d t}&=(j+1){\rm e}^{-it}\bar c_{j+1}-j\bar a_j
\qquad\text{for $j=1,\ldots,k-1$,}
\label{ev12eq14}\\
\frac{\d c_j}{\d t}&=-(j+1)(i{\rm e}^{it}\bar a_{j+1}+{\rm e}^{-it}
\bar b_{j+1})\quad\text{for $j=1,\ldots,k-1$,}
\label{ev12eq15}\\
\frac{\d a_k}{\d t}&=k\bar b_k,\qquad
\frac{\d b_k}{\d t}=-k\bar a_k\qquad\text{and}\qquad
\frac{\d c_k}{\d t}=0.
\label{ev12eq16}
\end{align}
Furthermore, equations \eq{ev12eq4} and\/ \eq{ev12eq5} are equivalent to
$\Im c_k=0$ and
\e
j\Im({\rm e}^{-it}a_j-i{\rm e}^{it}b_j)+(j-1)\Im c_{j-1}=0
\quad\text{for $j=1,\ldots,k$.}
\label{ev12eq17}
\e
\label{ev12prop3}
\end{prop}

The best way to solve equations \eq{ev12eq12}--\eq{ev12eq16} is in
reverse order. That is, we begin by solving \eq{ev12eq16} for
$a_k,b_k,c_k$. Then we inductively solve equations
\eq{ev12eq13}--\eq{ev12eq15} for $a_j,b_j,c_j$ with $j=k-1,k-2,\ldots,1$,
treating $a_{j+1},b_{j+1}$ and $c_{j+1}$ as known. Note that we can
combine \eq{ev12eq13} and \eq{ev12eq14} to get
\e
\frac{\d^2a_j}{\d t^2}+j^2a_j=j(j+1){\rm e}^{it}c_{j+1}+
(j+1)\frac{\d}{\d t}\bigl(i{\rm e}^{it}\bar c_{j+1}\bigr),
\label{ev12eq18}
\e
which is a linear second-order o.d.e.\ for $a_j$ with prescribed
right hand side, and can be solved by standard techniques.
Finally we solve \eq{ev12eq12} for $a_0,b_0,c_0$, which is just a
matter of integration. Here are the first four steps in this process.

The general solutions of \eq{ev12eq16} are easily shown to be
\e
\begin{gathered}
a_k(t)=A_k{\rm e}^{ikt}+B_k{\rm e}^{-ikt},\qquad
b_k(t)=i\bar B_k{\rm e}^{ikt}-i\bar A_k{\rm e}^{-ikt}\\
\text{and}\qquad c_k(t)=C_k,
\end{gathered}
\label{ev12eq19}
\e
for constants $A_k,B_k,C_k\in\C$. The condition $\Im c_k=0$ gives
$\Im C_k=0$. But applying Proposition \ref{ev12prop1} with $a=c=1$,
$b=d_0=\cdots=d_{k-2}=0$ and $d_{k-1}=-\Re C_k$ shows that we can
use the symmetries of the construction to set $\Re C_k=0$ as well.
So let us set~$C_k=0$.

With these values for $a_k,b_k,c_k$ we can solve equations
\eq{ev12eq13}--\eq{ev12eq15} for $a_{k-1}$, $b_{k-1}$ and
$c_{k-1}$. We get
\begin{align}
a_{k-1}(t)&=A_{k-1}{\rm e}^{i(k-1)t}+B_{k-1}{\rm e}^{-i(k-1)t},
\label{ev12eq20}\\
b_{k-1}(t)&=i\bar B_{k-1}{\rm e}^{i(k-1)t}-i\bar A_{k-1}{\rm e}^{-i(k-1)t},
\label{ev12eq21}\\
\begin{split}
\text{and}\quad
c_{k-1}(t)&=-{\ts\frac{k}{k-1}}A_k{\rm e}^{i(k-1)t}
+{\ts\frac{k}{k-1}}\bar A_k{\rm e}^{-i(k-1)t}\\
&\quad -{\ts\frac{k}{k+1}}B_k{\rm e}^{-i(k+1)t}
-{\ts\frac{k}{k+1}}\bar B_k{\rm e}^{i(k+1)t}+C_{k-1},
\end{split}
\label{ev12eq22}
\end{align}
for $A_{k-1},B_{k-1},C_{k-1}\in\C$, provided $k\ne 1$. When $k=1$
the terms in $A_k,\bar A_k$ in \eq{ev12eq22} are replaced by
$-2i\Re(A_1)t$. The case $j=k$ of \eq{ev12eq17} reduces to
$\Im C_{k-1}=0$, and as above we can use a symmetry involving
the variable $d_{k-2}$ in Proposition \ref{ev12prop1} to set
$\Re C_{k-1}=0$ too. So fix~$C_{k-1}=0$.

Repeating the same process for $a_{k-2},b_{k-2}$ and $c_{k-2}$, we obtain
\begin{align}
\begin{split}
a_{k-2}(t)&=A_{k-2}{\rm e}^{i(k-2)t}+B_{k-2}{\rm e}^{-i(k-2)t}\\
&\quad +{\ts\frac{k}{2}}A_k{\rm e}^{ikt}
+{\ts\frac{k(k-1)}{2(k+1)}}B_k{\rm e}^{-ikt}
-{\ts\frac{(k-1)}{2(k+1)}}\bar B_k{\rm e}^{i(k+2)t},
\end{split}
\label{ev12eq23}\\
\begin{split}
b_{k-2}(t)&=i\bigl(\bar B_{k-2}-{\ts\frac{k}{k-2}}A_k\bigr)
{\rm e}^{i(k-2)t}-i\bar A_{k-2}{\rm e}^{-i(k-2)t}\\
&\quad -{\ts\frac{i(k-1)}{2(k+1)}}B_k{\rm e}^{-i(k+2)t}
+{\ts\frac{ik(k-1)}{2(k+1)}}\bar B_k{\rm e}^{ikt}
-{\ts\frac{ik}{2}}\bar A_k{\rm e}^{-ikt},
\end{split}
\label{ev12eq24}\\
\begin{split}
c_{k-2}(t)&=-{\ts\frac{k-1}{k-2}}A_{k-1}{\rm e}^{i(k-2)t}
+{\ts\frac{k-1}{k-2}}\bar A_{k-1}{\rm e}^{-i(k-2)t}\\
&\quad -{\ts\frac{k-1}{k}}B_{k-1}{\rm e}^{-ikt}
-{\ts\frac{k-1}{k}}\bar B_{k-1}{\rm e}^{ikt}+C_{k-2},
\end{split}
\label{ev12eq25}
\end{align}
for $A_{k-2},B_{k-2},C_{k-2}\in\C$, provided $k\ne 2$. When $k=2$
the terms in $A_{k-1},\bar A_{k-1}$ in \eq{ev12eq25} should be
$-2i\Re(A_1)t$. Putting $j=k-1$ in \eq{ev12eq17} gives $\Im C_{k-2}=0$, and
we can again use symmetry to set $\Re C_{k-2}=0$, so we fix~$C_{k-2}=0$.

There is no term in $\bar A_k$ in \eq{ev12eq23}, because the $\bar A_k$
terms on the right hand side of \eq{ev12eq18} with $j=k-2$ cancel out.
If this had not happened, then \eq{ev12eq23} and \eq{ev12eq24} would have
included multiples of $t{\rm e}^{\pm i(k-2)t}$. This will be significant
later, when we consider periodicity of the functions~$a_j,b_j,c_j$.

Applying the same process for $a_{k-3},b_{k-3}$ and $c_{k-3}$, we obtain
\begin{align}
\begin{split}
a_{k-3}(t)&=A_{k-3}{\rm e}^{i(k-3)t}+B_{k-3}{\rm e}^{-i(k-3)t}
+{\ts\frac{k-1}{2}}A_{k-1}{\rm e}^{i(k-1)t}\\
&\quad
+{\ts\frac{(k-1)(k-2)}{2k}}B_{k-1}{\rm e}^{-i(k-1)t}
-{\ts\frac{(k-2)}{2k}}\bar B_{k-1}{\rm e}^{i(k+1)t},
\end{split}
\label{ev12eq26}\\
\begin{split}
b_{k-3}(t)&=i\bigl(\bar B_{k-3}-{\ts\frac{k-1}{k-3}}A_{k-1}\bigr)
{\rm e}^{i(k-3)t}-i\bar A_{k-3}{\rm e}^{-i(k-3)t}\\
-{\ts\frac{i(k\!-\!2)}{2k}}&B_{k\!-\!1}{\rm e}^{-i(k\!+\!1)t}
\!+\!{\ts\frac{i(k\!-\!1)(k\!-\!2)}{2k}}\bar B_{k\!-\!1}{\rm e}^{i(k\!-\!1)t}
\!-\!{\ts\frac{i(k\!-\!1)}{2}}\bar A_{k\!-\!1}{\rm e}^{-i(k\!-\!1)t},
\end{split}
\label{ev12eq27}\\
\begin{split}
c_{k-3}(t)&=-{\ts\frac{k-2}{k-3}}A_{k-2}{\rm e}^{i(k-3)t}
+{\ts\frac{k-2}{k-3}}\bar A_{k-2}{\rm e}^{-i(k-3)t}\\
-{\ts\frac{k-2}{k-1}}&
\bigl(B_{k-2}\!-\!{\ts\frac{k^2}{2(k-2)}}\bar A_k\bigr){\rm e}^{-i(k-1)t}
\!-\!{\ts\frac{k-2}{k-1}}
\bigl(\bar B_{k-2}\!+\!{\ts\frac{k}{2}}A_k\bigr){\rm e}^{i(k-1)t}\\
-&{\ts\frac{(k\!-\!1)(k\!-\!2)}{2(k+1)}}B_k{\rm e}^{-i(k\!+\!1)t}
-{\ts\frac{(k\!-\!1)(k\!-\!2)}{2(k+1)}}\bar B_k{\rm e}^{i(k\!+\!1)t}
\!+\!C_{k-3},
\end{split}
\label{ev12eq28}
\end{align}
for $A_{k-3},B_{k-3},C_{k-3}\in\C$, provided $k\ne 3$. When $k=3$
the terms in $A_{k-2},\bar A_{k-2}$ in \eq{ev12eq28} should be
$-2i\Re(A_1)t$. As above we may fix~$C_{k-3}=0$. 

The reader may readily carry on calculating $a_j,b_j,c_j$ for 
$j=k-4,k-5,\ldots$ by this method. The expressions get steadily 
longer and more complicated, so we shall stop at this point. 
In Theorem \ref{ev12thm4} we will give the general form of the 
solutions $a_j,b_j,c_j$ for all $j$, without the constant factors
depending on $j$ and~$k$.

As an example, let $k=4$. Then equations \eq{ev12eq19}--\eq{ev12eq28} 
give $a_j,b_j,c_j$ for $j=1,2,3,4$ in terms of complex constants
$A_1,\ldots,A_4$ and $B_1,\ldots,B_4$, and we get $a_0,b_0,c_0$ by 
integrating \eq{ev12eq12}. Equation \eq{ev12eq17} holds for $j=2,3,4$
by construction. However, we still have to satisfy \eq{ev12eq17}
for $j=1$. This is $\Im({\rm e}^{-it}a_1-i{\rm e}^{it}b_1)=0$, 
which simplifies to $2\Im A_1=0$ after substituting in for $a_1,b_1$. 
So $A_1\in\R$, and from Theorem \ref{ev12thm1} we deduce:

\begin{thm} Define functions ${\bf p}_0,\ldots,{\bf p}_4,
{\bf q}_1,{\bf q}_2:\R\ra\C^3$ by 
\begin{align}
\begin{split}
{\bf p}_0(t)=\bigl(&
-2i\bar A_2t+A_2{\rm e}^{2it}-{\ts\frac{1}{6}}(\bar B_2-4A_4){\rm e}^{4it}\\
&\quad
+{\ts\frac{1}{3}}(B_2+2\bar A_4){\rm e}^{-2it}
-{\ts\frac{1}{10}}\bar B_4{\rm e}^{6it}+{\ts\frac{3}{20}}B_4{\rm e}^{-4it},\\
&2A_2t-i\bar A_2{\rm e}^{-2it}+{\ts\frac{i}{3}}(\bar B_2-4A_4){\rm e}^{2it}\\
&\quad
-{\ts\frac{i}{6}}(B_2+2\bar A_4){\rm e}^{-4it}
-{\ts\frac{i}{10}}B_4{\rm e}^{-6it}+{\ts\frac{3i}{20}}\bar B_4{\rm e}^{4it},\\
&-2iA_1t-{\ts\frac{1}{2}}(B_1-{\ts\frac{9}{2}}\bar A_3){\rm e}^{-2it}\\
&\quad
-{\ts\frac{1}{2}}(\bar B_1+{\ts\frac{3}{2}}A_3){\rm e}^{2it}
-{\ts\frac{1}{4}}B_3{\rm e}^{-4it}-{\ts\frac{1}{4}}\bar B_3{\rm e}^{4it}\bigr),
\end{split}
\label{ev12eq29}
\displaybreak[0]\\
\begin{split}
{\bf p}_1(t)=\bigl(&
A_1{\rm e}^{it}+B_1{\rm e}^{-it}+{\ts\frac{3}{2}}A_3{\rm e}^{3it}
+{\ts\frac{3}{4}}B_3{\rm e}^{-3it}-{\ts\frac{1}{4}}\bar B_3{\rm e}^{5it},\\
&i(\bar B_1\!-\!3A_3){\rm e}^{it}\!-\!iA_1{\rm e}^{-it}
\!-\!{\ts\frac{i}{4}}B_3{\rm e}^{-5it}
\!+\!{\ts\frac{3i}{4}}\bar B_3{\rm e}^{3it}
\!-\!{\ts\frac{3i}{2}}\bar A_3{\rm e}^{-3it},\\
&-2A_2{\rm e}^{it}+2\bar A_2{\rm e}^{-it}
-{\ts\frac{2}{3}}(B_2-4\bar A_4){\rm e}^{-3it}\\
&\quad
-{\ts\frac{2}{3}}(\bar B_2+2A_4){\rm e}^{3it}-{\ts\frac{3}{5}}B_4{\rm e}^{-5it}
-{\ts\frac{3}{5}}\bar B_4{\rm e}^{5it}\bigr),
\end{split}
\label{ev12eq30}
\displaybreak[0]\\
\begin{split}
{\bf p}_2(t)=\bigl(&
A_2{\rm e}^{2it}+B_2{\rm e}^{-2it}+2A_4{\rm e}^{4it}
+{\ts\frac{6}{5}}B_4{\rm e}^{-4it}-{\ts\frac{3}{10}}\bar B_4{\rm e}^{6it},\\
&i(\bar B_2\!-\!2A_4){\rm e}^{2it}\!-\!i\bar A_2{\rm e}^{-2it}
\!-\!{\ts\frac{3i}{10}}B_4{\rm e}^{-6it}
\!+\!{\ts\frac{6i}{5}}\bar B_4{\rm e}^{4it}
\!-\!2i\bar A_4{\rm e}^{-4it},\\
&-{\ts\frac{3}{2}}A_3{\rm e}^{2it}
+{\ts\frac{3}{2}}\bar A_3{\rm e}^{-2it}
-{\ts\frac{3}{4}}B_3{\rm e}^{-4it}
-{\ts\frac{3}{4}}\bar B_3{\rm e}^{4it}\bigr),
\end{split}
\label{ev12eq31}\\
\begin{split}
{\bf p}_3(t)=\bigl(&A_3{\rm e}^{3it}+B_3{\rm e}^{-3it},
i\bar B_3{\rm e}^{3it}-i\bar A_3{\rm e}^{-3it},\\
&-{\ts\frac{4}{3}}A_4{\rm e}^{3it}
+{\ts\frac{4}{3}}\bar A_4{\rm e}^{-3it}
-{\ts\frac{4}{5}}B_4{\rm e}^{-5it}
-{\ts\frac{4}{5}}\bar B_4{\rm e}^{5it}\bigr),
\end{split}
\label{ev12eq32}
\displaybreak[0]\\
{\bf p}_4(t)=\bigl(&A_4{\rm e}^{4it}+B_4{\rm e}^{-4it},
i\bar B_4{\rm e}^{4it}-i\bar A_4{\rm e}^{-4it},0\bigr),
\label{ev12eq33}\\
{\bf q}_1(t)=(&{\rm e}^{it},i{\rm e}^{-it},0) \qquad\text{and}\qquad
{\bf q}_2(t)=(0,0,1),
\end{align}
where $A_1\in\R$ and\/ $A_2,A_3,A_4,B_1,\ldots,B_4\in\C$. Define a 
subset\/ $N$ of\/ $\C^3$ to be
\e
\bigl\{{\bf p}_0(t)\!+\!x\,{\bf p}_1(t)\!+\!\cdots\!+\!x^4\,{\bf p}_4(t)
\!+\!y\,{\bf q}_1(t)\!+\!xy\,{\bf q}_2(t):x,y,t\in\R\bigr\}.
\label{ev12eq34}
\e
Then $N$ is a special Lagrangian $3$-fold in $\C^3$ wherever
it is nonsingular.
\label{ev12thm2}
\end{thm}

The expression \eq{ev12eq29} for ${\bf p}_0(t)$ could also include 
three complex constants of integration, but we have set them to 
zero for simplicity. For general $k\ge 1$, one can prove the
following result.

\begin{thm} In the situation above, for each\/ $k\ge 1$ there exist
solutions $a_j,b_j,c_j$ to equations \eq{ev12eq12}--\eq{ev12eq17}
depending on $A_1\in\R$, $A_2,\ldots,A_k\in\C$ and\/ $B_1,\ldots,B_k
\in\C$, such that
\begin{itemize} 
\item[{\rm(i)}] For $1\le j\le k$, $a_j$ is a real linear combination 
of terms $A_{j+2l}{\rm e}^{i(j+2l)t}$, $B_{j+2l}{\rm e}^{-i(j+2l)t}$,
$\bar A_{j+2l+4}{\rm e}^{-i(j+2l+2)t}$ 
and\/~$\bar B_{j+2l+2}{\rm e}^{i(j+2l+4)t}$;
\item[{\rm(ii)}] For $1\le j\le k$, $b_j$ is a real linear combination of 
terms $iA_{j+2l+2}{\rm e}^{i(j+2l)t}$, $iB_{j+2l+2}{\rm e}^{-i(j+2l+4)t}$,
$i\bar A_{j+2l}{\rm e}^{-i(j+2l)t}$ and\/~$i\bar B_{j+2l}{\rm e}^{i(j+2l)t}$;
\item[{\rm(iii)}] For $1\le j\le k$, $c_j$ is a real linear combination of 
terms $A_{j+2l+1}{\rm e}^{i(j+2l)t}$, $B_{j+2l+1}{\rm e}^{-i(j+2l+2)t}$,
$\bar A_{j+2l+1}{\rm e}^{-i(j+2l)t}$
and\/~$\bar B_{j+2l+1}{\rm e}^{i(j+2l+2)t}$;
\item[{\rm(iv)}] $a_0$ is a real linear combination of terms $i\bar A_2t$, 
$A_{2l\!+\!2}{\rm e}^{i(2l\!+\!2)t}$, $B_{2l\!+\!2}{\rm e}^{-i(2l\!+\!2)t}$, 
$\bar A_{2l+4}{\rm e}^{-i(2l+2)t}$ and\/~$\bar B_{2l+2}{\rm e}^{i(2l+4)t}$;
\item[{\rm(v)}] $b_0$ is a real linear combination of terms $A_2t$, 
$iA_{2l\!+\!4}{\rm e}^{i(2l\!+\!2)t}$, $iB_{2l\!+\!2}{\rm e}^{-i(2l\!+\!4)t}$,
$i\bar A_{2l+2}{\rm e}^{-i(2l+2)t}$ and\/~$i\bar B_{2l+2}{\rm e}^{i(2l+2)t}$;
and
\item[{\rm(vi)}] $c_0$ is a real linear combination of terms $iA_1t$, 
$A_{2l\!+\!3}{\rm e}^{i(2l\!+\!2)t}$, $B_{2l\!+\!1}{\rm e}^{-i(2l\!+\!2)t}$,
$\bar A_{2l+3}{\rm e}^{-i(2l+2)t}$ and\/~$\bar B_{2l+1}{\rm e}^{i(2l+2)t}$.
\end{itemize}
Here in each case $l\ge 0$ is an integer, and we take $A_j=B_j=0$ for~$j>k$. 
\label{ev12thm3}
\end{thm}

One surprising thing about this theorem is that the solutions contain 
no terms like $t^a{\rm e}^{ibt}$ for $a>0$ and $b$ integers, except
for the $t$ terms in $a_0,b_0$ and $c_0$. Here is the reason why. Let 
$1\le j<k$, and suppose by induction that we know $a_i,b_i,c_i$ for 
$i=j+1,\ldots,k$, and that parts (i)--(iii) hold for them. To find $a_j$ 
we must solve equation \eq{ev12eq18}. By part (iii) above, the right 
hand side of \eq{ev12eq18} is a linear combination of exponentials 
${\rm e}^{int}$ for various integers~$n$. 

The corresponding terms in $a_j$ are multiples of ${\rm e}^{int}$ if 
$n\ne\pm j$, but multiples of $t{\rm e}^{int}$ if $n=\pm j$. Thus for 
$a_j$ to be of the form given in part (i), we need the right hand side 
of \eq{ev12eq18} to contain no multiples of ${\rm e}^{\pm ijt}$. By part
(iii), we expect to get a multiple of $\bar A_{j+2}{\rm e}^{-ijt}$.
However, this multiple is zero, as we saw in \eq{ev12eq23} above,
and so $a_j$ satisfies (i). It easily follows that $b_j$ and $c_j$ 
satisfy (ii) and (iii), and the inductive step is complete.

\subsection{Periodicity}
\label{ev124}

In Theorems \ref{ev12thm2} and \ref{ev12thm3}, the only terms which are 
not periodic in $t$ with common period $2\pi$ are those in $\bar A_2t$, 
$A_2t$ and $A_1t$ in equation \eq{ev12eq29} and parts (iv)--(vi) of
Theorem \ref{ev12thm3}. So setting $A_1=A_2=0$ gives $\Phi(x,y,t+2\pi)
=\Phi(x,y,t)$ for all $x,y,t$, where $\Phi$ is defined in~\eq{ev12eq10}. 

Furthermore, Theorem \ref{ev12thm3} shows that ${\bf p}_j(t+\pi)=
(-1)^j{\bf p}_j(t)$, and clearly ${\bf q}_j(t+\pi)=(-1)^j{\bf q}_j(t)$. 
This gives $\Phi(x,y,t+\pi)=\Phi(-x,-y,t)$ for all $x,y,t\in\R$.
We may therefore regard $\Phi$ as mapping $\R^3/\Z\ra\C^3$, where 
$\Z$ acts on $\R^3$ by $(x,y,t)\,{\buildrel n\over\longmapsto}\,
\bigl((-1)^nx,(-1)^ny,y+n\pi\bigr)$, for~$n\in\Z$. 

For $k\ge 3$ one can show that $\Phi$ is an immersion for
generic values of $A_3,\ldots,B_k$. Then $N$ is a nonsingular 
immersed 3-submanifold diffeomorphic to $\R^3/\Z$, or equivalently 
to ${\cal S}^1\t\R^2$. We have proved:

\begin{thm} For each\/ $k\ge 3$, the construction above with\/
$A_1\!=\!A_2\!=\!0$ gives a family of special Lagrangian $3$-folds in $\C^3$ 
depending on $A_3,\ldots,A_k$ and\/ $B_1,\ldots,B_k\in\C$, with 
generic member a closed, nonsingular, immersed\/ $3$-submanifold 
diffeomorphic to~${\cal S}^1\t\R^2$.
\label{ev12thm4}
\end{thm}

Here is a parameter count for this family. It depends upon 
$A_3,\ldots,A_k$ and $B_1,\ldots,B_k\in\C$, which is $4k-4$ real 
parameters. It can be shown that the only symmetry left to take into 
account is translation in time, $t\mapsto t+c$. Subtracting one 
for this, the family of distinct SL 3-folds in $\C^3$ in Theorem 
\ref{ev12thm4}, up to automorphisms of $\C^3$, has dimension~$4k-5$.

Next we discuss the {\it asymptotic behaviour} of the SL 3-folds $N$ of
Theorem \ref{ev12thm4} near infinity in $\C^3$. Suppose for simplicity 
that $A_k$ and $B_k$ are not both zero, since if they are we can reduce 
$k$ to $k-1$. Then ${\bf p}_k$ is nonzero for all $t$, by \eq{ev12eq19}.
It is clear from \eq{ev12eq10} that $\Phi(x,y,t)\ra\iy$ in $\C^3$ 
as $(x,y)\ra\iy$ in $\R^2$. To describe the asymptotic behaviour
at infinity to leading order, we need to decide which terms in 
\eq{ev12eq10} are dominant when $x,y$ are large, and neglect the
other terms.

Obviously, when $\md{x}$ is large the term $x^k\,{\bf p}_k(t)$ 
dominates the terms $x^j\,{\bf p}_j(t)$ for $j<k$, as $p_k(t)$ 
is always nonzero. Thus the three terms in \eq{ev12eq10} which 
may dominate are $x^k\,{\bf p}_k(t)$, $y\,{\bf q}_1(t)$ and 
$xy\,{\bf q}_2(t)$. Careful thought shows that for fixed $t$ 
there are the following four different asymptotic regimes, 
depending on the relative sizes of $x$ and~$y$: 
\begin{itemize}
\item[(i)] $x\gg 0$, $x=O(r^{1/k})$, $y=O(r^{(k-1)/k})$, 
with~$\Phi(x,y,t)\approx x^k\,{\bf p}_k(t)+xy\,{\bf q}_2(t)$.
\item[(ii)] $y\gg 0$, $x=O(1)$, $y=O(r)$, 
with~$\Phi(x,y,t)\approx y\,{\bf q}_1(t)+xy\,{\bf q}_2(t)$.
\item[(iii)] $x\ll 0$, $x=O(r^{1/k})$, $y=O(r^{(k-1)/k})$, 
with~$\Phi(x,y,t)\approx x^k\,{\bf p}_k(t)+xy\,{\bf q}_2(t)$.
\item[(iv)] $y\ll 0$, $x=O(1)$, $y=O(r)$, 
with~$\Phi(x,y,t)\approx y\,{\bf q}_1(t)+xy\,{\bf q}_2(t)$.
\end{itemize}
Here $r$ is the radius function on~$\C^3$.

These four regimes join onto each other in a cyclic fashion,
with the single dominant term $xy\,{\bf q}_2(t)$ at the junction,
so that for instance the junction between (i) and (ii) we have 
$x\gg 0$, $y\gg 0$ and $\Phi(x,y,t)\approx xy\,{\bf q}_2(t)$. 
When $t$ varies as well we can identify (i) with (iii), and (ii) 
with (iv), since $\Phi(x,y,t+\pi)=\Phi(-x,-y,t)$. So there are only 
really two kinds of asymptotic behaviour to consider.

Using the expressions above for ${\bf p}_k$ and ${\bf q}_2$,
we find that in case (i) we have
\begin{equation*}
\Phi(x,y,t)\approx
\bigl(x^k(A_k{\rm e}^{ikt}+B_k{\rm e}^{-ikt}),
x^k(i\bar B_k{\rm e}^{ikt}-i\bar A_k{\rm e}^{-ikt}),xy\bigr),
\end{equation*}
with $x>0$. This sweeps out the special Lagrangian 3-plane
\begin{equation*}
L_1=\ban{(A_k+B_k,i\bar B_k-i\bar A_k,0),
(-iA_k+iB_k,\bar B_k+\bar A_k,0),(0,0,1)}_{\sst\mathbb R}
\end{equation*}
in $\C^3$. Considering $t$ to be a cyclic coordinate with period
$2\pi$, we see that to leading order $\Phi$ is a $k$-{\it fold 
branched cover} of $L_1$, branched along the real 
line~$\ban{(0,0,1)}_{\sst\mathbb R}$.

Similarly, in case (ii) we have
\begin{equation*}
\Phi(x,y,t)\approx (y{\rm e}^{it},iy{\rm e}^{-it},xy),
\end{equation*}
with $y>0$. This sweeps out the special Lagrangian 3-plane
\begin{equation*}
L_2=\ban{(1,i,0),(i,1,0),(0,0,1)}_{\sst\mathbb R}
\end{equation*}
in $\C^3$. As $t$ is cyclic with period $2\pi$ this is a 1-1 
correspondence, rather than a branched cover. If $A_k=0$ then 
$L_1=L_2$. We shall assume $A_k\ne 0$ for simplicity, though 
it doesn't make much difference.

Thus we arrive at the following description of the SL 3-fold $N$ 
at infinity in $\C^3$. It is weakly asymptotic to the union of
two special Lagrangian 3-planes $L_1,L_2$ in $\C^3$, which 
intersect in the real line $\ban{(0,0,1)}_{\sst\mathbb R}$.
Along $L_1$ it converges to a $k$-fold branched cover, branched 
along $\ban{(0,0,1)}_{\sst\mathbb R}$, so that $k$ points of
$N$ `converge' to 1 point of $L_1$ at infinity. Along $L_2$
the convergence is~1-1.

The boundary of $N\cong{\cal S}^1\t\R^2$ at infinity is $T^2$,
whereas the boundary of $L_1\cup L_2$ at infinity is two
copies of ${\cal S}^2$, intersecting in two points. The $T^2$
wraps itself round these two ${\cal S}^2$, so as to cover the
first $k$ times and the second once. The order of convergence 
is $O(r^{(k-1)/k})$, which is rather weak, as it means that
towards infinity $N$ gets further away from $L_1\cup L_2$ 
rather than closer.

\end{document}